\documentclass{scrartcl}
\usepackage{amsmath}
\usepackage{amsfonts}
\usepackage{amssymb}
\usepackage{dsfont}
\usepackage[T1]{fontenc}
\usepackage[latin1]{inputenc}
\usepackage{epsfig}
\usepackage{graphicx}

\usepackage[boxed, lined]{algorithm2e}
\usepackage{float}
\usepackage{comment}
\usepackage{mathtools}

\usepackage{bm}

\newcommand{\norm}[1]{\left\lVert#1\right\rVert}

\newcommand{\subalign}[1]{%
   \vcenter{%
     \Let@ \restore@math@cr \default@tag
     \baselineskip\fontdimen10 \scriptfont\tw@
     \advance\baselineskip\fontdimen12 \scriptfont\tw@
     \lineskip\thr@@\fontdimen8 \scriptfont\thr@@
     \lineskiplimit\lineskip
     \ialign{\hfil$\m@th\scriptstyle##$&$\m@th\scriptstyle{}##$\hfil\crcr
       #1\crcr
     }%
   }%
}

\newcommand{\bll}{\mathbf{l}}

\newcommand{\bk}{\mathbf{k}}
\newcommand{\bj}{\mathbf{j}}

\newcommand{\bx}{x}
\newcommand{\bX}{X}

\newcommand{\bz}{\mathbf{z}}
\newcommand{\bw}{\mathbf{w}}
\newcommand{\bv}{\mathbf{v}}
\newcommand{\bL}{\mathbf{L}}
\newcommand{\bM}{\mathbf{M}}

\newcommand{\btheta}{\mathbf{\vartheta}}
\newcommand{\bTheta}{\mathbf{\Theta}}

\newcommand{\D}{{\mathcal{D}}}

\newcommand{\W}{{\mathcal{W}}}

\newcommand{\Mu}{{\mathcal{M}}}

\newcommand{\N}{\mathbb{N}}

\newcommand{\R}{\mathbb{R}}
\newcommand{\Z}{\mathbb{Z}}
\newcommand{\Rd}{\mathbb{R}^d}
\newcommand{\IND}{\mathbb{I}}
\newcommand{\balpha}{\mathbf{\alpha}}

\newcommand{\beq}{\begin{eqnarray*}}

\newcommand{\eeq}{\end{eqnarray*}}

\newcommand{\beqm}{\begin{eqnarray}}

\newcommand{\eeqm}{\end{eqnarray}}

\newtheorem{theorem}{Theorem}

\newtheorem{lemma}{Lemma}
\newtheorem{definition}{Definition}

\DeclareOldFontCommand{\bf}{\normalfont\bfseries}{\mathbf}
\DeclareOldFontCommand{\it}{\normalfont\itshape}{\mathit}

\newcommand{\EXP}{{\mathbf E}}
\newcommand{\PROB}{{\mathbf P}}

\renewcommand{\P}{{\cal P}}

\allowdisplaybreaks
\begin{document}
\renewcommand{\thefootnote}{\fnsymbol{footnote}}
\newcommand{\F}{{\cal F}}
\newcommand{\Sp}{{\cal S}}
\newcommand{\G}{{\cal G}}
\newcommand{\HH}{{\cal H}}

\begin{center}

  {\LARGE \bf
    Learning of deep convolutional network image classifiers via 
stochastic gradient descent and over-parametrization
  }
\footnote{
Running title: {\it Deep network classifiers}}
\vspace{0.5cm}

Michael Kohler $^{1}$,
Adam Krzy\.zak$^{2}$
and Alisha S\"anger $^{1,}$\footnote{Corresponding author. Tel:+49 6151 16-23386}

{\it $^1$
Fachbereich Mathematik, Technische Universit\"at Darmstadt,
Schlossgartenstr. 7, 64289 Darmstadt, Germany,
email: kohler@mathematik.tu-darmstadt.de, saenger@mathematik.tu-darmstadt.de
}

{\it $^2$ Department of Computer Science and Software Engineering, Concordia University, 1455 De Maisonneuve Blvd. West, Montreal, Quebec, Canada H3G 1M8, email: krzyzak@cs.concordia.ca}

\end{center}
\vspace{0.5cm}

\begin{center}
December 12, 2024
\end{center}
\vspace{0.5cm}

\noindent
    {\bf Abstract}\\
Image classification from independent and identically
distributed random variables is considered. 
Image classifiers are 
defined which are 
based on a linear combination of deep convolutional networks with max-pooling layer. Here
all the weights are learned by stochastic gradient descent. A general
result is presented which shows that the image classifiers are able
to approximate the best possible deep convolutional network. In case that
the a posteriori probability satisfies a suitable hierarchical composition
model it is shown that the corresponding deep convolutional neural
network image classifier achieves a rate of convergence which is independent of the dimension of the
images.

    \vspace*{0.2cm}


\vspace*{0.2cm}


\section{Introduction}
\label{se1}


\subsection{Scope of the paper}
\label{se1sub1}

In image classification the task is to learn the functional relationship between input and output, where the input consists of observed images and the output represents classes of the corresponding images that describe what kind of objects are present in the images.
Since many years the
most successful approaches in the area of image classification
are based on deep convolutional neural networks (CNNs), see, e.g., Krizhevsky, Sutskever and Hinton (2012), LeCun, Bengio and Hinton (2015) and Rawat and Wang (2017). Recently, it has been shown that CNN image classifiers that minimize empirical risk are able to achieve dimension reduction (see Kohler, Krzy\.zak and Walter (2022), Kohler and Langer (2025), Walter (2021) and Kohler and Walter (2023)). However, in practice, it is not  possible to compute the empirical risk minimizer. Instead, a gradient descent approach based on smooth surrogate losses and over-parameterized networks having many more trainable parameters than training samples is used.

In Kohler, Krzy\.zak and Walter (2023) a plug-in classifier based on
convolutional networks, which is learned by gradient descent, has been
analyzed. The main result there was that this classifier achieves
a dimension reduction in an average-pooling model, however in contrast
to the results above for the estimates based on empirical risk
minimization the model there does neither use a (more realistic) max-pooling
model nor any kind of hierarchical structure.

In the present paper we consider the case of large datasets such as
ImageNet, which make use of gradient descent prohibitively expensive.
To be able to deal with such large data sets, we define the estimate
by using stochastic gradient descent. In addition, we consider surrogate logistic loss and hierarchical models with max-pooling. Here we show dimensionality reduction and independence of rates from image dimensions. 

\subsection{Pattern recognition}
\label{se1sub2}
We study image classifiers in the context of pattern recognition.
 Let
$d_1,d_2 \in \N$ and let
$(\bX,Y)$, $(\bX_1,Y_1)$, \dots, $(\bX_n,Y_n)$
be independent and identically distributed random variables
with values in
\[
  [0,1]^{d_1 \times d_2}
 \times \{-1,1\}.
\]
Here we use the notation
\[
  [0,1]^{d_1 \times d_2}
  =
[0,1]^{
\{1, \dots, d_1\} \times \{1, \dots, d_2\}
  }
\]
and
\[
[0,1]^J
= \left\{
(a_j)_{j \in J}
\, : \,
a_j \in [0,1] \quad (j \in J)
\right\}
\]
for a nonempty and finite index set $J$, and we describe a (random)
image from (random) class $Y \in \{-1,1\}$ by a (random) matrix $X$
with $d_1$ columns and $d_2$ rows, which contains at position $(i,j)$
the grey scale value of the pixel of the
image at the corresponding position. Our aim is to predict $Y$ given $X$.
More precisely, given the data set
\[
\D_n = \{ (X_1,Y_1), \dots, (X_n,Y_n)\}
\]
the goal is to construct a classifier
\[
\hat{C}_n(\cdot) = \hat{C}_n(\cdot,\D_n):[0,1]^{d_1 \times d_2} \rightarrow \{-1,1\}
\]
such that the misclassification probability
\[
\PROB \{ \hat{C}_n(X) \neq Y | \D_n\}
\]
is as small as possible.

Let
\begin{equation}
\label{se1eq1}
\eta(x) = \PROB\{ Y=1|X=x\}
\quad
( x \in [0,1]^{d_1 \times d_2})
\end{equation}
be the so--called a posteriori probability of class 1. Then
\[
f^*(x)=
\begin{cases}
  1, & \mbox{if } \eta(x) \geq \frac{1}{2} \\
  -1, & \mbox{elsewhere}
  \end{cases}
\]
is the so--called Bayes classifier, i.e., it satisfies
\[
\PROB\{ f^*(X) \neq Y\}
=
\min_{f : [0,1]^{d_1 \times d_2} \rightarrow \{0,1\}}
\PROB\{f(X) \neq Y\}
\]
(cf., e.g., Theorem 2.1 in Devroye, Gy\"orfi and Lugosi (1996)).

In this paper we derive upper bounds on
\begin{eqnarray}
  \label{se1eq2}
  &&
  \EXP \left\{
\PROB \{ \hat{C}_n(X) \neq Y | \D_n\}
-
\PROB\{f^*(X) \neq Y\}
\right\}
\nonumber \\
&&
=
\PROB \{ \hat{C}_n(X) \neq Y \}
-
\min_{f : [0,1]^{d_1 \times d_2} \rightarrow \{0,1\}}
\PROB\{f(X) \neq Y\}
.
\end{eqnarray}

\subsection{Main results}
\label{se1sub3}

We define deep convolutional neural network estimates by minimizing
the empirical logistic loss of a linear combination of networks via stochastic
gradient descent. Here we use a projection step on the weights
in order to ensure that we can control the over-parametrization
of the estimate. We use this estimate to define an image
classifier $\hat{C}_n$.

We show, that in case that the a posteriori probability $\eta(\bx) = \PROB\{Y=1|X = \bx\}$ satisfies a $(p,C)$--smooth hierarchical max-pooling model of finite level $l$ and $supp(\PROB_{\bX}) \subseteq [0,1]^{d_1 \times d_2}$,
we have
\[
  \PROB \left\{
Y \neq \hat{C}_n(\bX) \right\}
-
\PROB \left\{
Y \neq f^*(\bX) \right\}
\leq
c_1\cdot (\log n)^2 \cdot n^{-\min\{\frac{p}{4p+8},
   \frac{1}{8} \}}. 
\]
And if, in addition, 
\[
\PROB \left\{
\max\left\{
\frac{\eta(X)}{1-\eta(X)},
\frac{1-\eta(X)}{\eta(X)}
\right\}
 > n^{\frac{1}{4}}\right\} \geq
1 - \frac{1}{n^{\frac{1}{4}}}
\quad (n \in \N)
\]
holds (which implies that with high probability $\eta(X)=\PROB\{Y=1|X\}$
is either close to one or close to zero),
then we show that the estimates achieve the improved rate of convergence
  \[
  \PROB \left\{
Y \neq \hat{C}_n(\bX)\right\}
-
\PROB \left\{
Y \neq f^*(\bX) \right\}
\leq
c_2 \cdot (\log n)^{4} \cdot n^{-\min\{\frac{p}{2p+4},
   \frac{1}{4} \}}. 
  \]
  In order to prove these results we derive a general result
  which gives an upper bound on the expected logistic loss of an
  over-parametrized linear combination of deep convolutional
  networks learned
  by minimizing an empirical logistic loss via stochastic gradient descent.

\subsection{Discussion of related results}
\label{se1sub4}

Stochastic gradient descent has been proposed in Robbins and Monroe (1951) and further discussed in Nemirovsky et al. (2008), Polyak and Yuditsky (1992), Spall (2003) and  Kushner and Yin (2003). It is an efficient alternative to the standard batch gradient descent (GD) which has high computational complexity owing it to using all large-scale
data stored in memory in each iteration and not allowing online updates. In SGD one sample is used to randomly update 
the gradient in each iteration, instead of directly calculating the exact value of the gradient. SGD is an unbiased estimate of the real gradient, its cost does not depend on the number of samples and it converges with sublinear rate, but achieves the optimal rate for convex problems, cf., e.g., Nemirovsky
et al. (2008). SGD algorithms have been used in many classical machine learning problems such as perceptron, k-means, SVM and lasso, see Bottou (2012). Many improvements of classical SGD have been introduced over the years. They include momentum, Nesterov Accelerated GD, Adaptive Learning Rate Method, Adaptive Moment Estimation (ADAM), Stochastic Average Gradient, Stochastic Variance Reduction Gradient and Altering Direction Method of Multipliers, see Sun et al. (2019)
for a comprehensive survey.  Asymptotic and finite-sample properties of estimators based on stochastic gradients
were investigated by Toulis and Airaldi (2017). Statistical inference for model parameters in SGD has been discussed in Chen et al. (2020). An excellent survey of optimization methods for large-scale machine learning including SGD is provided in Bottou, Curtis and Nocedal (2018).

In recent years much attention has been devoted to properties of deep neural network estimates.
There exist quite a few approximation
results for neural networks (cf., e.g.,
Yarotsky (2018),
Yarotsky and Zhevnerchute (2019),
Lu et al. (2020), Langer (2021) and the literature cited therein). 
Generalization abilities of deep neural networks can either
be analyzed within the framework of the classical VC theory
(using e.g. the result of
Bartlett et al. (2019) to bound the VC dimension of classes
of neural networks) or in case of over-parametrized deep
neural networks
(where the number of free parameters adjusted
to the observed data set is much larger than the sample size)
by using bounds on the Rademacher complexity
(cf., e.g., Liang, Rakhlin and Sridharan (2015), Golowich, Rakhlin and
Shamir (2019), Lin and Zhang (2019),
Wang and Ma (2022)
and the literature cited therein).

Combining such results leads to a rich theory
showing that owing to the network structure the least squares
neural network estimates can achieve suitable dimension
reduction in hierarchical composition models for the function to be
estimated. For a simple model this was first shown by Kohler and Krzy\.zak
(2017) for H\"older smooth function and later extended to arbitrary
smooth functions by Bauer and Kohler (2019). For a more complex
hierarchical composition model and the ReLU activation function this
was shown in Schmidt-Hieber (2020) under the assumption that the networks
satisfy some sparsity constraint. Kohler and Langer (2021) showed
that this also possible for fully connected neural networks, i.e.,
without imposing a sparsity constraint on the network.
Adaptation of deep neural network to especially
weak smoothness assumptions was shown in
Imaizumi and Fukamizu (2018), Suzuki (2018) and Suzuki and Nitanda (2019).

Less well understood is the optimization of deep neural networks.
As was shown, e.g., in Zou et al. (2018), Du et al. (2019),
Allen-Zhu, Li and Song (2019) and Kawaguchi and Huang (2019) the
 application of gradient descent to over-parameterized deep
neural networks 
leads to a neural network which (globally) minimizes the empirical risk
considered. However, as was shown in Kohler and Krzy\.zak (2021),
the corresponding estimates do not behave well on new independent
data. So the main question is why gradient descent (and its variants
like stochastic gradient descent) can be used to fit a neural
network to observed data in such a way that the resulting estimate
achieves good results on new independent data.
 The challenge here is not only to analyze optimization
but to consider it simultaneously with approximation and generalization.

In case of shallow neural networks (i.e., neural networks with only
one hidden layer) this has been done successfully in Braun et al. (2024).
Here it was possible to show that the classical dimension free
rate of convergence of Barron (1994) for estimation of a regression function
where its Fourier transform has a finite moment can also  be achieved by
shallow neural networks learned by gradient descent. The main idea
here is that the gradient descent  selects a subset of the neural
network where random initialization of the inner weights leads
to values with good approximation properties, and that it adjusts
the outer weights for these neurons properly. A similar idea
was also applied in Gonon (2021).
Kohler and Krzy\.zak
(2022) applied this idea in the context of over-parametrized deep
neural networks where a linear combination of a huge number of 
deep neural networks of fixed size are computed in parallel.
Here the gradient descent selects again a subset of the
neural networks computed in parallel and chooses a proper linear
combination of the networks. By using metric entropy bounds
(cf., e.g.,
Birman and Solomnjak (1967) and
Li, Gu and Ding (2021)) it is possible to control generalization of
the over-parametrized neural networks, and as a result the rate of
convergence of order close to $n^{-1/(1+d)}$ (or $n^{1/(1+d^*)}$ in case
of interaction models, where it is assumed that the regression function
is a sum of functions applied to only $d^*$ of the $d$ components
of the predictor variable) can be shown
for H\"older-smooth regression functions
with H\"older exponent $p \in [1/2,1]$. Universal consistency of
such estimates for bounded $X$ was shown in Drews and Kohler (2022).

In all those results adjusting the inner weights with gradient
descent is not important. In fact, Gonon (2021)
does not do this at all, while Braun et al. (2024) and Kohler and Krzy\.zak
(2022) use the fact that the relevant inner weights do not move too far away from
their starting values during gradient descent.
Similar ideas have also been applied in
Andoni et al. (2014) and Daniely (2017).
This whole approach is related to random feature networks
(cf., e.g., Huang, Chen and Siew (2006) and Rahimi and Recht (2008a, 2008b, 2009)),
where the inner weights are chosen randomly and
only the outer weights are learned during gradient descent.
Yehudai and Shamir (2022) present a lower bound which implies
that either the number of neurons or the absolute value
of the coefficients must grow exponential
in the dimension in order to learn a single ReLU neuron with
random feature networks. But since Braun et al. (2024) was able to
prove a useful rate of convergence result for networks similar
to random feature networks, the practical relevance of this lower
bound is not clear.

The estimates in Kohler and Krzy\.zak (2022) use a $L_2$
regularization on the outer weights during gradient descent.
As was shown in Drews and Kohler (2023), it is possible to
achieve similar results without $L_2$ regularization.

Often gradient descent in neural networks is studied in the
neural tangent kernel setting proposed by Jacot, Gabriel and Hongler (2020),
where instead of a neural network
estimate a kernel estimate is studied and its error is used
to bound the error of the neural network estimate. For further
results in this context see Hanin and Nica (2019) and the literature
cited therein.
Suzuki and Nitanda (2019) were able to analyze the global error of 
an over-parametrized shallow neural network
learned by gradient descent based on this approach. However, due to the
use of the neural tangent kernel, also the smoothness assumption
of the function to be estimated has to be defined with the aid of
a norm involving the kernel, which does not lead to the classical
smoothness conditions of our paper. Another approach where the estimate
is studied  in some asymptotically equivalent model
is the  mean field approach, cf., Mei, Montanari, and Nguyen (2018), Chizat and Bach (2018) or Nguyen and Pham (2020).
A survey of various results on over-parametrized deep neural network
estimates learned by gradient descent can be found in
Bartlett, Montanari and Rakhlin (2021).

In recent years deep transformer networks became very popular in research and applications. 
They have been introduced by Vaswani et al. (2017) and their approximation and generalization 
properties have been investigated by Gurevych et al. (2022). The rates of convergence of over-parametrized transformer classifiers learned by gradient descent have been studied by Kohler and Krzy\.zak (2023).

\subsection{Notation}
\label{se1sub5}
  The sets of natural numbers, real numbers and nonnegative real numbers
  are denoted by $\N$, $\R$ and $\R_+$, respectively. We define
  furthermore $\bar{\R}=$ $\R \cup\{-\infty, \infty\}$. For $z \in \R$, we denote
the smallest integer greater than or equal to $z$ by
$\lceil z \rceil$,
the largest integer less than or equal to $z$ by
$\lfloor z \rfloor$,
and we set $z_+=\max\{z,0\}$.
The Euclidean norm of $x \in \Rd$
is denoted by $\|x\|$. For a closed and convex set $A \subseteq \R^d$
we denote by $Proj_A x$ that element $Proj_A x \in A$ with
\[
\|x-Proj_A x\|= \min_{z \in A} \|x-z\|.
\]
For $f:\R^d \rightarrow \R$
\[
\|f\|_\infty = \sup_{x \in \R^d} |f(x)|
\]
is its supremum norm, and we set
\[
\|f\|_{\infty,A}
= \sup_{x \in A} |f(x)|
\]
for $A \subseteq \Rd$.

For $\bj=(j^{(1)},\dots,j^{(d)}) \in \N_0^d$ we write
\[
\|\bj\|_1=j^{(1)}+ \dots + j^{(d)}
\]
and for $f:\R^d \rightarrow \R$ we set
\[
\partial^{\bj} f =
\frac{\partial^{\|\bj\|_1} f}{(\partial x^{(1)})^{j^{(1)}} \dots
  (\partial x^{(d)})^{j^{(d)}}}.
\]
Let $\F$ be a set of functions $f:\Rd \rightarrow \R$,
let $x_1, \dots, x_n \in \Rd$, set $x_1^n=(x_1,\dots,x_n)$ and let
$p \geq 1$.
A finite collection $f_1, \dots, f_N:\Rd \rightarrow \R$
  is called an $L_p$ $\varepsilon$--packing in $\F$ on $x_1^n$
  if $f_1, \dots, f_N \in \F$ and
  \[
  \min_{1 \leq i<j\leq N}
  \left(
  \frac{1}{n} \sum_{k=1}^n |f_i(x_k)-f_j(x_k)|^p
  \right)^{1/p} \geq \varepsilon
  \]
  hold.
  The $L_p$ $\varepsilon$--packing number of $\F$ on $x_1^n$
  is the  size $N$ of the largest $L_p$ $\varepsilon$--packing
  of $\F$ on $x_1^n$ and is denoted by $\Mu_p(\varepsilon,\F,x_1^n)$.

For $z \in \R$ and $\beta>0$ we define
$T_\beta z = \max\{-\beta, \min\{\beta,z\}\}$. If $f:\R^d \rightarrow
\R$
is a function  then we set
$
(T_{\beta} f)(x)=
T_{\beta} \left( f(x) \right)$. And $sign(z)$ is the sign of $z \in \bar{\R}$.

\subsection{Outline}
\label{se1sub6}
 
A general result, namely a bound on the logistic risk of an
over-parametrized deep convolutional network fitted to data via
stochastic
gradient descent, is presented in Section \ref{se2}. The over-parametrized
deep convolutional neural network classifiers considered
in this paper are introduced in Section \ref{se3} and a bound on their
misclassification probability
is also presented in this section. Section \ref{se4} contains the proofs.

\section{A general result}
\label{se2}

Let $\bTheta$ be a closed and convex set of parameter values
(weights) for
a deep convolutional network of a given topology.
In the sequel we assume that our aim is to learn the parameter
$\btheta \in \bTheta$ (vector of weights) for a deep convolutional  network
\[
f_\btheta: [0,1]^{d_1 \times d_2} \rightarrow \R
\]
from the data $\D_n$ such that
\[
sign(f_\btheta(x))
\]
is a good classifier. We do this by considering linear
combinations
\begin{equation}
  \label{se2eq1}
f_{(\bw,\btheta)}(x)= \sum_{k=1}^{K_n} w_k \cdot T_{\beta_n}(f_{\btheta_k}(x))
  \end{equation}
of truncated versions of deep convolutional networks $f_{\btheta_k}(x)$ $(k=1, \dots, K_n)$,
where $\bw=(w_k)_{k=1, \dots, K_n}$ satisfies
\begin{equation}
  \label{se2eq2}
  w_k \geq 0 \quad (k=1, \dots, K_n), \quad 
  \sum_{k=1}^{K_n} w_k \leq 1
\quad \mbox{and} \quad
  \sum_{k=1}^{K_n} w_k^2 \leq \alpha_n
  \end{equation}
for some $\alpha_n \in [0,1]$,
where $\btheta=(\btheta_1, \dots, \btheta_{K_n}) \in \bTheta^{K_n}$
and where $\beta_n=c_3 \cdot \log n$.
Observe that by choosing $\alpha_n=\frac{1}{N_n}$, $w_j=\frac{1}{N_n}$
  $(j=1, \dots, N_n)$, $\vartheta_j=\vartheta$
  $(j=1, \dots, N_n)$
and $w_k=0$ for $k>N_n$ we get
\[
f_{(\bw,\btheta)}(x)=T_{\beta_n}(f_{\btheta}(x))
\]
and in this way we can construct an estimate which satisfies
\[
sign(f_{(\bw,\btheta)}(x))
=
sign(f_{\btheta}(x))
\]
for any $\btheta \in \bTheta$. And by choosing $K_n$ very large
our estimate will be over-parametrized in the sense that the number
of parameters of the estimate is much larger than the sample size.

Let
\[
\varphi(z)=\log( 1 + \exp(-z))
\]
be the logistic loss (or cross entropy loss). Our aim in choosing
$(\bw,\btheta)$ is the minimization of the logistic risk
\[
F((\bw,\btheta))
=
\EXP\left\{
\varphi(Y \cdot f_{(\bw,\btheta)}(X))
\right\}.
\]
In order to achieve this, we start with a random
initialization of $(\bw,\btheta)$: We choose
\begin{equation}
  \label{se2eq3}
\btheta_1^{(0)}, \dots, \btheta_{K_n}^{(0)}
\end{equation}
uniformly from some closed and convex set $\bTheta^0 \subseteq \bTheta$ such that
the random variables in (\ref{se2eq3}) are independent and also
independent from $(X,Y)$, $(X_1,Y_1)$, \dots, $(X_n,Y_n)$, and we set
\[
w_k^{(0)}=0 \quad (k=1, \dots, K_n).
\]
Then we perform $t_n \in \N$ stochastic gradient descent steps starting with
\[
\btheta^{(0)}=(\btheta_1^{(0)}, \dots, \btheta_{K_n}^{(0)})
\quad \mbox{and} \quad \bw^{(0)}=(w_1^{(0)}, \dots, w_{K_n}^{(0)}).
\]
Here we assume that $t_n/n$ is a natural number, and for
$s \in \{1, \dots, t_n/n\}$ we let
\[
j_{(s-1) \cdot n}, \dots, j_{s \cdot n-1}
\]
be an arbitrary permutation of $1, \dots, n$,
we choose a stepsize $\lambda_n>0$ and we set
\begin{eqnarray*}
  \bw^{(t+1)}
  &=&
  Proj_A \left(
 \bw^{(t)} - \lambda_n \cdot \nabla_\bw \varphi \left( Y_{j_t} \cdot 
f_{( \bw^{(t)}, \vartheta^{(t)})}(X_{j_t})\right)
 \right), \\
 \btheta^{(t+1)}
 &=&
 Proj_B\left(
 \btheta^{(t)} - \lambda_n \cdot \nabla_\btheta \varphi \left( Y_{j_t} \cdot 
f_{( \bw^{(t)}, \vartheta^{(t)})}(X_{j_t})\right)
 \right) \\
\end{eqnarray*}
for $t=0, \dots, t_n-1$.
Here $A$ is the set of all $\bw$ which satisfy (\ref{se2eq2}), and
\[
B= \left\{
\btheta \in \bTheta^{K_n} \, : \, \|\btheta-\btheta^{(0)}\| \leq 1
\right\},
\]
and $Proj_A$ and $Proj_B$ is the $L_2$ projection on $A$ and $B$.
Our estimate is then defined by
\begin{equation}
  \label{se2eq4}
f_n(x)= f_{(\hat{\bw},\vartheta^{(t_n)})}(x)
  \end{equation}
where
\begin{equation}
\label{se2eq5}
  \hat{\bw} = \frac{1}{t_n} \cdot \sum_{t=0}^{t_n-1} \bw^{(t)}.
\end{equation}

Our main result in this general setting is the following bound on the
logistic risk of the above estimate.
\begin{theorem}
  \label{th1}
  Let $(X,Y)$, $(X_1,Y_1), \ldots, (X_n,Y_n)$
  be independent and identically distributed random variables
  with values in $[0,1]^{d_1 \times d_2} \times \{-1,1\}$. 
Let
$N_n,I_n, t_n \in \N$ and let $C_n,D_n \geq 0$.
Set
  $\beta_n = c_3 \cdot \log n$,
  \[
\alpha_n=\frac{1}{N_n}, \quad
\lambda_n=\frac{1}{t_n}, \quad K_n = N_n \cdot I_n
\]
and define the estimate $f_n$ as above.

Assume that there exists $\tilde\vartheta\in \bTheta^0$, such that $f_{\tilde\vartheta}(X)=0$, let
$\bTheta^* \subset \bTheta^0$ and set
\[
\bar{\bTheta}= \left\{ \btheta \in \bTheta \; : \,
\inf_{\tilde{\btheta} \in \bTheta^0} \|\btheta - \tilde{\btheta}\| \leq 1 \right\}.
\]
Assume
\begin{equation}
  \label{th1eq1}
  \| f_\btheta - f_{\bar{\btheta}} \|_{\infty,supp(X)}
  \leq
  C_n \cdot \| \btheta - \bar{\btheta} \|_\infty
  \quad
  \mbox{for all }
  \btheta, \bar{\btheta} \in \bar{\bTheta},
\end{equation}
\begin{equation}
  \label{th1eq2}
  \kappa_n = \PROB \left\{
\btheta^{(0)}_1 \in \bTheta^*
\right\}
>0,
\end{equation}
\begin{equation}
  \label{th1eq3}
N_n \cdot (1- \kappa_n)^{I_n} \leq \frac{1}{n} 
\end{equation}
and
\begin{equation}
  \label{th1eq4}
  \| \nabla_{\bw} \varphi( y \cdot f_{(\bw,\btheta)}) \| \leq D_n
\end{equation}
for all 
  $x \in [0,1]^{d_1 \times d_2}$, $y \in \{-1,1\}$,
$\bw \in A$, $\btheta \in \bar{\bTheta}$, $t \in \{0, \dots, t_n-1\}$.

Then we have
\begin{eqnarray*}
  &&
  \EXP \left\{ \varphi(Y \cdot f_n(X)) \right\}
  -
  \min_{f : [0,1]^{d_1 \times d_2} \rightarrow \bar{\R}}
  \EXP \left\{ \varphi(Y \cdot f(X)) \right\}
  \\
  &&
  \leq
  c_4 \cdot  \Bigg(
  \frac{\log n}{n}
  +
  \beta_n \cdot
  \sup_{x_1, \dots, x_n \in [0,1]^{d_1 \times d_2}}
  \EXP \left\{
  \left|
  \sup_{\btheta \in \bar{\bTheta}}
  \frac{1}{n} \sum_{i=1}^n \epsilon_i \cdot T_{\beta_n} f_\btheta(x_i)
  \right|
  \right\}
 +
  \frac{C_n +1}{\sqrt{N_n}}
  +
  \frac{D_n^2}{t_n}
\\
&&
\quad
+
\frac{n \cdot \left(
6 \cdot K_n \cdot  \beta_n^2 
+
4 \cdot (\beta_n +1) \cdot C_n \cdot
\sup_{{(\bw,\vartheta) \in}\W, y \in \{-1,1\}, x \in [0,1]^{d_1 \times d_2}} \| \nabla_{\vartheta}
   \varphi( y \cdot f_{(\bw,\vartheta)}(x)) \|_\infty
\right)
}{t_n}
  \\
  &&
  \quad
  +
  \sup_{\btheta \in \bTheta^*}
  \EXP \left\{ \varphi(Y \cdot T_{\beta_n} f_\btheta(X)) \right\}
  -
  \min_{f : [0,1]^{d_1\times d_2} \rightarrow \bar{\R}}
  \EXP \left\{ \varphi(Y \cdot f(X)) \right\}
  \Bigg),
\end{eqnarray*}
where $\epsilon_1$, \dots, $\epsilon_n$ are independent and
uniformly distributed
on $\{-1,1\}$ (so-called Rademacher random variables).
  \end{theorem}

\noindent
{\bf Remark 1.}
In Theorem 2 we use the Rademacher complexity
$$
\EXP\left\{\sup _{\vartheta \in \boldsymbol{\Theta}}\left|\frac{1}{n} \sum_{i=1}^n \epsilon_i \cdot T_{\beta_n}\left(f_{\vartheta}\left(x_i\right)\right)\right|\right\}
$$
 to control the generalization error of the estimate. The approximation error is measured by the term
$$
\sup _{\vartheta \in \boldsymbol{\Theta}^*} \EXP\left\{\varphi\left(Y \cdot f_{\vartheta}(X)\right)\right\}-\min _{f: \mathbb{R}^{d \cdot l} \rightarrow \overline{\mathbb{R}}} \EXP\{\varphi(Y \cdot f(X))\},
$$
which describes the maximal error occuring in the set $\Theta^*$.
The last three terms are used to bound the error which occurs during optimization, i.e. due to stochastic gradient descent.\\
\noindent
~\\
    {\bf Remark 2.}
    Our result above extends Theorem 2 in Kohler and Krzy\.zak (2023)
    from gradient descent to the case of stochastic gradient descent.
    To be able to do this we need in the definition of the
    estimate an additional $L_2$ penalty
    on the weights in the linear combination of the
    networks (depending on $\alpha_n$). Furthermore, assumption (\ref{th1eq1})
    is substantially stronger than the corresponding assumption
    in Theorem 2 in Kohler and Krzy\.zak (2023), because there
    it is only required that (\ref{th1eq1}) holds for networks
    which have good approximation properties (which is because
    of the maximal attention used in Transformer networks crucial
    for Transformer networks).

\section{Image classification using deep convolutional neural networks}
\label{se3}

\subsection{Convolutional neural network classifiers}
\label{se3sub2}
We aim to learn feature representations of the inputs by means of $L$ (hidden) convolutional layers. Each of these $r\in \{1, \ldots, L\}$ feature maps
consists of $k_r$ channels.
The input image is considered as layer $0$ with only one channel, i.e. $k_0 = 1$.\\
A convolution in layer $r$ is performed by using a window of values of the previous layer $r-1$ of size $M_r$. The window has to fit within the dimensions of the input image, i.e. $M_r\leq \min\{d_1,d_2\}$.  It relies on so-called filters, i.e. a weight matrix that determines how a neuron is computed from a weighted sum of neighboring neurons from the previous layer.
The weight matrix is defined by 
\begin{align*}
\bold{w} = \left(w_{i,j, s_1, s_2}^{(r)}\right)_{1 \leq i, j \leq M_r, s_1 \in \{1, \dots, k_{r-1}\}, s_2 \in \{1, \dots, k_r\}, r \in \{1, \dots, L\}}.
\end{align*}
Furthermore we need some weights
\begin{align*}
\bold{w}_{bias} = (w_{s_2}^{(r)})_{s_2 \in \{1, \dots, k_r\}, r \in \{1, \dots, L\}}
\end{align*}
for the bias in each channel and output weights 
\begin{align*}
\bold{w}_{out} = (w_{s})_{s \in \{1, \dots, k_L\}},
\end{align*}
which are required for the max-pooling layer defined below.\\
 In the following the ReLU function $\sigma(x) = \max\{x, 0\}$ is chosen as activation function.
The value of a feature map in the $s_2$-th channel of layer $r$ at the position $(i,j)$ is recursively defined by:
\begin{align}
\label{se3eq2}
o_{(i,j), s_2}^{(r)} = \sigma\left(\sum_{s_1=1}^{k_{r -1}} \sum_{\substack{t_1, t_2 \in \{1, \dots, M_r\} \\ (i+t_1-1, j+t_2-1) \in D}} w_{t_1, t_2, s_1, s_2}^{(r)} o_{(i+t_1-1, j+t_2-1), s_1}^{(r-1)} + w_{s_2}^{(r)}\right), 
\end{align}
where  $(i,j)\in D=\{1, \dots, d_1\} \times \{1, \dots, d_2\}$, $s_2\in\{1,\ldots,k_r\}$ and $r\in\{1, \ldots, L\}$.\\
The anchor case $r=0$ of this recursion reflects the values of the input image
\begin{align*}
o_{(i,j), 1}^{(0)} = x_{i,j} \quad \text{for} \ i \in \{1, \dots, d_1\} \ \text{and} \ j \in \{1, \dots, d_2\}.
\end{align*}
In definition \eqref{se3eq2} above we see that weights generating the feature map $o^{(r)}_{(:,:),s2}$ are shared. Weight sharing is used to reduce model complexity, thereby increasing the network's computational efficiency. 
In the last step a max-pooling layer is applied to the values in the $k_L$ channels of the last convolutional layer $L$, such that the output of the network is given by a real-valued function on $[0,1]^{\{1, \dots, d_1\} \times \{1, \dots, d_2\}}$ of the form
\begin{align*}
f_{\bold{w}, \bold{w}_{bias}, \bold{w}_{out}}(x) = &\max\bigg\{\sum_{s_2=1}^{k_L} w_{s_2} \cdot o_{(i,j), s_2}^{(L)}: i \in \{1, \dots, d_1-M_L+1\},\\
& \hspace*{7cm}  j \in \{1, \dots, d_2-M_L+1\}\bigg\}.
\end{align*}
Our class of convolutional neural networks with parameters $L$, $\bold{k} = (k_1, \dots, k_L)$ and $\bold{M} = (M_1, \dots, M_L)$ is defined by 
$\mathcal{F}_{L, \bold{k}, \bold{M}}^{CNN}$. 
As in Kohler, Krzy\.zak and Walter (2022), the definition of the summation index over $t_1, t_2 \in \{1, \ldots,M_r\}$, such that $1\leq i+t_1-1\leq d_1$ and $1\leq j+t_2-1\leq d_2$, corresponds to zero padding to the left and to the bottom of the image. 
Thus, the size of a channel is the same as in the previous layer (see Kohler, Krzy\.zak and Walter (2022) for a further illustration).
Our final estimate is a composition of a convolutional neural network
out of the class $\mathcal{F}_{L, \bold{k}, \bold{M}}^{CNN}$ and a
shallow neural network, which is defined as follows: The output of this network is produced by a function $g:\R \rightarrow\R$ of the form
\begin{equation}
g(x)=\sum_{i=1}^{L_n^{(2)}} w_{i}^{(1) } \sigma \left(
w_{i,1}^{(0)} \cdot x + w_{i,0}^{(0)}
\right)+w_{0}^{(1)},
\label{eq1}
\end{equation}
where $w_{0}^{(1)}, w_{1}^{(1)}, w_{1,0}^{(1)}, w_{1,1}^{(1)}, \dots,
w_{L_n^{(2)}}^{(1)},
w_{L_n^{(2)},0}^{(0)}, w_{L_n^{(2)},1}^{(0)}\in\R$ denote the weights
of this network and $\sigma(z)=\max\{z,0\}$ is again the ReLU activation function.
We define the function class of all real-valued functions on $\R$ of the form \eqref{eq1} with parameter $L_n^{(2)}$ by $\mathcal{F}^{FNN}_{ L_n^{(2)}}$.
\\
\\
Our final function class $\mathcal{F}_n$ is then of the form
\begin{equation}
  \label{se3eq*}
  \mathcal{F}_n = \left\{g \circ f: g \in \mathcal{F}^{FNN}_{L_n^{(2)}},
  f \in \mathcal{F}^{CNN}_{L_n^{(1)}, \bold{k}, \bold{M}} \right\}, 
\end{equation}
which depends on the parameters 
\[\bL=(L_n^{(1)},L_n^{(2)}),
  ~\bk=\left(k_1,\dots,k_{L_n^{(1)}}\right),
~ \bM=(M_1,\dots,M_{L_n^{(1)}}).\]

\subsection{Definition of the estimate}
\label{se3sub3}
 Let $\Theta$ be the set of all weights of the function class 
\[
\mathcal{F}_n=\{f_\btheta \, : \, \theta \in \Theta\}
\]
 introduced in the previous subsection. In the sequel
we fit a linear combination
\[
f_{(\bw, \btheta)}(x)=
\sum_{k=1}^{K_n} w_k \cdot T_{\beta_n} f_{\btheta_k}
\]
to the data where $\bw$ satisfies the assumption
(\ref{se2eq2}) and $\btheta
=(\theta_1, \dots, \theta_{K_n}) \in \Theta^{K_n}$.

Depending on some $B_n>0$, which will be defined in Theorem \ref{th2}
below, we choose
\begin{equation}
  \label{se3eq3}
\btheta_1^{(0)}, \dots, \btheta_{K_n}^{(0)}
\end{equation}
uniformly from  
\[
\bTheta^0 = \left\{ \btheta \in \bTheta \quad : \quad
\|\btheta\|_\infty \leq B_n 
\right\}
\] 
such that
the random variables in (\ref{se3eq3}) are independent and also
independent from $(X,Y)$, $(X_1,Y_1)$, \dots, $(X_n,Y_n)$ and we set
\[
w_k^{(0)}=0 \quad (k=1, \dots, K_n).
\]
Then we perform $t_n \in \N$ stochastic gradient descent steps starting with
\[
\btheta^{(0)}=(\btheta_1^{(0)}, \dots, \btheta_{K_n}^{(0)})
\quad \mbox{and} \quad \bw^{(0)}=(w_1^{(0)}, \dots, w_{K_n}^{(0)}).
\]
As in the previous section
we assume that $t_n/n$ is a natural number and for
$s \in \{1, \dots, t_n/n\}$ we let
\[
j_{(s-1) \cdot n}, \dots, j_{s \cdot n-1}
\]
be an arbitrary permutation of $1, \dots, n$.
We choose a stepsize $\lambda_n>0$ and set
\begin{eqnarray*}
  \bw^{(t+1)}
  &=&
  Proj_A \left(
 \bw^{(t)} - \lambda_n \cdot \nabla_\bw \varphi \left( Y_{j_t} \cdot 
f_{( \bw^{(t)}, \vartheta^{(t)})}(X_{j_t})\right)
 \right), \\
 \btheta^{(t+1)}
 &=&
 Proj_B\left(
 \btheta^{(t)} - \lambda_n \cdot \nabla_\btheta \varphi \left( Y_{j_t} \cdot 
f_{( \bw^{(t)}, \vartheta^{(t)})}(X_{j_t})\right)
 \right) \\
\end{eqnarray*}
for $t=0, \dots, t_n-1$.
Here $A$ is the set of all $\bw$ which satisfy (\ref{se2eq2}), and
\[
B=\left\{
\btheta \in (\bTheta^{(0)})^{K_n} \, : \, \|\btheta-\btheta^{(0)}\| \leq 1
\right\},
\]
and $Proj_A$ and $Proj_B$ is the $L_2$ projection on $A$ and $B$.
In order to compute the gradient with respect to the inner weights
we use the following convention: We set
\[
\sigma^\prime(z)
=
\frac{\partial}{\partial z} \max\{z,0\}
=
\begin{cases}
1, & \mbox{if } z \geq 0 \\
0, & \mbox{else}
\end{cases}
\]
\[
\frac{\partial}{\partial z}
T_{\beta_n} z
=
\frac{\partial}{\partial z} \max\{-\beta_n, \min\{\beta_n,z\}\}
=
\begin{cases}
1, & \mbox{if } |z| \leq \beta_n \\
0, & \mbox{else}
\end{cases}
\]
and
\[
\frac{\partial}{\partial \vartheta_j}
\max \left\{
f_{\vartheta_1}(x),
\dots,
f_{\vartheta_L}(x)
\right\}
=
\frac{\partial}{\partial \vartheta_j}
f_{ \vartheta_l}(x)
\]
where $l \in \{1, \dots, L\}$      satisfies
\[
\max \left\{
f_{\vartheta_1}(x),
\dots,
f_{\vartheta_L}(x)
\right\} =
f_{ \vartheta_l}(x)
\]
and
\[
l=1
\quad \mbox{or} \quad
\max \left\{
f_{\vartheta_1}(x),
\dots,
f_{\vartheta_{l-1}}(x)
\right\} <
f_{ \vartheta_l}(x).
\]
Our  classifier $\hat{C}_n(x)$ is then defined by
\begin{equation}
  \label{se3eq5}
\hat{C}_n(x)= sign (f_n(x)),
  \end{equation}
where
\begin{equation}
  \label{se3eq6}
  f_n(x)=
  f_{(\hat{\bw},\vartheta^{(t_n)})}(x)
  \quad \mbox{and} \quad
  \hat{\bw} = \frac{1}{t_n} \cdot \sum_{t=0}^{t_n-1} \bw^{(t)}.
\end{equation}

\subsection{Main result}

It is well known that one needs smoothness assumptions on the
a posteriori probability in order to derive non-trivial rate
of convergence results for the difference between the
misclassification risk of any estimate and the optimal
misclassification risk (cf., e.g., Cover (1968) and
Section 3 in Devroye and Wagner (1982)). For this we will
use our next definition.

\begin{definition}
Let $p=q+s$ for some $q \in \N_0$ and $0< s \leq 1$.
A function $f:\R^d \rightarrow \R$ is called
$(p,C)$-smooth, if for every $\balpha=(\alpha_1, \dots, \alpha_d) \in
\N_0^d$
with $\sum_{j=1}^d \alpha_j = q$ the partial derivative
$\frac{
\partial^q f
}{
\partial x_1^{\alpha_1}
\dots
\partial x_d^{\alpha_d}
}$
exists and satisfies
\[
\left|
\frac{
\partial^q f
}{
\partial x_1^{\alpha_1}
\dots
\partial x_d^{\alpha_d}
}
(\bx)
-
\frac{
\partial^q f
}{
\partial x_1^{\alpha_1}
\dots
\partial x_d^{\alpha_d}
}
(\bz)
\right|
\leq
C
\cdot
\| \bx-\bz \|^s
\]
for all $\bx,\bz \in \R^d$.
\end{definition}

Furthermore we will use a model from Kohler, Krzy\.zak and Walter (2022)
to be able to derive rates of convergence which do not depend
on the dimension $d_1 \cdot d_2$ of the images. In order to be
able to introduce this model, we need the following notation:
For $M \subseteq \R^d$ and $\bx \in \R^d$ we define
\[
\bx+M = \{\bx+\bz \, : \, \bz \in M \}.
\]
For $I \subseteq \{1, \dots, d_1\} \times \{1, \dots, d_2\}$
and
$\bx=(x_i)_{ i \in \{1, \dots, d_1\} \times \{1, \dots, d_2\}}
\in  [0,1]^{\{1, \dots, d_1\} \times \{1, \dots, d_2\}}$ we set
\[
\bx_I =(x_i)_{i \in I}.
\]

The basic idea behind the next definition is that the a posteriori
probability is a maximum of probabilities that special objects
occur in subparts of the image, and that the decision about the
latter events are hierarchically decided.\\
This notion is inspired by the way humans would generally proceed with
a visual recognition task, specifically the task of deciding whether an image contains a
certain object or not. The question whether or not an object can be detected in an image
is usually solved in several sub-tasks in the sense that an individual mentally divides the
image hierarchically into several parts, scanning each of them for the desired object and thereby estimating the likelihood of the object being present in each part. This
idea leads to the hierarchical model introduced in part b) of the definition below.\\
Naturally, this idea then leads to the assumption that the probability that an image contains the
desired object is simply the maximum of the probabilities estimated in each sub-task, i.e. in each subpart of the image, which motivates the introduction the max-pooling model in part a) of definition \ref{de1}.

\begin{definition}
  \label{de1}
	Let $d_1,d_2\in\N$ with $d_1,d_2>1$ and $m: [0,1]^{\{1, \dots, d_1\} \times \{1, \dots, d_2\}} \rightarrow \R$.

\noindent
{\bf a)}
We say that $m$
satisfies a {\bf max-pooling model with index set}
\[
I \subseteq \{0, \dots, d_1-1\} \times \{0, \dots, d_2-1\},
\]
if there exists a function $f:[0,1]^{(1,1)+I}\rightarrow \R$ such that
\[
m(\bx)=
\max_{
  (i,j) \in \Z^2 \, : \,
  (i,j)+I \subseteq \{1, \dots, d_1\} \times \{1, \dots, d_2\}
}
f\left(
\bx_{(i,j)+I}
\right)
\quad
(\bx \in [0,1]^{\{1, \dots, d_1\} \times \{1,
  \dots, d_2\}}).
\]

\noindent
    {\bf b)}
    Let $I=\{0, \dots, 2^l-1\} \times \{0, \dots, 2^l-1\}$
    for some $l \in \N$.
    We say that
\[
f:[0,1]^{\{1, \dots, 2^l\} \times \{1, \dots, 2^l\}} \rightarrow \R
\]
 satisfies a
    {\bf hierarchical model of level $l$},
    if there exist functions
    \[
    g_{k,s}: \R^4 \rightarrow [0,1]
    \quad (k=1, \dots, l, s=1, \dots, 4^{l-k} )
    \]
    such that we have
    \[
f=f_{l,1}
    \]
    for some
    $f_{k,s} :[0,1]^{\{1, \dots, 2^k\} \times \{1, \dots, 2^k\}} \rightarrow \R$ recursively defined by
    \begin{eqnarray*}
    f_{k,s}(\bx)&=&g_{k,s} \big(
    f_{k-1,4 \cdot (s-1)+1}(\bx_{
\{1, \dots, 2^{k-1}\} \times \{1, \dots, 2^{k-1}\}
    })
    , \\
        &&
        \hspace*{1cm}
        f_{k-1,4 \cdot (s-1)+2}(\bx_{
\{1, \dots, 2^{k-1}\} \times \{2^{k-1}+1, \dots, 2^k\}
        }), \\
        &&
        \hspace*{1cm}
        f_{k-1,4 \cdot (s-1)+3}(\bx_{
\{2^{k-1}+1, \dots, 2^k\} \times \{1, \dots, 2^{k-1}\}
        }), \\
        &&
        \hspace*{1cm}
        f_{k-1,4 \cdot s}(\bx_{
\{2^{k-1}+1, \dots, 2^k\} \times \{2^{k-1}+1, \dots, 2^k\}
        })
    \big)
\\
&&
\hspace*{6cm}
\left(
\bx \in
[0,1]^{
\{
1, \dots, 2^k
\}
\times
\{ 1, \dots, 2^k
\}
}
\right)
    \end{eqnarray*}
    for $k=2, \dots, l, s=1, \dots,4^{l-k}$,
    and
    \[
 f_{1,s}(
x_{1,1},x_{1,2},x_{2,1},x_{2,2}
)= g_{1,s}(x_{1,1},x_{1,2},x_{2,1},x_{2,2})
\quad
( x_{1,1},x_{1,2},x_{2,1},x_{2,2} \in [0,1])
 \]
 for $s=1, \dots, 4^{l-1}$.

 \noindent
     {\bf c)}
     We say that
     $m: [0,1]^{\{1, \dots, d_1\} \times \{1, \dots, d_2\}} \rightarrow \R$
     satisfies a {\bf hierarchical max-pooling model of level $l$
      } (where $2^l \leq \min\{ d_1,d_2\}$),
     if $m$ satisfies a max-pooling model with index set
     \[
I=\left\{0, \dots, 2^{{l}}-1\right\} \times \left\{0, \dots, 2^{{l}}-1\right\}
     \]
     and the function
     $f:[0,1]^{(1,1)+I} \rightarrow \R$ in the definition of this
     max-pooling model satisfies a hierarchical model
     with level ${l}$.

    \noindent
  {\bf d)} We say that the hierarchical max-pooling model $m: [0,1]^{\{1, \dots, d_1\} \times \{1, \dots, d_2\}} \rightarrow \R$ of level $l$ is $(p,C)$--smooth if all functions $g_{k,s}$ in the definition of the function $m$ are $(p,C)$--smooth for some $C>0$.
  \end{definition}

Our main result is the following theorem, in which bounds on the difference
between the misclassification probability of our classifier and the
optimal misclassification probability are derived.

\begin{theorem}
  \label{th2}Let $(X,Y)$, $(X_1,Y_1), \ldots, (X_n,Y_n)$
  be independent and identically distributed random variables
  with values in $[0,1]^{d_1 \times d_2} \times \{-1,1\}$. 
       Let $p \geq 1$ and $C>0$ be arbitrary.
Assume that the a posteriori probability $\eta(\bx) = \PROB\{Y=1|X = \bx\}$ satisfies a $(p,C)$--smooth hierarchical
max-pooling model of finite level $l$ and set
$\beta_n = c_3 \cdot \log n$,
\[
L_n^{(1)}= \frac{4^l-1}{3} \cdot \lceil c_5 \cdot n^{2/(2p+4)} \rceil+l
\quad \mbox{and} \quad
L_n^{(2)}= \lceil c_6 \cdot n^{1/4} \rceil, 
\]
\[
M_s = 2^{\pi(s)} \quad (s=1, \dots, L_n^{(1)}),
\]
where the function
$\pi:\{1, \dots, L_n^{(1)} \}
\rightarrow \{1, \dots, l\}$ is defined by
\[
\pi(s)=\sum_{i=1}^l
\mathds{1}_{
  \{
s \geq i+\sum_{r=l-i+1}^{l-1} 4^r \cdot \lceil c_5 \cdot n^{2/(2p+4)} \rceil
  \}
},
\]
choose $\bold{k} =(c_7, \dots, c_7) \in \N^{L_n^{(1)}}$ and
 set
\[
B_n = e^{\sqrt{n}},
\]
assume that $K_n \in \N$ satisfies
\[
\frac{K_n}{e^{2 \cdot n^{1.5}}} \rightarrow \infty \quad (n \rightarrow \infty)
\]
and set
\[
\alpha_n = \frac{1}{n^2 \cdot e^{2 \cdot n}}
\quad \mbox{and} \quad
t_n = \left\lceil n^2 \cdot K_n  \right\rceil.
\]
Define the classifier $\hat{C}_n$
as in Section \ref{se3sub3}.
Assume that the constants $c_3, c_5, c_6, c_7$ are sufficiently large.

\noindent
{\bf a)} There exists a constant $c_8 >0$ such that we have for $n$
sufficiently large
\begin{eqnarray*}
&&
\PROB \left\{
Y \neq \hat{C}_n(\bX) \right\}
-
\PROB \left\{
Y \neq f^*(\bX) \right\}
\leq
c_8 \cdot (\log n)^2 \cdot n^{-\min\{\frac{p}{4p+8},
   \frac{1}{8} \}}. 
\end{eqnarray*}

\noindent
    {\bf b)}
    If, in addition,
    \begin{equation}
      \label{th2eq1}
      \PROB \left\{
\max\left\{
\frac{\eta(X)}{1-\eta(X)},
\frac{1-\eta(X)}{\eta(X)}
\right\}
 > n^{\frac{1}{4}}\right\} \geq
1 - \frac{1}{n^{\frac{1}{4}}}
\quad (n \in \N)
\end{equation}
    holds, then there exists a constant $c_9 > 0$ such that we have
    for 
$n$ sufficiently large
\begin{eqnarray*}
&&
\PROB \left\{
Y \neq \hat{C}_n(\bX)\right\}
-
\PROB \left\{
Y \neq f^*(\bX) \right\}
\leq
c_9 \cdot (\log n)^{4} \cdot n^{-\min\{\frac{p}{2p+4},
   \frac{1}{4} \}}. 
\end{eqnarray*}    

  \end{theorem}

\noindent
    {\bf Remark 3.}
    The rates of convergence above do not depend on the
    dimension $d_1 \cdot d_2$ of the image, hence in case that
    the a posteriori distribution satisfies a hierarchical
    composition model, our estimate is able to circumvent
    the curse of dimensionality.
~\\
\noindent 
{\bf Remark 4}
 In \cite{KoLa20} the authors have shown that CNN
image classifiers from the function class $\mathcal{F}_n $ defined in \eqref{se3eq*} that minimize empirical risk are able to achieve dimension reduction and achieved the same convergence rates as in Theorem \ref{th2} above. However, in practice, it is not possible to compute the empirical
risk minimizer. Instead, we apply stochastic gradient descent to over-parametrized linear combinations of functions from \eqref{se3eq*} and show that this estimate achieves the same rate of convergence as the empirical risk minimizer of Kohler and Langer (2025).
    \section{Proofs}
\label{se4}

\subsection{Proof of Theorem \ref{th1}}

In the proof of Theorem \ref{th1} we will need the
following auxiliary result.

\begin{lemma}
  \label{le1}
  Let $l_1, l_2, t_n \in \N$, let $D_n \geq 0$,
  let $A \subset \R^{l_1}$ 
  be closed and convex, let $B \subseteq \R^{l_2}$ and let $F_t,
  F:\R^{l_1} \times \R^{l_2} \rightarrow \R_+$ $(t=0, \dots, t_n-1)$
  be functions such that for all $t \in \{0, \dots, t_n-1\}$
  \[
u \mapsto F(u,v) \quad \mbox{is differentiable and convex for all } v \in \R^{l_2},
\]
\[
u \mapsto F_t(u,v) \quad \mbox{is differentiable for all  } v \in \R^{l_2},
\]
and
  \begin{equation}
    \label{le1eq1}
\| (\nabla_{u} F_t)(u,v) \| \leq D_n
  \end{equation}
  for all $(u,v) \in A \times B$. 
  Choose $(u_0,v_0) \in A \times B$, let $v_1, \dots, v_{t_n} \in B$ and set
  \[
  u_{t+1} = Proj_A \left(
u_t - \lambda \cdot \left( \nabla_u F_t \right)(u_t,v_t)
  \right) \quad (t=0, \dots, t_n-1),
  \]
where
  \[
\lambda = \frac{1}{t_n}.
\]
Let    $u^* \in A$. Then it holds:
  \begin{eqnarray*}
\frac{1}{t_n} \sum_{t=0}^{t_n-1}   F(u_t,v_t)
    &\leq&
    F(u^*,v_0)
    +
    \frac{1}{t_n} \sum_{t=1}^{t_n-1} | F(u^*,v_t) - F(u^*,v_0)|
    +
    \frac{\|u^*-u_0\|^2}{2}
    +
    \frac{D_n^2}{2 \cdot t_n}
    \\
    &&
    +
    \frac{1}{t_n} \sum_{t=0}^{t_n-1}
    < \left( \nabla_u F \right)(u_t,v_t) - \left( \nabla_u F_t \right)(u_t,v_t),
    u_t-u^*>
    .
    \end{eqnarray*}
\end{lemma}

\noindent
    {\bf Proof.}
    By convexity of $u \mapsto F(u,v_t)$ and because of
    $u^* \in A$ we have
    \begin{eqnarray*}
      &&
      F(u_t,v_t) - F(u^*,v_t)
      \\
      &&
      \leq
      \,
      < (\nabla_u F)(u_t,v_t), u_t-u^* >
      \\
      &&
      =
      \frac{1}{2 \cdot \lambda}
      \cdot 2
      \cdot
      < \lambda \cdot (\nabla_u F_t)(u_t,v_t), u_t-u^* >
          +
    < \left( \nabla_u F \right)(u_t,v_t) - \left( \nabla_u F_t \right)(u_t,v_t),
    u_t-u^*>
      \\
      &&
      =
      \frac{1}{2 \cdot \lambda}
      \cdot
      \left(
      - \| u_t - u^* - \lambda \cdot (\nabla_u F_t)(u_t,v_t)\|^2
      +
      \|u_t-u^*\|^2
      +
      \|  \lambda \cdot (\nabla_u F_t)(u_t,v_t) \|^2
      \right)\\
      && \quad
                +
    < \left( \nabla_u F \right)(u_t,v_t) - \left( \nabla_u F_t \right)(u_t,v_t),
    u_t-u^*>
\\
      &&
      \leq
       \frac{1}{2 \cdot \lambda}
      \cdot
      \left(
      - \| Proj_A (u_t  - \lambda \cdot (\nabla_u F_t)(u_t,v_t)) - u^*\|^2
      +
      \|u_t-u^*\|^2
      +
      \lambda^2 \cdot \|  (\nabla_u F_t)(u_t,v_t) \|^2
      \right)\\
      && \quad
                +
    < \left( \nabla_u F \right)(u_t,v_t) - \left( \nabla_u F_t \right)(u_t,v_t),
    u_t-u^*>
      \\
      &&
      =
       \frac{1}{2 \cdot \lambda}
      \cdot
      \left(
      \|u_t-u^*\|^2
      - \| u_{t+1} - u^*\|^2
      +
      \lambda^2 \cdot \|  (\nabla_u F_t)(u_t,v_t) \|^2
      \right)\\
      && \quad
                +
   < \left( \nabla_u F \right)(u_t,v_t) - \left( \nabla_u F_t \right)(u_t,v_t),
    u_t-u^*>.
      \end{eqnarray*}
    This implies
    \begin{eqnarray*}
&&      \frac{1}{t_n}
      \sum_{t=0}^{t_n-1} F(u_t,v_t)
      -
          \frac{1}{t_n}
      \sum_{t=0}^{t_n-1} F(u^*,v_t)
      \\
      &&
      =
        \frac{1}{t_n}
      \sum_{t=0}^{t_n-1} \left( F(u_t,v_t)
      -
      F(u^*,v_t) \right)
      \\
      &&
      \leq
      \frac{1}{t_n}
      \sum_{t=0}^{t_n-1}
       \frac{1}{2 \cdot \lambda}
      \cdot
      \left(
      \|u_t-u^*\|^2
      - \| u_{t+1} - u^*\|^2
      \right)
      +
      \frac{1}{t_n}
      \sum_{t=0}^{t_n-1}
      \frac{\lambda}{2} \cdot \|  (\nabla_u F_t)(u_t,v_t) \|^2
      \\
      && \quad           +
    \frac{1}{t_n} \sum_{t=0}^{t_n-1}
    < \left( \nabla_u F \right)(u_t,v_t) - \left( \nabla_u F_t \right)(u_t,v_t),
    u_t-u^*>
      \\
      &&
      =
        \frac{1}{2}
\cdot
        \sum_{t=0}^{t_n-1}
      \left(
      \|u_t-u^*\|^2
      - \| u_{t+1} - u^*\|^2
      \right)
      +
      \frac{1}{2 \cdot t_n^2}
      \sum_{t=0}^{t_n-1}
      \|  (\nabla_u F_t)(u_t,v_t) \|^2
            \\
      && \quad           +
    \frac{1}{t_n} \sum_{t=0}^{t_n-1}
    < \left( \nabla_u F \right)(u_t,v_t) - \left( \nabla_u F_t \right)(u_t,v_t),
    u_t-u^*>
      \\
      &&
      \leq
       \frac{ \|u_0-u^*\|^2}{2}
       +
      \frac{1}{2 \cdot t_n^2}
      \sum_{t=0}^{t_n-1}
      \|  (\nabla_u F_t)(u_t,v_t) \|^2       \\
      && \quad
      +
    \frac{1}{t_n} \sum_{t=0}^{t_n-1}
    < \left( \nabla_u F \right)(u_t,v_t) - \left( \nabla_u F_t \right)(u_t,v_t),
    u_t-u^*>. 
      \end{eqnarray*}
    Using the above result and (\ref{le1eq1}) we get
    \begin{eqnarray*}
&&
      \frac{1}{t_n}
      \sum_{t=0}^{t_n-1} F(u_t,v_t)
      \\
      &&
      \leq
          \frac{1}{t_n}
      \sum_{t=0}^{t_n-1} F(u^*,v_t)
      +
      \frac{\|u^*-u_0\|^2}{2}
      +
      \frac{1}{2 \cdot t_n^2}
      \sum_{t=0}^{t_n-1}
      \| (\nabla_{u}F_t)(u_t,v_t) \|^2\\
              &&
      \quad
      +
    \frac{1}{t_n} \sum_{t=0}^{t_n-1}
    < \left( \nabla_u F \right)(u_t,v_t) - \left( \nabla_u F_t \right)(u_t,v_t),
    u_t-u^*>\\
      &&
      \leq
      F(u^*,v_0)
      +
         \frac{1}{t_n}
         \sum_{t=0}^{t_n-1}
         |F(u^*,v_t)-F(u^*,v_0)|
      +
      \frac{\|u^*-u_0\|^2}{2}
      +
      \frac{D_n^2}{2 \cdot t_n}
      \\
      &&
      \quad
      +
    \frac{1}{t_n} \sum_{t=0}^{t_n-1}
    < \left( \nabla_u F \right)(u_t,v_t) - \left( \nabla_u F_t \right)(u_t,v_t),
    u_t-u^*>.
      \end{eqnarray*}
    \hfill $\Box$

    \noindent
        {\bf Proof of Theorem \ref{th1}.}
        Let $E_n$ be the event that there exist pairwise distinct
        $j_1, \dots, j_{N_n} \in \{1, \dots, K_n\}$ such that
        \[
\btheta_{j_i}^{(0)} \in \bTheta^*
       \]
       holds for all $i=1, \dots, N_n$. If $E_n$ holds set
       \[
       w_{j_i}^*= \frac{1}{N_n} \quad (i=1, \dots, N_n)
       \quad \mbox{and} \quad
       w_k^*=0 \quad (k \in \{1, \dots, K_n\} \setminus \{j_1, \dots, j_{N_n} \})
       \]
       and $\bw^*=(w_k^*)_{k=1, \dots, K_n}$, otherwise set $\bw^*=0$.
       
       We will use the following error decomposition:
       \begin{eqnarray*}
         &&
           \EXP \left\{ \varphi(Y \cdot f_n(X)) \right\}
  -
  \min_{f : [0,1]^{d_1 \times d_2} \rightarrow \bar{\R} }
  \EXP \left\{ \varphi(Y \cdot f(X)) \right\}
  \\
  &&
  = \EXP \left\{ \varphi(Y \cdot f_n(X)) \cdot 1_{E_n^c}\right\}
\\
&&
\quad
+
\EXP \left\{ 
\EXP \left\{  
\varphi(Y \cdot f_{(\hat{\bw}, \vartheta^{(t_n)})}(X)) 
\big| \btheta^{(0)}, \D_n \right\}
\cdot 1_{E_n}
\right\}
\\
&&
\hspace*{3cm}
-
\EXP \left\{ 
\frac{1}{t_n} \sum_{t=0}^{t_n-1}
\EXP \left\{  
\varphi(Y \cdot f_{(\bw^{(t)}, \vartheta^{(t_n)})}(X)) 
\big| \btheta^{(0)}, \D_n \right\}
\cdot 1_{E_n}
\right\}
\\
&&
\quad
+
\EXP \left\{ 
\frac{1}{t_n} \sum_{t=0}^{t_n-1}
\EXP \left\{  
\varphi(Y \cdot f_{(\bw^{(t)}, \vartheta^{(t_n)})}(X)) 
\big| \btheta^{(0)}, \D_n \right\}
\cdot 1_{E_n}
\right\}
\\
&&
\hspace*{3cm}
-
\EXP \left\{ 
\frac{1}{t_n} \sum_{t=0}^{t_n-1}
\EXP \left\{  
\varphi(Y \cdot f_{(\bw^{(t)}, \vartheta^{(t)})}(X)) 
\big| \btheta^{(0)}, \D_n \right\}
\cdot 1_{E_n}
\right\}
\\
&&
\quad
+
\EXP \left\{ 
\frac{1}{t_n} \sum_{t=0}^{t_n-1}
\EXP \left\{  
\varphi(Y \cdot f_{(\bw^{(t)}, \vartheta^{(t)})}(X)) 
\big| \btheta^{(0)}, \D_n \right\}
\cdot 1_{E_n}
\right\}
-
  \min_{f : [0,1]^{d_1 \times d_2} \rightarrow \bar{\R} }
  \EXP \left\{ \varphi(Y \cdot f(X)) \right\}
  \\
  &&
  =: T_{1,n}+ T_{2,n}+ T_{3,n}+ T_{4,n}.
       \end{eqnarray*}
       In the {\it first step of the proof} we show
       \begin{equation}
         \label{pth1eq1}
         \PROB\{E_n^c\} \leq \frac{1}{n}.
         \end{equation}
       To do this we consider a sequential choice of the
       initial weights $\btheta_1^{(0)}$, \dots, $\btheta_{K_n}^{(0)}$.
       By definition of $\kappa_n$ we know that the probability
       that none of $\btheta_1^{(0)}$, \dots, $\btheta_{I_n}^{(0)}$
       is contained in $\bTheta^*$ is given by
       \[
(1-\kappa_n)^{I_n}.
       \]
       This implies that the probability that there exists $l \in \{1, \dots, N_n\}$
       such that  none of $\btheta_{(l-1)\cdot I_n+1}^{(0)}$,
       \dots, $\btheta_{l \cdot I_n}^{(0)}$
       is contained in $\bTheta^*$ is upper bounded by
       \[
N_n \cdot (1-\kappa_n)^{I_n}.
       \]
       (\ref{th1eq3}) implies
       \[
\PROB\{ E_n^c \} \leq N_n \cdot (1-\kappa_n)^{I_n} \leq \frac{1}{n}.
\]

In the {\it second step of the proof} we show
\[
T_{1,n} \leq c_{10} \cdot \frac{(\log n)}{n}.
\]
To do this, we observe that for $|z| \leq \beta_n$ we have
\[
\varphi(z)=\log( 1 + \exp(-z))
\leq (\log 4) \cdot I_{\{z>-1\}} + \log( 2 \cdot \exp(-z)) \cdot I_{\{z \leq -1\}}
\leq
3 + |z| \leq c_{11} \cdot \log n,
\]
from which we can conclude by the first step of the proof
\[
T_{1,n} \leq c_{11} \cdot (\log n) \cdot \PROB\{E_n^c\}
\leq c_{11} \cdot \frac{\log n}{n}.
\]

In the
{\it third step of the proof} we show
\begin{eqnarray*}
T_{2,n}
&\leq & 0   
    .
\end{eqnarray*}
This follows from the convexity of the logistic loss, which implies
\begin{eqnarray*}
&&
\EXP \left\{  
\varphi(Y \cdot f_{(\hat{\bw}, \vartheta^{(t_n)})}(X)) 
\big| \btheta^{(0)}, \D_n \right\}
\\
&&
=
\EXP \left\{  
\varphi(Y \cdot  \frac{1}{t_n} \sum_{t=0}^{t_n-1}
f_{(\bw^{(t)}, \vartheta^{(t_n)})}(X)) 
\big| \btheta^{(0)}, \D_n \right\}
\\
&&
\leq
\EXP \left\{  
\frac{1}{t_n} \sum_{t=0}^{t_n-1}
\varphi(Y \cdot  
f_{(\bw^{(t)}, \vartheta^{(t_n)})}(X)) 
\big| \btheta^{(0)}, \D_n \right\}
\\
&&
=
\frac{1}{t_n} \sum_{t=0}^{t_n-1}
\EXP \left\{  
\varphi(Y \cdot  
f_{(\bw^{(t)}, \vartheta^{(t_n)})}(X)) 
\big| \btheta^{(0)}, \D_n \right\}.
\end{eqnarray*}

  In the {\it fourth step of the proof} we show
   \begin{eqnarray*}
     &&
T_{3,n} \leq 2 \cdot \frac{C_n}{\sqrt{N}_n}.
\end{eqnarray*}
Due to the fact that the logistic loss is Lipschitz continuous
with Lipschitz constant $1$
and by assumptions (\ref{se2eq2}) and
(\ref{th1eq1})
we have
\begin{eqnarray*}
&&
T_{3,n}
\\
&&
\leq
\frac{1}{t_n} \sum_{t=0}^{t_n-1}
\EXP \left\{  
\left|
\varphi( Y \cdot f_{(\bw^{(t)},\vartheta^{(t_n)})}(X))
-
\varphi( Y \cdot f_{(\bw^{(t)},\vartheta^{(t)})}(X))
\right|
\right\}
\\
&&
\leq
\frac{1}{t_n} \sum_{t=0}^{t_n-1}
\EXP \left\{  
|
f_{(\bw^{(t)}, \vartheta^{(t_n)})}(X) 
-
f_{(\bw^{(t)}, \vartheta^{(t)})}(X) 
|
\right\}
\\
&&
\leq
\frac{1}{t_n} \sum_{t=0}^{t_n-1}
\EXP \left\{  
\sum_{k=1}^{K_n}
w^{(t)}_k
\cdot
|
T_{\beta_n}
 f_{\vartheta_k^{(t_n)}}(X)
-
T_{\beta_n}
 f_{\vartheta_k^{(t)}}(X)
|
\right\}
\\
&&
\leq
\frac{1}{t_n} \sum_{t=0}^{t_n-1}
\EXP \left\{  
\sum_{k=1}^{K_n}
w_k^{(t)}
\cdot
|
 f_{\vartheta_k^{(t_n)}}(X)
-
 f_{\vartheta_k^{(t)}}(X)
 |
\right\}
\\
&&
\leq
\frac{1}{t_n} \sum_{t=0}^{t_n-1}
\EXP \left\{  
\sqrt{\sum_{k=1}^{K_n}
(w_k^{(t)})^2}
\cdot
\sqrt{
\sum_{k=1}^{K_n}
|
 f_{\vartheta_k^{(t_n)}}(X)
-
 f_{\vartheta_k^{(t)}}(X)
|^2
}
\right\}
\\
&&
\leq
\frac{1}{t_n} \sum_{t=0}^{t_n-1}
\EXP \left\{  
\sqrt{\alpha_n}
\cdot
\sqrt{
\sum_{k=1}^{K_n}
C_n^2
\cdot
\|\vartheta_k^{(t_n)}
-
 \vartheta_k^{(t)}
\|_\infty^2
}
\right\}
\\
&&
=
\frac{C_n}{\sqrt{N_n}}\cdot
\frac{1}{t_n} \sum_{t=0}^{t_n-1}
\EXP \left\{  
\sqrt{
\|\vartheta^{(t_n)}
-
 \vartheta^{(t)}
\|^2
}
\right\}
\\
&&
\leq
2 \cdot \frac{C_n}{\sqrt{N_n}}.
\end{eqnarray*}

    In the {\it fifth step of the proof} we apply Lemma \ref{le1}
    to $T_{4,n}$. 
    Set
    \[
    F((\bw,\vartheta))=\EXP\{ \varphi(Y \cdot f_{(\bw,\vartheta)}(X))\big| \btheta^{(0)}, \D_n \}
    \mbox{ and }
      F_t((\bw,\vartheta))=
       \varphi(Y_{j_t} \cdot f_{(\bw,\vartheta)}(X_{j_t})).
    \]
    Then Lemma \ref{le1} implies
    \begin{eqnarray*}
      T_{4,n} &\leq&
      \EXP \left\{
      \EXP \left\{  
\varphi(Y \cdot f_{(\bw^{*}, \vartheta^{(0)})}(X)) 
\big| \btheta^{(0)}, \D_n \right\}
\cdot 1_{E_n}
\right\}
-
  \min_{f : [0,1]^{d_1\times d_2} \rightarrow \bar{\R} }
  \EXP \left\{ \varphi(Y \cdot f(X)) \right\}
  \\
  &&
  +
  \frac{1}{t_n} \sum_{t=1}^{t_n-1}
  \EXP\Bigg\{
  \Bigg|
\EXP\left\{
  \varphi(Y \cdot f_{(\bw^{*}, \vartheta^{(t)})}(X))
\big|
\D_n, \vartheta^{(0)} \right\} 
\\
&&
\hspace*{2.5cm}
  -
\EXP\left\{
  \varphi(Y \cdot f_{(\bw^{*}, \vartheta^{(0)})}(X)) 
\big|
\D_n, \vartheta^{(0)} \right\}
  \Bigg| \cdot 1_{E_n}
  \Bigg\}
  +
\frac{1}{2} \cdot
  \frac{1}{N_n}
  +
  \frac{D_n^2}{2 \cdot t_n}
  \\
  &&
      +
    \frac{1}{t_n} \sum_{t=0}^{t_n-1}
    \EXP \left\{
    < \left( \nabla_\bw F \right)(\bw^{(t)},\vartheta^{(t)}) - \left( \nabla_\bw F_t \right)(\bw^{(t)},\vartheta^{(t)}),
    \bw^{(t)}-\bw^*> 
    \cdot 1_{E_n}
    \right\}
   .
      \end{eqnarray*}

    In the {\it sixth step of the proof} we show
    \begin{eqnarray*}
      &&
      \EXP \left\{
      \EXP \left\{  
\varphi(Y \cdot f_{(\bw^{*}, \vartheta^{(0)})}(X)) 
\big| \btheta^{(0)}, \D_n \right\}
\cdot 1_{E_n}
\right\}
-
  \min_{f : [0,1]^{d_1\times d_2} \rightarrow \bar{\R} }
  \EXP \left\{ \varphi(Y \cdot f(X)) \right\}
  \\
  &&
  \leq
    \sup_{\btheta \in \bTheta^*}
  \EXP \left\{ \varphi(Y \cdot T_{\beta_n} f_\btheta(X)) \right\}
  -
  \min_{f : [0,1]^{d_1\times d_2} \rightarrow \bar{\R}}
  \EXP \left\{ \varphi(Y \cdot f(X)) \right\}
  \Bigg).
    \end{eqnarray*}
    This follows from the convexity of the logistic loss, which implies
    \begin{eqnarray*}
      &&
      \EXP \left\{
      \EXP \left\{  
\varphi(Y \cdot f_{(\bw^{*}, \vartheta^{(0)})}(X)) 
\big| \btheta^{(0)}, \D_n \right\}
\cdot 1_{E_n}
\right\}
\\
&&
=
      \EXP \left\{
      \EXP \left\{  
      \varphi(
      \frac{1}{N_n}
      \sum_{k=1}^{N_n}
      Y \cdot T_{\beta_n} f_{\vartheta_{j_k}^{(0)}}(X)) 
\big| \btheta^{(0)}, \D_n \right\}
\cdot 1_{E_n}
\right\}
\\
&&
\leq
     \frac{1}{N_n}
     \sum_{k=1}^{N_n}
\EXP \left\{
      \EXP \left\{  
      \varphi(
      Y \cdot T_{\beta_n} f_{\vartheta_{j_k}^{(0)})}(X)) 
\big| \btheta^{(0)}, \D_n \right\}
\cdot 1_{E_n}
\right\}
\\
&&
\leq
  \sup_{\btheta \in \bTheta^*}
  \EXP \left\{ \varphi(Y \cdot T_{\beta_n} f_\btheta(X)) \right\}.
      \end{eqnarray*}

    In the {\it seventh step of the proof} we show
    \begin{eqnarray*}
      &&
\frac{1}{t_n} \sum_{t=1}^{t_n-1}
  \EXP\left\{
  \left|
\EXP\left\{
  \varphi(Y \cdot f_{(\bw^{*}, \vartheta^{(t)})}(X))
\big|
\D_n, \vartheta^{(0)} \right\}
  -
\EXP\left\{
  \varphi(Y \cdot f_{(\bw^{*}, \vartheta^{(0)})}(X)) 
\big|
\D_n, \vartheta^{(0)} \right\}
  \right| \cdot 1_{E_n}
  \right\}
\\
&&
  \leq
  \frac{C_n}{\sqrt{N_n}}.
      \end{eqnarray*}
    The Lipschitz continuity of the logistic loss and
    assumption (\ref{th1eq1}) imply
    \begin{eqnarray*}
&&
\frac{1}{t_n} \sum_{t=1}^{t_n-1}
  \EXP\left\{
  \left|
\EXP\left\{
  \varphi(Y \cdot f_{(\bw^{*}, \vartheta^{(t)})}(X))
\big|
\D_n, \vartheta^{(0)} \right\}
  -
\EXP\left\{
  \varphi(Y \cdot f_{(\bw^{*}, \vartheta^{(0)})}(X)) 
\big|
\D_n, \vartheta^{(0)} \right\}
  \right| 
  \cdot 1_{E_n}
  \right\}
\\
      &&
\leq
\frac{1}{t_n} \sum_{t=1}^{t_n-1}
  \EXP\left\{
  \left|
  \varphi(Y \cdot f_{(\bw^{*}, \vartheta^{(t)})}(X))
  -
  \varphi(Y \cdot f_{(\bw^{*}, \vartheta^{(0)})}(X)) 
  \right|
  \right\}
  \\
  &&
  \leq
  \frac{1}{t_n} \sum_{t=1}^{t_n-1}
  \EXP\left\{
  \left|
  f_{(\bw^{*}, \vartheta^{(t)})}(X)
  -
  f_{(\bw^{*}, \vartheta^{(0)})}(X) 
  \right|
  \right\}
  \\
  &&
  =
   \frac{1}{t_n} \sum_{t=1}^{t_n-1}
  \EXP\left\{
  \left|
  \sum_{k=1}^{K_n} w_k^* \cdot
  (
  T_{\beta_n} f_{\vartheta^{(t)}_k}(X)
  -
  T_{\beta_n} f_{\vartheta^{(0)}_k}(X)
  )
  \right|
  \right\}
  \\
  &&
  \leq
    \frac{1}{t_n} \sum_{t=1}^{t_n-1}
  \EXP\left\{
  \sqrt{
    \sum_{k=1}^{K_n} |w_k^*|^2
  }
  \cdot
  \sqrt{ \sum_{k=1}^{K_n}
  (
T_{\beta_n}   f_{\vartheta^{(t)}_k}(X)
  -
 T_{\beta_n}  f_{\vartheta^{(0)}_k}(X)
  )^2
  }
  \right\}
  \\
  &&
  \leq
    \frac{1}{t_n} \sum_{t=1}^{t_n-1}
  \EXP\left\{
  \sqrt{
    \sum_{k=1}^{K_n} |w_k^*|^2
  }
  \cdot
  \sqrt{ \sum_{k=1}^{K_n}
  (
 f_{\vartheta^{(t)}_k}(X)
  -
   f_{\vartheta^{(0)}_k}(X)
  )^2
  }
  \right\}
  \\
  &&
  \leq
  \frac{1}{t_n} \sum_{t=1}^{t_n-1}
  \EXP\left\{
\frac{1}{\sqrt{N_n}}
  \sqrt{ \sum_{k=1}^{K_n}
    C_n^2 \cdot
    \|
  \vartheta^{(t)}_k
  -
  \vartheta^{(0)}_k
  \|_\infty^2
  }
  \right\}
  \\
  &&
  =
   \frac{1}{t_n} \sum_{t=1}^{t_n-1}
   \EXP\left\{
\frac{C_n}{\sqrt{N_n}} \cdot  \|
  \vartheta^{(t)}
  -
  \vartheta^{(0)}
  \|
  \right\}
  \leq
  \frac{C_n}{\sqrt{N_n}}.
      \end{eqnarray*}

 Let $\W$ be the set of all weight vectors $\bw=( (w_k)_{k=1, \dots, K_n},
    (\btheta_k)_{k=1, \dots, K_n})$
    which satisfy $\btheta=(\btheta_k)_{k=1, \dots, K_n} \in \bar{\bTheta}^{K_n}$
    and (\ref{se2eq2}). Let $(X_1^\prime, Y_1^\prime)$, \dots, $(X_n^\prime,Y_n^\prime)$,
    $\epsilon_1$, \dots, $\epsilon_n$ be random variables independent of $\vartheta^{(0)}$ with $\PROB\{ \epsilon_j=1\}=1/2=\PROB\{ \epsilon_j=-1\}$ $(j=1, \dots, n)$
    such that 
    $(X_1, Y_1)$, \dots, $(X_n,Y_n)$,
    $\epsilon_1$, \dots, $\epsilon_n$, $(X_1^\prime, Y_1^\prime)$, \dots, $(X_n^\prime,Y_n^\prime)$
are independent and \\
    $(X_1, Y_1)$, \dots, $(X_n,Y_n)$,
    $(X_1^\prime, Y_1^\prime)$, \dots, $(X_n^\prime,Y_n^\prime)$ are identically distributed. 
    
            In the {\it eighth step of the proof} we show
   \begin{eqnarray*}
     &&
  \frac{1}{t_n} \sum_{t=0}^{t_n-1}
    \EXP \left\{
    < \left( \nabla_\bw F \right)(\bw^{(t)},\vartheta^{(t)}) - \left( \nabla_\bw F_t \right)(\bw^{(t)},\vartheta^{(t)}),
    \bw^{(t)}-\bw^*> 
    {}\cdot{} 1_{E_n} 
    \right\}
\\ 
&&
\leq
4 \cdot
\EXP \Bigg\{
\sup_{{(\bw,\vartheta) \in}\W, \atop k \in \{1, \dots, K_n\}} \left(
  \frac{1}{n}
  \sum_{i=1}^n
  \epsilon_i \cdot 
  \frac{1}{1+\exp(Y_i \cdot f_{ (\bw,\vartheta) }(X_i))}
  \cdot Y_i \cdot T_{\beta_n} f_{{\vartheta}_k}(X_i)
  \right) {}\cdot{} 1_{E_n}\Bigg\}
\\
&&
\quad
+
\frac{n \cdot \left(
6 \cdot K_n \cdot  \beta_n^2 
+
4 \cdot (\beta_n +1) \cdot C_n \cdot
\sup_{{(\bw,\vartheta) \in}\W, y \in \{-1,1\}, x \in [0,1]^{d_1 \times d_2}} \| \nabla_{\vartheta}
   \varphi( y \cdot f_{(\bw,\vartheta)}(x)) \|_\infty
\right)
}{t_n}
.
   \end{eqnarray*}
   We have
   \[
   \varphi^\prime(z)
   =
   \frac{1}{1+\exp(-z)} \cdot (- \exp(-z))
   =
   \frac{-1}{1+\exp(z)},
   \]
   which implies
   \[
   \frac{\partial}{\partial w_j}
   \varphi \left(
Y \cdot \sum_{k=1}^{K_n} w_k \cdot T_{\beta_n} f_{\vartheta_k}(X)
\right)
=
\frac{-Y \cdot T_{\beta_n} f_{\vartheta_j}(X) }{
  1+ \exp \left(
Y \cdot \sum_{k=1}^{K_n} w_k \cdot T_{\beta_n} f_{\vartheta_k}(X)
  \right)
}
   \]
   and
   \begin{eqnarray*}
     &&
    \frac{\partial}{\partial w_j}
    \EXP \left\{
   \varphi \left(
Y \cdot \sum_{k=1}^{K_n} w_k \cdot T_{\beta_n} f_{\vartheta_k}(X)
\right)
\right\}
\\
&&
=
\lim_{h \rightarrow 0}
\EXP \left\{
\frac{
\varphi \left(
Y \cdot \sum_{k=1}^{K_n} (w_k+h \cdot I_{\{k=j\}}) \cdot T_{\beta_n} f_{\vartheta_k}(X)
\right)
-
   \varphi \left(
Y \cdot \sum_{k=1}^{K_n} w_k \cdot T_{\beta_n} f_{\vartheta_k}(X)
\right)
}{
h
}
\right\}
\\
&&
=
\EXP \left\{
\lim_{h \rightarrow 0}
\frac{
\varphi \left(
Y \cdot \sum_{k=1}^{K_n} (w_k+h \cdot I_{\{k=j\}}) \cdot T_{\beta_n} f_{\vartheta_k}(X)
\right)
-
   \varphi \left(
Y \cdot \sum_{k=1}^{K_n} w_k \cdot T_{\beta_n} f_{\vartheta_k}(X)
\right)
}{
h
}
\right\}
\\
&&
=
\EXP \left\{
\frac{- Y \cdot T_{\beta_n} f_{\vartheta_j}(X) }{
  1+ \exp \left(
Y \cdot \sum_{k=1}^{K_n} w_k \cdot T_{\beta_n} f_{\vartheta_k}(X)
  \right)
}
\right\},
     \end{eqnarray*}
   where we have used
   \begin{eqnarray*}
     &&
     \left|
\frac{
\varphi \left(
Y \cdot \sum_{k=1}^{K_n} (w_k+h \cdot I_{\{k=j\}}) \cdot T_{\beta_n} f_{\vartheta_k}(X)
\right)
-
   \varphi \left(
Y \cdot \sum_{k=1}^{K_n} w_k \cdot T_{\beta_n} f_{\vartheta_k}(X)
\right)
}{
h
}
\right|
\\
&&
=
\beta_n \cdot | \varphi^\prime (\xi)| \leq \beta_n
   \end{eqnarray*}
   and the dominated convergence theorem in order to interchange
   limites and expectations above.
   
   Consequently,
   \begin{eqnarray*}
&&
 \frac{1}{t_n} \sum_{t=0}^{t_n-1}
    \EXP \left\{
    < \left( \nabla_\bw F \right)(\bw^{(t)},\vartheta^{(t)}) - \left( \nabla_\bw F_t \right)(\bw^{(t)},\vartheta^{(t)}),
    \bw^{(t)}-\bw^*>  
    {}\cdot{} 1_{E_n}\right\}
 \\
&&
=
\EXP \Bigg\{
\frac{1}{t_n} \sum_{t=0}^{t_n-1}
\sum_{k=1}^{K_n}
\Bigg(
\EXP \left\{
\frac{-
Y \cdot T_{\beta_n} f_{\vartheta_k^{(t)}}(X)
}{
1 + \exp \left(
Y \cdot f_{(\bw^{(t)},\vartheta^{(t)})}(X)
\right)
}
\bigg|
\D_n, \vartheta^{(0)}
\right\}
\\
&&
\hspace*{4cm}
-
\frac{-
Y_{j_t} \cdot T_{\beta_n} f_{\vartheta_k^{(t)}}(X_{j_t})
}{
1 + \exp \left(
Y_{j_t} \cdot f_{(\bw^{(t)},\vartheta^{(t)})}(X_{j_t})
\right)
}
\Bigg)
\cdot
(w_k^{(t)}-w_k^*)
 {}\cdot{} 1_{E_n}
\Bigg\}
\\
&&
=
\frac{1}{t_n/n} \cdot \sum_{s=1}^{t_n/n}
\EXP \Bigg\{
\frac{1}{n} \sum_{t=(s-1) \cdot n}^{s \cdot n-1}
\sum_{k=1}^{K_n}
\Bigg(
\EXP \left\{
\frac{-
Y \cdot T_{\beta_n} f_{\vartheta_k^{(t)}}(X)
}{
1 + \exp \left(
Y \cdot f_{(\bw^{(t)},\vartheta^{(t)})}(X)
\right)
}
\bigg|
\D_n, \vartheta^{(0)}
\right\}
\\
&&
\hspace*{4cm}
-
\frac{-
Y_{j_t} \cdot T_{\beta_n} f_{\vartheta_k^{(t)}}(X_{j_t})
}{
1 + \exp \left(
Y_{j_t} \cdot f_{(\bw^{(t)},\vartheta^{(t)})}(X_{j_t})
\right)
}
\Bigg)
\cdot
(w_k^{(t)}-w_k^*)
 {}\cdot{} 1_{E_n}
\Bigg\}.
\end{eqnarray*}
During $n$ gradient descent steps the parameter $(\bw,\vartheta)$
changes in supremum norm at most by
\begin{eqnarray*}
&&
n \cdot \lambda_n \cdot \max \Bigg\{
 \sup_{{(\bw,\vartheta) \in}\W, y \in \{-1,1\}, x \in [0,1]^{d_1 \times d_2}} \| \nabla_{\bw}
   \varphi( y \cdot f_{(\bw,\vartheta)}(x)) \|_\infty,\\
&&
\hspace*{4cm}
\sup_{{(\bw,\vartheta) \in}\W, y \in \{-1,1\}, x \in [0,1]^{d_1 \times d_2}} \| \nabla_{\vartheta}
   \varphi( y \cdot f_{(\bw,\vartheta)}(x)) \|_\infty
\Bigg\}
\\
&&
\leq
n \cdot \lambda_n \cdot \left(
 \beta_n 
+
\sup_{{(\bw,\vartheta) \in}\W, y \in \{-1,1\}, x \in [0,1]^{d_1 \times d_2}} \| \nabla_{\vartheta}
   \varphi( y \cdot f_{(\bw,\vartheta)}(x)) \|_\infty
\right)
\\
&&
=
\frac{
n \cdot \left(
 \beta_n 
+
\sup_{{(\bw,\vartheta) \in}\W, y \in \{-1,1\}, x \in [0,1]^{d_1 \times d_2}} \| \nabla_{\vartheta}
   \varphi( y \cdot f_{(\bw,\vartheta)}(x)) \|_\infty
\right)
}{t_n}.
\end{eqnarray*}
Using
\begin{eqnarray*}
  &&
\sum_{k=1}^{K_n}
\Bigg(
\EXP \left\{
\frac{-
Y \cdot T_{\beta_n} f_{\vartheta_k^{(t)}}(X)
}{
1 + \exp \left(
Y \cdot f_{(\bw^{(t)},\vartheta^{(t)})}(X)
\right)
}
\bigg|
\D_n, \vartheta^{(0)}
\right\}
\\
&&
\hspace*{3cm}
-\EXP \left\{
\frac{-
Y \cdot T_{\beta_n} f_{\vartheta_k^{(s \cdot n -1)}}(X)
}{
1 + \exp \left(
Y \cdot f_{(\bw^{(s \cdot n -1)},\vartheta^{(s \cdot n -1)})}(X)
\right)
}
\bigg|
\D_n, \vartheta^{(0)}
\right\}
\\
&&
\hspace*{2cm}
-
\frac{-
Y_{j_t} \cdot T_{\beta_n} f_{\vartheta_k^{(t)}}(X_{j_t})
}{
1 + \exp \left(
Y_{j_t} \cdot f_{(\bw^{(t)},\vartheta^{(t)})}(X_{j_t})
\right)
}
\\
&&
\hspace*{3cm}
+
\frac{-
Y_{j_t} \cdot T_{\beta_n} f_{\vartheta_k^{(s \cdot n -1)}}(X_{j_t})
}{
1 + \exp \left(
Y_{j_t} \cdot f_{(\bw^{(s \cdot n -1)},\vartheta^{(s \cdot n -1)})}(X_{j_t})
\right)
}
\Bigg)
\cdot
(w_k^{(t)}-w_k^*)
\\
  &&
  \leq
\sum_{k=1}^{K_n}
\Bigg(
\EXP \left\{
\left|
\frac{-
Y \cdot T_{\beta_n} f_{\vartheta_k^{(t)}}(X)
+
Y \cdot T_{\beta_n} f_{\vartheta_k^{(s \cdot n -1)}}(X)
}{
1 + \exp \left(
Y \cdot f_{(\bw^{(t)},\vartheta^{(t)})}(X)
\right)
}
\right|
\right.
\\
&&
\hspace*{1.5cm}
+ 
\left|
Y \cdot T_{\beta_n} f_{\vartheta_k^{(s \cdot n -1)}}(X)
\right|
\cdot
\left|
\frac{1
}{
1 + \exp \left(
Y \cdot f_{(\bw^{(t)},\vartheta^{(t)})}(X)
\right)
}\right.
\\
&&
\hspace{5.5cm}
\left.
\left.
-
\frac{1
}{
1 + \exp \left(
Y \cdot f_{(\bw^{(s \cdot n -1)},\vartheta^{(s \cdot n -1)})}(X)
\right)
}
\right|
\bigg|
\D_n, \vartheta^{(0)}
\right\}
\\
&&
\hspace*{1cm}
+
\left|
\frac{-
Y_{j_t} \cdot T_{\beta_n} f_{\vartheta_k^{(s \cdot n -1)}}(X_{j_t})
+
Y_{j_t} \cdot T_{\beta_n} f_{\vartheta_k^{(t)}}(X_{j_t})}{
1 + \exp \left(
Y_{j_t} \cdot f_{(\bw^{(s \cdot n -1)},\vartheta^{(s \cdot n -1)})}(X_{j_t})
\right)
}
\right|
\\
&&
\hspace*{1cm}
+ 
\left|
Y_{j_t} \cdot T_{\beta_n} f_{\vartheta_k^{(t)}}(X_{j_t})
\right|
\cdot 
\left|
\frac{1
}{
1 + \exp \left(
Y_{j_t} \cdot f_{(\bw^{(t)},\vartheta^{(t)})}(X_{j_t})
\right)
}
\right.
\\
&&
\hspace*{5cm}
\left.
-
\frac{1
}{
1 + \exp \left(
Y_{j_t} \cdot f_{(\bw^{(s \cdot n -1)},\vartheta^{(s \cdot n -1)})}(X_{j_t})
\right)
}
\right|
\Bigg)
\cdot
(w_k^{(t)}-w_k^*)
\\
&&
\leq
\sum_{k=1}^{K_n}
\Bigg(
\EXP \left\{
\left|
{
 T_{\beta_n} f_{\vartheta_k^{(s \cdot n -1)}}(X)
-
 T_{\beta_n} f_{\vartheta_k^{(t)}}(X)
}
\right|
+ 
\beta_n
\cdot
\left|
 f_{(\bw^{(t)},\vartheta^{(t)})}(X)
-
 f_{(\bw^{(s \cdot n -1)},\vartheta^{(s \cdot n -1)})}(X)
\right|
\bigg|
\D_n, \vartheta^{(0)}
\right\}
\\
&&
\hspace*{1.5cm}
+
\left|
{
 T_{\beta_n} f_{\vartheta_k^{(t)}}(X_{j_t})}
 -
 T_{\beta_n} f_{\vartheta_k^{(s \cdot n -1)}}(X_{j_t})
\right|
\\
&&
\hspace*{3cm}
+ 
\beta_n
\cdot 
\left|
 f_{(\bw^{(t)},\vartheta^{(t)})}(X_{j_t})
-
 f_{(\bw^{(s \cdot n -1)},\vartheta^{(s \cdot n -1)})}(X_{j_t})
\right|
\Bigg)
\cdot
|w_k^{(t)}-w_k^*|
\\
&&
\leq
\sum_{k=1}^{K_n}
\Bigg(
\EXP \left\{
\left|
{
  f_{\vartheta_k^{(s \cdot n -1)}}(X)
-
 f_{\vartheta_k^{(t)}}(X)
}
\right|
+ 
\beta_n
\cdot
\left|
 f_{(\bw^{(t)},\vartheta^{(t)})}(X)
-
 f_{(\bw^{(s \cdot n -1)},\vartheta^{(s \cdot n -1)})}(X)
\right|
\bigg|
\D_n, \vartheta^{(0)}
\right\}
\\
&&
\hspace*{2cm}
+
\left|
{
  f_{\vartheta_k^{(t)}}(X_{j_t})}
 -
f_{\vartheta_k^{(s \cdot n -1)}}(X_{j_t})
\right|
\\
&&
\hspace*{3cm}
+ 
\beta_n
\cdot 
\left|
 f_{(\bw^{(t)},\vartheta^{(t)})}(X_{j_t})
-
 f_{(\bw^{(s \cdot n -1)},\vartheta^{(s \cdot n -1)})}(X_{j_t})
\right|
\Bigg)
\cdot
|w_k^{(t)}-w_k^*|
\\
&&
\leq
2\cdot
\sum_{k=1}^{K_n}
\Bigg(
C_n \cdot \|\vartheta_k^{(t)}- \vartheta_k^{(s \cdot n -1)}\|_\infty
+
 \beta_n \cdot
  C_n
\cdot
\max_{j} \|\vartheta_j^{(t)}- \vartheta_j^{(s \cdot n -1)}\|_\infty +K_n\cdot\frac{\beta_n^2\cdot n}{t_n}
\Bigg)
\cdot
|w_k^{(t)}-w_k^*|
\\
&&
\leq 
4\cdot
C_n \cdot \|\vartheta_k^{(t)}- \vartheta_k^{(s \cdot n -1)}\|_\infty
+
 4\cdot\beta_n \cdot C_n
\cdot
\max_{j} \|\vartheta_j^{(t)}- \vartheta_j^{(s \cdot n -1)}\|_\infty
+4\cdot K_n\cdot\frac{\beta_n^2\cdot n}{t_n}\\
&&
\leq
\frac{
4 \cdot (\beta_n+1) \cdot C_n \cdot n \cdot
\sup_{{(\bw,\vartheta)\in} \W, y \in \{-1,1\}, x \in [0,1]^{d_1 \times d_2}} \| \nabla_{\vartheta}
\varphi( y \cdot f_{(\bw,\vartheta)}(x)) \|_\infty
+4\cdot K_n\cdot\beta_n^2\cdot n}{t_n},
\end{eqnarray*}
where the fourth inequality follows from 
\begin{eqnarray*}
  &&  \left|
 f_{(\bw^{(t)},\vartheta^{(t)})}(X)
-
 f_{(\bw^{(s \cdot n -1)},\vartheta^{(s \cdot n -1)})}(X)
\right|\\
&&=
\left|
\sum_{k=1}^{K_n}
w_k^{(t)}
T_{\beta_n}f_{\vartheta^{(t)}}(X)
-
w_k^{(s \cdot n -1)}T_{\beta_n} f_{\bw^{(s \cdot n -1)}}(X)
\right|\\
&&\leq 
\sum_{k=1}^{K_n}
\left|w_k^{(t)}
\right|\cdot 
\left|
T_{\beta_n}f_{\vartheta^{(t)}}(X)
-
T_{\beta_n} f_{\bw^{(s \cdot n -1)}}(X)
\right|
+
\left|
w_k^{(t)}
-
w_k^{(s \cdot n -1)}
\right|\cdot 
\left|
T_{\beta_n} f_{\bw^{(s \cdot n -1)}}(X)
\right|\\
&&\leq 
C_n\cdot \max_{j} \|\vartheta_j^{(t)}- \vartheta_j^{(s \cdot n -1)}\|_\infty
+ K_n\cdot \frac{n\cdot \beta_n}{t_n}\cdot \beta_n
\end{eqnarray*}

and
\begin{eqnarray*}
  &&
  \sum_{k=1}^{K_n}
\Bigg(
\EXP \left\{
\frac{-
Y \cdot T_{\beta_n} f_{\vartheta_k^{(s \cdot n-1)}}(X)
}{
1 + \exp \left(
Y \cdot f_{(\bw^{(s \cdot n-1)},\vartheta^{(s \cdot n-1)})}(X)
\right)
}
\bigg|
\D_n, \vartheta^{(0)}
\right\}
\\
&&
\hspace*{3cm}
-
\frac{-
Y_{j_t} \cdot T_{\beta_n} f_{\vartheta_k^{(s \cdot n-1)}}(X_{j_t})
}{
1 + \exp \left(
Y_{j_t} \cdot f_{(\bw^{(s \cdot n-1)},\vartheta^{(s \cdot n-1)})}(X_{j_t})
\right)
}
 \Bigg)
 \cdot
 (w_k^{(t)}-w_k^{(s \cdot n -1)})
 \\
 &&
 \leq
 2 \cdot\beta_n\cdot \sum_{k=1}^{K_n} |w_k^{(t)}-w_k^{(s \cdot n -1)}|
 \leq
 2 \cdot K_n \cdot n \cdot \frac{\beta_n^2}{t_n}
  \end{eqnarray*}
we can conclude
\begin{eqnarray*}
&&
\frac{1}{t_n/n} \cdot \sum_{s=1}^{t_n/n}
\EXP \Bigg\{
\frac{1}{n} \sum_{t=(s-1) \cdot n}^{s \cdot n-1}
\sum_{k=1}^{K_n}
\Bigg(
\EXP \left\{
\frac{-
Y \cdot T_{\beta_n} f_{\vartheta_k^{(t)}}(X)
}{
1 + \exp \left(
Y \cdot f_{(\bw^{(t)},\vartheta^{(t)})}(X)
\right)
}
\bigg|
\D_n, \vartheta^{(0)}
\right\}
\\
&&
\hspace*{4cm}
-
\frac{-
Y_{j_t} \cdot T_{\beta_n} f_{\vartheta_k^{(t)}}(X_{j_t})
}{
1 + \exp \left(
Y_{j_t} \cdot f_{(\bw^{(t)},\vartheta^{(t)})}(X_{j_t})
\right)
}
\Bigg)
\cdot
(w_k^{(t)}-w_k^*)
\cdot 1_{E_n}
\Bigg\}
\\
&&
\leq
\frac{1}{t_n/n} \cdot \sum_{s=1}^{t_n/n}
\EXP \Bigg\{
\frac{1}{n} \sum_{t=(s-1) \cdot n}^{s \cdot n-1}
\sum_{k=1}^{K_n}
\Bigg(
\EXP \left\{
\frac{-
Y \cdot T_{\beta_n} f_{\vartheta_k^{(s \cdot n-1)}}(X)
}{
1 + \exp \left(
Y \cdot f_{(\bw^{(s \cdot n-1)},\vartheta^{(s \cdot n-1)})}(X)
\right)
}
\bigg|
\D_n, \vartheta^{(0)}
\right\}
\\
&&
\hspace*{4cm}
-
\frac{-
Y_{j_t} \cdot T_{\beta_n} f_{\vartheta_k^{(s \cdot n-1)}}(X_{j_t})
}{
1 + \exp \left(
Y_{j_t} \cdot f_{(\bw^{(s \cdot n-1)},\vartheta^{(s \cdot n-1)})}(X_{j_t})
\right)
}
\Bigg)
\cdot
(w_k^{(s \cdot n-1)}-w_k^*)
\cdot 1_{E_n}
\Bigg\}
\\
&&
\quad
+ 
\frac{n \cdot \left(
6 \cdot K_n \cdot  \beta_n^2 
+
4 \cdot (\beta_n +1) \cdot C_n \cdot
\sup_{{(\bw,\vartheta) \in}\W, y \in \{-1,1\}, x \in [0,1]^{d_1 \times d_2}} \| \nabla_{\vartheta}
   \varphi( y \cdot f_{(\bw,\vartheta)}(x)) \|_\infty
\right)
}{t_n}.
\end{eqnarray*}
We continue by deriving an upper bound on the first term of the sum
of the right-hand side above. We have
\begin{eqnarray*}
&&
\frac{1}{t_n/n} \cdot \sum_{s=1}^{t_n/n}
\EXP \Bigg\{
\frac{1}{n} \sum_{t=(s-1) \cdot n}^{s \cdot n-1}
\sum_{k=1}^{K_n}
\Bigg(
\EXP \left\{
\frac{-
Y \cdot T_{\beta_n} f_{\vartheta_k^{(s \cdot n-1)}}(X)
}{
1 + \exp \left(
Y \cdot f_{(\bw^{(s \cdot n-1)},\vartheta^{(s \cdot n-1)})}(X)
\right)
}
\bigg|
\D_n, \vartheta^{(0)}
\right\}
\\
&&
\hspace*{4cm}
-
\frac{-
Y_{j_t} \cdot T_{\beta_n} f_{\vartheta_k^{(s \cdot n-1)}}(X_{j_t})
}{
1 + \exp \left(
Y_{j_t} \cdot f_{(\bw^{(s \cdot n-1)},\vartheta^{(s \cdot n-1)})}(X_{j_t})
\right)
}
\Bigg)
\cdot
(w_k^{(s \cdot n-1)}-w_k^*)
\cdot 1_{E_n}
\Bigg\}
\\
&&
\leq
\frac{1}{t_n/n} \cdot \sum_{s=1}^{t_n/n}
\EXP \Bigg\{ \sup_{{(\bw,\vartheta) \in}\W}
\frac{1}{n} \sum_{t=(s-1) \cdot n}^{s \cdot n-1}
\sum_{k=1}^{K_n}
\Bigg(
\EXP \left\{
\frac{-
Y \cdot T_{\beta_n} f_{\vartheta_k}(X)
}{
1 + \exp \left(
Y \cdot f_{(\bw,\vartheta)}(X)
\right)
}
\right\}
\\
&&
\hspace*{4cm}
-
\frac{-
Y_{j_t} \cdot T_{\beta_n} f_{\vartheta_k}(X_{j_t})
}{
1 + \exp \left(
Y_{j_t} \cdot f_{(\bw,\vartheta)}(X_{j_t})
\right)
}
\Bigg)
\cdot
(w_k-w_k^*)
\cdot 1_{E_n}
\Bigg\}
\\
&&
=
\EXP \Bigg\{ \sup_{{(\bw,\vartheta) \in}\W}
\frac{1}{n} \sum_{j=1}^{n}
\sum_{k=1}^{K_n}
\Bigg(
\EXP \left\{
\frac{-
Y \cdot T_{\beta_n} f_{\vartheta_k}(X)
}{
1 + \exp \left(
Y \cdot f_{(\bw,\vartheta)}(X)
\right)
}
\right\}
\\
&&
\hspace*{5cm}
-
\frac{-
Y_{j} \cdot T_{\beta_n} f_{\vartheta_k}(X_{j})
}{
1 + \exp \left(
Y_j \cdot f_{(\bw,\vartheta)}(X_{j})
\right)
}
\Bigg)
\cdot
(w_k-w_k^*)
\cdot 1_{E_n}
\Bigg\}
\\
&&
=
\EXP \Bigg\{ \sup_{{(\bw,\vartheta) \in}\W}
\frac{1}{n} \sum_{j=1}^{n}
\sum_{k=1}^{K_n}
\Bigg(
\EXP \left\{\frac{-
Y_{j}^\prime \cdot T_{\beta_n} f_{\vartheta_k}(X_{j}^\prime)
}{
1 + \exp \left(
Y_{j}^\prime \cdot f_{(\bw,\vartheta)}(X_{j}^\prime)
\right)
}
\bigg| \D_n, \vartheta^{(0)} \right\}
\\
&&
\hspace*{5cm}
-
\frac{-
Y_{j} \cdot T_{\beta_n} f_{\vartheta_k}(X_{j})
}{
1 + \exp \left(
Y_j \cdot f_{(\bw,\vartheta)}(X_{j})
\right)
}
\Bigg)
\cdot
(w_k-w_k^*)
\cdot 1_{E_n}
\Bigg\}
\\
&&
\leq
\EXP \Bigg\{ \sup_{{(\bw,\vartheta) \in}\W}
\frac{1}{n} \sum_{j=1}^{n}
\sum_{k=1}^{K_n}
\Bigg(
\frac{-
Y_{j}^\prime \cdot T_{\beta_n} f_{\vartheta_k}(X_{j}^\prime)
}{
1 + \exp \left(
Y_{j}^\prime \cdot f_{(\bw,\vartheta)}(X_{j}^\prime)
\right)
}
\\
&&
\hspace*{4cm}
-
\frac{-
Y_{j} \cdot T_{\beta_n} f_{\vartheta_k}(X_{j})
}{
1 + \exp \left(
Y_j \cdot f_{(\bw,\vartheta)}(X_{j})
\right)
}
\Bigg)
\cdot
(w_k-w_k^*)
\cdot 1_{E_n}
\Bigg\}
\\
&&
=
\EXP \Bigg\{ \sup_{{(\bw,\vartheta) \in}\W}
\frac{1}{n} \sum_{j=1}^{n}
\sum_{k=1}^{K_n}
\epsilon_j \cdot
\Bigg(
\frac{-
Y_{j}^\prime \cdot T_{\beta_n} f_{\vartheta_k}(X_{j}^\prime)
}{
1 + \exp \left(
Y_{j}^\prime \cdot f_{(\bw,\vartheta)}(X_{j}^\prime)
\right)
}
\\
&&
\hspace*{4cm}
-
\frac{-
Y_{j} \cdot T_{\beta_n} f_{\vartheta_k}(X_{j})
}{
1 + \exp \left(
Y_j \cdot f_{(\bw,\vartheta)}(X_{j})
\right)
}
\Bigg)
\cdot
(w_k-w_k^*)
\cdot 1_{E_n}
\Bigg\}
\\
&&
\leq
\EXP \Bigg\{ \sup_{{(\bw,\vartheta) \in}\W}
\frac{1}{n} \sum_{j=1}^{n}
\sum_{k=1}^{K_n}
\epsilon_j \cdot
\Bigg(
\frac{-
Y_{j}^\prime \cdot T_{\beta_n} f_{\vartheta_k}(X_{j}^\prime)
}{
1 + \exp \left(
Y_{j}^\prime \cdot f_{(\bw,\vartheta)}(X_{j}^\prime)
\right)
}
\Bigg)
\cdot
(w_k-w_k^*)
\cdot 1_{E_n}
\Bigg\}
\\
&&
\quad
+
\EXP \Bigg\{ \sup_{{(\bw,\vartheta) \in}\W}
\frac{1}{n} \sum_{j=1}^{n}
\sum_{k=1}^{K_n}
\epsilon_j \cdot
\Bigg(
\frac{
Y_{j} \cdot T_{\beta_n} f_{\vartheta_k}(X_{j})
}{
1 + \exp \left(
Y_j \cdot f_{(\bw,\vartheta)}(X_{j})
\right)
}
\Bigg)
\cdot
(w_k-w_k^*)
\cdot 1_{E_n}
\Bigg\}
\\
&&
=
2 \cdot
\EXP \Bigg\{ \sup_{{(\bw,\vartheta) \in}\W}
\frac{1}{n} \sum_{j=1}^{n}
\sum_{k=1}^{K_n}
\epsilon_j \cdot
\Bigg(
\frac{
Y_{j} \cdot T_{\beta_n} f_{\vartheta_k}(X_{j})
}{
1 + \exp \left(
Y_j \cdot f_{(\bw,\vartheta)}(X_{j})
\right)
}
\Bigg)
\cdot
(w_k-w_k^*)
\cdot 1_{E_n}
\Bigg\}
\\
&&
\leq
2 \cdot
\EXP \Bigg\{ \sup_{{(\bw,\vartheta) \in}\W}
\frac{1}{n} \sum_{j=1}^{n}
\sum_{k=1}^{K_n}
\epsilon_j \cdot \frac{
Y_{j} \cdot T_{\beta_n} f_{\vartheta_k}(X_{j})
}{
1 + \exp \left(
Y_j \cdot f_{(\bw,\vartheta)}(X_{j})
\right)
}
\cdot
(w_k-w_k^*)_+
\cdot 1_{E_n}
\Bigg\}
\\
&&
\quad
+
2 \cdot
\EXP \Bigg\{ \sup_{{(\bw,\vartheta) \in}\W}
\frac{1}{n} \sum_{j=1}^{n}
\sum_{k=1}^{K_n}
(-\epsilon_j) \cdot \frac{
Y_{j} \cdot T_{\beta_n} f_{\vartheta_k}(X_{j})
}{
1 + \exp \left(
Y_j \cdot f_{(\bw,\vartheta)}(X_{j})
\right)
}
\cdot
(w_k^*-w_k)_+
\cdot 1_{E_n}
\Bigg\}
\\
&&
\leq
2 \cdot
\EXP \Bigg\{
\sup_{{(\bw,\vartheta) \in}\W}
\sum_{k=1}^{K_n}
\sup_{(\bar{\bw},\bar{\vartheta}) \in \W,\atop \bar{k} \in \{1, \dots, K_n\}}
\frac{1}{n} \sum_{j=1}^{n}
\epsilon_j \cdot \frac{
Y_{j} \cdot T_{\beta_n} f_{\bar{\vartheta}_{\bar{k}}}(X_{j})
}{
1 + \exp \left(
Y_j \cdot f_{(\bar{\bw},\bar{\vartheta})}(X_{j})
\right)
}
\cdot
(w_k-w_k^*)_+
\cdot 1_{E_n}
\Bigg\}
\\
&&
\quad
+
2 \cdot
\EXP \Bigg\{
\sup_{{(\bw,\vartheta) \in}\W}
\sum_{k=1}^{K_n}
\sup_{(\bar{\bw},\bar{\vartheta}) \in \W,\atop \bar{k} \in \{1, \dots, K_n\}}
\frac{1}{n} \sum_{j=1}^{n}
(-\epsilon_j) \cdot \frac{
Y_{j} \cdot T_{\beta_n} f_{\bar{\vartheta}_{\bar{k}}}(X_{j})
}{
1 + \exp \left(
Y_j \cdot f_{(\bar{\bw},\bar{\vartheta})}(X_{j})
\right)
}
\cdot (w_k^*-w_k)_+
\cdot 1_{E_n}
\Bigg\}
\\
&&
\leq
2 \cdot
\EXP \Bigg\{
\sup_{{(\bw,\vartheta) \in}\W, \atop k \in \{1, \dots, K_n\}}
\frac{1}{n} \sum_{j=1}^{n}
\epsilon_j \cdot \frac{
Y_{j} \cdot T_{\beta_n} f_{\vartheta_k}(X_{j})
}{
1 + \exp \left(
Y_j \cdot f_{(\bw,\vartheta)}(X_{j})
\right)
}
\cdot
\sup_{{(\bw,\vartheta) \in}\W}
\sum_{k=1}^{K_n}
(w_k-w_k^*)_+
\cdot 1_{E_n}
\Bigg\}
\\
&&
\quad
+
2 \cdot
\EXP \Bigg\{ \sup_{{(\bw,\vartheta) \in}\W,\atop k \in \{1, \dots, K_n\}}
\frac{1}{n} \sum_{j=1}^{n}
(-\epsilon_j) \cdot \frac{
Y_{j} \cdot T_{\beta_n} f_{\vartheta_k}(X_{j})
}{
1 + \exp \left(
Y_j \cdot f_{(\bw,\vartheta)}(X_{j})
\right)
}
\cdot
\sup_{{(\bw,\vartheta) \in}\W}
\sum_{k=1}^{K_n}
(w_k^*-w_k)_+
\cdot 1_{E_n}
\Bigg\}
\\
&&
\leq
4 \cdot
\EXP \Bigg\{ \sup_{{(\bw,\vartheta) \in}\W,\atop k \in \{1 , \dots, K_n\}}
\frac{1}{n} \sum_{j=1}^{n}
\epsilon_j \cdot \frac{
Y_{j} \cdot T_{\beta_n} f_{\vartheta_k}(X_{j})
}{
1 + \exp \left(
Y_j \cdot f_{(\bw,\vartheta)}(X_{j})
\right)
}
\Bigg\}.
\end{eqnarray*}
Where in the second to last inequality we have used that 
$$\sup_{(\bar{\bw}, \bar{\vartheta}) \in \W, \bar{k} \in \{1, \dots, K_n\}}
\frac{1}{n} \sum_{j=1}^{n}
(\pm\epsilon_j) \cdot \frac{
Y_{j} \cdot T_{\beta_n} f_{\bar{\vartheta}_{\bar{k}}}(X_{j})
}{
1 + \exp \left(
Y_j \cdot f_{(\bar{\bw},\bar{\vartheta})}(X_{j})
\right)
}
\geq 0$$
since there exists $(\mathbf{\tilde w}, \tilde{\vartheta})\in \W$  such that $f_{\tilde{\vartheta}}(X_{j})=0$ for $j\in \{1,\ldots, n\}$.

    In the {\it ninth step of the proof} we derive an upper bound on
\begin{eqnarray*}
  &&
\EXP \Bigg\{
\sup_{{(\bw,\vartheta) \in}\W, \atop k \in \{1, \dots, K_n\}} \left(
  \frac{1}{n}
  \sum_{i=1}^n
  \epsilon_i \cdot 
  \frac{1}{1+\exp(Y_i \cdot f_{ (\bw,\vartheta) }(X_i))}
  \cdot Y_i \cdot T_{\beta_n} f_{{\vartheta}_k}(X_i)
  \right)\Bigg\}.
  \end{eqnarray*}
   To do this we use a contraction style argument. Because of the independence
   of the random variables we can compute the expectation by first computing
   the expectation with respect to $\epsilon_1$ and then by computing
   the expectation with respect to all other random variables.
   This yields that the last term above is equal to
   \begin{eqnarray*}
     &&
      \frac{1}{2} \cdot
\EXP \Bigg\{
\sup_{{(\bw,\vartheta) \in}\W, \atop k \in \{1, \dots, K_n\}} \frac{1}{n}\Bigg(
  \sum_{i=2}^n
  \epsilon_i \cdot 
  \frac{1}{1+\exp(Y_i \cdot f_{ (\bw,\vartheta) }(X_i))}
  \cdot Y_i \cdot T_{\beta_n} f_{{\vartheta}_k}(X_i)
  \\
  &&
  \hspace*{6cm}
  +
   \frac{1}{1+\exp(Y_1 \cdot f_{ (\bw,\vartheta) }(X_1))}
  \cdot Y_1 \cdot T_{\beta_n} f_{{\vartheta}_k}(X_1)
  \Bigg) \Bigg\}
  \\
  &&
  \quad
  +
   \frac{1}{2} \cdot
\EXP \Bigg\{
\sup_{{(\bw,\vartheta) \in}\W, \atop k \in \{1, \dots, K_n\}}\frac{1}{n} \Bigg(
  \sum_{i=2}^n
  \epsilon_i \cdot 
  \frac{1}{1+\exp(Y_i \cdot f_{ (\bw,\vartheta) }(X_i))}
  \cdot Y_i \cdot T_{\beta_n} f_{{\vartheta}_k}(X_i)
  \\
  &&
  \hspace*{6cm}
  -
   \frac{1}{1+\exp(Y_1 \cdot f_{ (\bw,\vartheta) }(X_1))}
  \cdot Y_1 \cdot T_{\beta_n} f_{{\vartheta}_k}(X_1)
  \Bigg) 
  \Bigg\}
  \\
  &&
  =
  \frac{1}{2} \cdot \EXP \Bigg\{
\sup_{(\bw,\vartheta) , (\bar{\bw},\bar{\vartheta}) \in \W, \atop k, \bar{k} \in \{1, \dots, K_n\}}  \frac{1}{n}\Bigg(
  \sum_{i=2}^n
  \epsilon_i \cdot 
  \frac{1}{1+\exp(Y_i \cdot f_{ (\bw,\vartheta) }(X_i))}
  \cdot Y_i \cdot T_{\beta_n} f_{{\vartheta}_k}(X_i)
  \\
  &&
  \hspace*{4cm}
  +
   \frac{1}{n}
  \sum_{i=2}^n
  \epsilon_i \cdot 
  \frac{1}{1+\exp(Y_i \cdot f_{ (\bar{\bw},\bar{\vartheta}) }(X_i))}
  \cdot Y_i \cdot T_{\beta_n} f_{\bar{\vartheta}_{\bar{k}}}(X_i)
 \\
  &&
  \hspace*{4cm}
  +
  \frac{1}{1+\exp(Y_1 \cdot f_{ (\bw,\vartheta) }(X_1))}
  \cdot Y_1 \cdot T_{\beta_n} f_{{\vartheta}_k}(X_1)
 \\
  &&
  \hspace*{4cm}
  -
    \frac{1}{1+\exp(Y_1 \cdot f_{ (\bar{\bw},\bar{\vartheta}) }(X_1))}
    \cdot Y_1 \cdot T_{\beta_n} f_{\bar{\vartheta}_{\bar{k}}}(X_1)
    \Bigg)
    \Bigg\}
    \\
    &&
  \leq
  \frac{1}{2} \cdot \EXP \Bigg\{
\sup_{(\bw,\vartheta) , (\bar{\bw},\bar{\vartheta}) \in \W, \atop k, \bar{k} \in \{1, \dots, K_n\}}  \frac{1}{n}\Bigg(
  \sum_{i=2}^n
  \epsilon_i \cdot 
  \frac{1}{1+\exp(Y_i \cdot f_{ (\bw,\vartheta) }(X_i))}
  \cdot Y_i \cdot T_{\beta_n} f_{{\vartheta}_k}(X_i)
  \\
  &&
  \hspace*{4cm}
  +
  \sum_{i=2}^n
  \epsilon_i \cdot 
  \frac{1}{1+\exp(Y_i \cdot f_{ (\bar{\bw},\bar{\vartheta}) }(X_i))}
  \cdot Y_i \cdot T_{\beta_n} f_{\bar{\vartheta}_{\bar{k}}}(X_i)
 \\
  &&
  \hspace*{2cm}
  +
  \beta_n \cdot |  f_{ (\bw,\vartheta) }(X_1) -  f_{ (\bar{\bw},\bar{\vartheta}) }(X_1)|
  +
  | T_{\beta_n} f_{{\vartheta}_k}(X_1) -  T_{\beta_n} f_{ \bar{\vartheta}_{\bar{k}} }(X_1)|
  \Bigg) 
  \Bigg\}
  .
   \end{eqnarray*}
   The sum inside the supremum above does not change its value if
$((\bw, \vartheta), k)$
 is interchanged with
   $((\bar{\bw},\bar{\vartheta}), \bar{k})$. Consequently we can assume without
     loss of generality that
     $ f_{ (\bw,\vartheta) }(X_1) -  f_{ (\bar{\bw},\bar{\vartheta})
     }(X_1)$  is positive or that it is negative.
     Set $\bar{\epsilon}_1=\bar{\epsilon}_1(\epsilon_1,X_1,Y_1, \dots,X_n,Y_n)=\epsilon_1$
     if the functions $ f_{ (\bw,\vartheta) }(X_1) -  f_{ (\bar{\bw},\bar{\vartheta}) }(X_1)$ and
     $ T_{\beta_n} f_{{\vartheta}_k}(X_1) -  T_{\beta_n} f_{ \bar{\vartheta}_{\bar{k}} }(X_1)$ which ''attain'' the above supremum have the same sign, and set it equal
     to $-\epsilon_1$ otherwise. Then the right-hand side
     above is equal to 
     \begin{eqnarray*}
       &&
  \frac{1}{2} \cdot \EXP \Bigg\{
\sup_{(\bw,\vartheta) , (\bar{\bw},\bar{\vartheta})\in \W, \atop k, \bar{k} \in \{1, \dots, K_n\}}  \frac{1}{n}
\Bigg(
  \sum_{i=2}^n
  \epsilon_i \cdot 
  \frac{1}{1+\exp(Y_i \cdot f_{ (\bw,\vartheta) }(X_i))}
  \cdot Y_i \cdot T_{\beta_n} f_{{\vartheta}_k}(X_i)
  \\
  &&
  \hspace*{4cm}
  +
  \sum_{i=2}^n
  \epsilon_i \cdot 
  \frac{1}{1+\exp(Y_i \cdot f_{ (\bar{\bw},\bar{\vartheta}) }(X_i))}
  \cdot Y_i \cdot T_{\beta_n} f_{\bar{\vartheta}_{\bar{k}}}(X_i)
 \\
  &&
  \hspace*{2cm}
  +
  \beta_n \cdot \epsilon_1 \cdot (  f_{ (\bw,\vartheta) }(X_1) -
  f_{ (\bar{\bw},\bar{\vartheta}) }(X_1))
  +
  \bar{\epsilon}_1 \cdot (T_{\beta_n} f_{{\vartheta}_k}(X_1) -  T_{\beta_n} f_{ \bar{\vartheta}_{\bar{k}} }(X_1))
  \Bigg)\Bigg\}\\
  &&
  \leq
\EXP \Bigg\{
\sup_{(\bw,\vartheta)  \in \W, \atop k \in \{1, \dots, K_n\}} \frac{1}{n}\Bigg(
  \sum_{i=2}^n
  \epsilon_i \cdot 
  \frac{1}{1+\exp(Y_i \cdot f_{ (\bw,\vartheta) }(X_i))}
  \cdot Y_i \cdot T_{\beta_n} f_{{\vartheta}_k}(X_i)
 \\
  &&
  \hspace*{4cm}
  +
  \beta_n \cdot \epsilon_1 \cdot   f_{ (\bw,\vartheta) }(X_1) 
  +
  \bar{\epsilon}_1 \cdot T_{\beta_n} f_{{\vartheta}_k}(X_1) 
  \Bigg) \Bigg\},
     \end{eqnarray*}
     where we have used that conditioned on $(X_1,Y_1)$, \dots, $(X_n,Y_n)$
     the random vector $(\epsilon_1, \bar{\epsilon}_1)$ has the
     same distribution as 
     the random vector $(-\epsilon_1, -\bar{\epsilon}_1)$.
     
   Arguing in the same way for $k=2, \dots, n$ we see that we can upper bound
   the term on the right-hand side above by
   \begin{eqnarray*}
     &&
  \EXP \Bigg\{
  \sup_{{(\bw,\vartheta) \in}\W, \atop k \in \{1, \dots, K_n\}}
\frac{1}{n}
\sum_{i=1}^n
  \Bigg(
       \beta_n \cdot \epsilon_i \cdot f_{ (\bw,\vartheta) }(X_i) 
  +
  \bar{\epsilon}_i \cdot T_{\beta_n} f_{{\vartheta}_k}(X_i) 
  \Bigg)
  \Bigg\}
  \\
  &&
  \leq
  \beta_n \cdot \EXP \Bigg\{
  \sup_{{(\bw,\vartheta) \in}\W}
\frac{1}{n}
\sum_{i=1}^n
   \epsilon_i \cdot  f_{ (\bw,\vartheta) }(X_i) 
   \Bigg\}
   \\
   &&
   \quad
   +
    \EXP \Bigg\{
  \sup_{{(\bw,\vartheta) \in}\W, \atop k \in \{1, \dots, K_n\}}
\frac{1}{n}
\sum_{i=1}^n
  \bar{\epsilon}_i \cdot T_{\beta_n} f_{{\vartheta}_k}(X_i)
  \Bigg\}
  \\
  &&
  \leq 
   \beta_n \cdot
  \sup_{x_1, \dots, x_n \in \R^{d_1 \times d_2}}
  \EXP \Bigg\{
  \sup_{{(\bw,\vartheta) \in}\W}
\frac{1}{n}
\sum_{i=1}^n
   \epsilon_i \cdot  f_{ (\bw,\vartheta) }(x_i) 
   \Bigg\}
   \\
   &&
   \quad
   +
    \sup_{x_1, \dots, x_n \in \R^{d_1 \times d_2}}
   \EXP \Bigg\{
  \sup_{{(\bw,\vartheta) \in}\W, \atop k \in \{1, \dots, K_n\}}
\frac{1}{n}
\sum_{i=1}^n
  {\epsilon}_i \cdot T_{\beta_n} f_{{\vartheta}_k}(x_i)
  \Bigg\}
\\
&&
=
   \beta_n \cdot
  \sup_{x_1, \dots, x_n \in \R^{d_1 \times d_2}}
  \EXP \Bigg\{
  \sup_{{(\bw,\vartheta) \in}\W}
\frac{1}{n}
\sum_{i=1}^n
   \epsilon_i \cdot  f_{ (\bw,\vartheta) }(x_i) 
   \Bigg\}
   \\
   &&
   \quad
   +
    \sup_{x_1, \dots, x_n \in \R^{d_1 \times d_2}}
   \EXP \Bigg\{
  \sup_{\vartheta \in \bar{\Theta}}
\frac{1}{n}
\sum_{i=1}^n
  {\epsilon}_i \cdot T_{\beta_n} f_{{\vartheta}}(x_i) 
  \Bigg\},
     \end{eqnarray*}
where the last equality follows from
\begin{equation}
\label{pth1eq2}
    \{ T_{\beta_n} f_{\btheta_k} \, : \, {(\bw,\vartheta) \in}\W, k \in \{1, \dots, K_n\}\}
    =
    \{ T_{\beta_n} f_{\btheta_1} \, : \, {(\bw,\vartheta) \in}\W \}. 
\end{equation}

    In the {\it tenth step of the proof} we show
for $x_1$, \dots, $x_n \in \R^{d_1 \times d_2}$ arbitrary
    \begin{eqnarray*}
      &&
      \EXP \Bigg\{
  \sup_{{(\bw,\vartheta) \in}\W}
\frac{1}{n}
\sum_{i=1}^n
   \epsilon_i \cdot  f_{ (\bw,\vartheta) }(x_i) 
   \Bigg\}
    \leq
       \EXP \Bigg\{
       \sup_{\btheta \in \bar{\Theta}}
       \left| \frac{1}{n} \sum_{i=1}^n \epsilon_i \cdot 
    (T_{\beta_n} f_{\btheta}(x_i))
    \right|
    \Bigg\}
.
    \end{eqnarray*}
We have
    \begin{eqnarray*}
      &&
            \EXP \Bigg\{
  \sup_{{(\bw,\vartheta) \in}\W}
\frac{1}{n}
\sum_{i=1}^n
   \epsilon_i \cdot  f_{ (\bw,\vartheta) }(x_i) 
   \Bigg\}
   \\
   &&
   =
       \EXP \Bigg\{
       \sup_{{(\bw,\vartheta) \in}\W}
    \frac{1}{n} \sum_{i=1}^n \epsilon_i \cdot 
    \sum_{j=1}^{K_n} w_j \cdot (T_{\beta_n} f_{\btheta_j}(x_i))
    \Bigg\}
    \\
    &&
   =
       \EXP \Bigg\{
       \sup_{{(\bw,\vartheta) \in}\W}
    \sum_{j=1}^{K_n} w_j \cdot  \frac{1}{n} \sum_{i=1}^n \epsilon_i \cdot 
   (T_{\beta_n} f_{\btheta_j}(x_i))
    \Bigg\}
    \\
   &&
   \leq
       \EXP \Bigg\{
       \sup_{{(\bw,\vartheta) \in}\W}
    \sum_{j=1}^{K_n} |w_j| \cdot  \left| \frac{1}{n} \sum_{i=1}^n \epsilon_i \cdot 
    (T_{\beta_n} f_{\btheta_j}(x_i))
    \right|
    \Bigg\}
    \\
   &&
   \leq
       \EXP \Bigg\{
       \sup_{{(\bw,\vartheta) \in}\W}
       \sum_{j=1}^{K_n} |w_j| \cdot
       \sup_{{(\bw,\vartheta) \in}\W, k \in \{1, \dots, K_n\}}
       \left| \frac{1}{n} \sum_{i=1}^n \epsilon_i \cdot 
    (T_{\beta_n} f_{\btheta_k}(x_i))
    \right|
    \Bigg\}
    \\
   &&
    \leq
    1 \cdot
       \EXP \Bigg\{
       \sup_{{(\bw,\vartheta) \in}\W, k \in \{1, \dots, K_n\}}
       \left| \frac{1}{n} \sum_{i=1}^n \epsilon_i \cdot 
    (T_{\beta_n} f_{\btheta_k}(x_i))
    \right|
    \Bigg\}
    \\
    &&
    =
       \EXP \Bigg\{
       \sup_{\btheta \in \bar{\Theta}}
       \left| \frac{1}{n} \sum_{i=1}^n \epsilon_i \cdot 
    (T_{\beta_n} f_{\btheta}(x_i))
    \right|
    \Bigg\}
,
      \end{eqnarray*}
    where the last equality followed from (\ref{pth1eq2}).

Summarizing the above results, the proof is complete.
        \hfill $\Box$

\subsection{Proof of Theorem \ref{th2}}
In the proof of Theorem \ref{th2} we will need the following auxiliary
results.

\subsubsection{Using the logistic risk for classification}

\begin{lemma}
  \label{le2}
  Let $\varphi$ be the logistic loss. Let $(X,Y), (X_1, Y_1), \dots, (X_n, Y_n)$ and $f^*$,
  $\mathcal{D}_n, f_n$ and $\hat{C}_n$ as in Sections \ref{se1} and \ref{se3}, and set
  \[
  f_{\varphi}^*=
  \arg \min_{f: [0,1]^{d_1 \times d_2} \rightarrow \bar{\R}}
  \EXP
  \left\{
\varphi( Y \cdot f(X) )
  \right\}
  .
  \]\\ 
{\bf a)} Then
\begin{eqnarray*}
&&
\PROB \left\{
Y \neq \hat{C}_n(X)| \D_n \right\}
-
\PROB \left\{
Y \neq f^*(X) \right\}
\\
&&
\leq
\frac{1}{\sqrt{2}}
\cdot
\left(
\EXP\left\{
\varphi( Y \cdot f_n(X))
| \D_n
\right\}
-
\EXP\left\{
\varphi( Y \cdot f_{\varphi}^*(X))
\right\}
\right)^{1/2}
\end{eqnarray*}
holds.\\
{\bf b)}
Then
\begin{eqnarray*}
&&
\PROB \left\{
Y \neq \hat{C}_n(X) | \D_n \right\}
-
\PROB \left\{
Y \neq f^*(X) \right\}
\\
&&
\leq
2
\cdot
\left(
\EXP\left\{
\varphi( Y \cdot f_n(X))
| \D_n
\right\}
-
\EXP\left\{
\varphi( Y \cdot f_{\varphi}^*(X))
\right\}
\right)
+
4 \cdot
\EXP\left\{
\varphi( Y \cdot f_{\varphi}^*(X))
\right\}.
\end{eqnarray*}
holds.
\\
{\bf c)} Assume that 
\[
\PROB \left\{ |f_{\varphi}^*(X)| > \tilde{F}_n\right\} \geq 1 - e^{-\tilde{F}_n}
\]
for a given sequence $\{\tilde{F}_n\}_{n \in \N}$ with $\tilde{F}_n \to \infty$.
Then
\[
\EXP\left\{
\varphi( Y \cdot f_{\varphi}^*(X))
\right\}
\leq c_{12} \cdot \tilde{F}_n \cdot e^{-\tilde{F}_n}
\]
holds.
  \end{lemma}

\noindent
    {\bf Proof.}
    {\bf a)} This result follows from Theorem 2.1 in Zhang (2004), where we choose $s=2$ and $c=2^{-1/2}$. \\
    {\bf b)} This result follows from Lemma 1 b)  in Kohler and Langer (2025). \\
    {\bf c)}  This result follows from Lemma 3 in Kim, Ohn and Kim (2019).
    \hfill $\Box$

\subsubsection{Lipschitz property of the networks}

\begin{lemma}
  \label{le3}
  Let $\F_n = \{ f_{\btheta} \, : \, \btheta \in \Theta \}$ be the class
  of deep convolutional neural networks introduced in Subsection
  \ref{se3sub2} (cf., (\ref{se3eq*})). Set
  \[
  M_{\max}=\max\{M_1, \dots, M_{L_n^{(1)}}\}
  \quad \mbox{and} \quad
  k_{max}=\max\{ k_1, \dots, k_{L_n^{(1)}} \}.
  \]
  Let $\btheta, \bar{\btheta} \in \Theta$ such that
  \begin{equation}
\label{weight-difference-bound}
    \norm{\btheta-\bar \btheta}_\infty\leq 1
  \end{equation}
  holds
and all weights in $f_{\btheta}$ are bounded in absolute value by $B_n \geq 0$.
Then 
\[
\|f_{\btheta}-f_{\bar{\btheta}}\|_{\infty,[0,1]^{d_1 \times d_2}}
\leq
  7 \cdot
  L_n^{(2)} \cdot L_n^{(1)} \cdot k_{max}^{L_n^{(1)}+1} \cdot (M_{max}^2+1)^{L_n^{(1)}} \cdot (B_n+1)^{L_n^{(1)}+3}
  \cdot   \|\btheta - \bar{\btheta}\|_\infty.
\]
\end{lemma}

\noindent
{\bf Proof.}
Let $o^{(l)}$, $g$ and $\bar{o}^{(l)}$, $\bar{g}$ be defined as in
Section \ref{se3sub2} using the weights in $\btheta$ and
$\bar{\btheta}$, resp., and set
\[
\|o^{(l)}(x)\|_\infty=\max_{(i,j),s_2}|o^{(l)}_{(i,j),s_2}(x)|.
\]
Then $|\sigma(z)| \leq |z|$ implies
\begin{eqnarray*}
|o_{(i,j),s_2}^{(r)}(x)|
&\leq& k_{max} \cdot (M_{max}^2+1) \cdot B_n \cdot
\max\{ \|o^{(r-1)}(x)\|_\infty,1\} \\
&\leq&
k_{max}^r \cdot (M_{max}^2+1)^r \cdot B_n^r. 
\end{eqnarray*}
Using this together with $|\sigma(z_1)-\sigma(z_2)| \leq |z_1-z_2|$
we conclude
\begin{eqnarray*}
  &&
  | o_{(i,j),s_2}^{(r)}(x) - \bar{o}_{(i,j),s_2}^{(r)}(x)|
  \\
  &&
  \leq
  \Bigg|
  \sum_{s_1=1}^{k_{r-1}}
  \sum_{t_1,t_2 \in \{1, \dots, M_r\}, \atop
    (i+t_1-1,j+t_2-1)}
  \left(w_{t_1,t_2,s_1,s_2}^{(r)} \cdot o_{i+t_1-1,j+t_2-1,s_1}^{(r-1)}(x)
  -
  \bar{w}_{t_1,t_2,s_1,s_2}^{(r)} \cdot \bar{o}_{i+t_1-1,j+t_2-1,s_1}^{(r-1)}(x)
  \right)
  \\
  &&
  \hspace*{4cm}
  + w_{s_2}^{(r)} - \bar{w}_{s_2}^{(r)}
  \Bigg|
  \\
  &&
  \leq
   \sum_{s_1=1}^{k_{r-1}}
  \sum_{t_1,t_2 \in \{1, \dots, M_r\}, \atop
    (i+t_1-1,j+t_2-1)}
  \left|w_{t_1,t_2,s_1,s_2}^{(r)} - \bar{w}_{t_1,t_2,s_1,s_2}^{(r)}\right| \cdot
  \left|o_{i+t_1-1,j+t_2-1,s_1}^{(r-1)}(x)
  \right|
   + |w_{s_2}^{(r)} - \bar{w}_{s_2}^{(r)}|
  \\
  &&
  \quad
  +
 \sum_{s_1=1}^{k_{r-1}}
  \sum_{t_1,t_2 \in \{1, \dots, M_r\}, \atop
    (i+t_1-1,j+t_2-1)}
  |\bar{w}_{t_1,t_2,s_1,s_2}^{(r)}| \cdot \left|o_{i+t_1-1,j+t_2-1,s_1}^{(r-1)}(x)
  - \bar{o}_{i+t_1-1,j+t_2-1,s_1}^{(r-1)}(x) \right|
  \\
  &&
  \leq
  k_{max} \cdot (M_{max}^2+1) \cdot \left(
  \|\btheta - \bar{\btheta}\|_\infty
  \cdot
k_{max}^{r-1} \cdot (M_{max}^2+1)^{r-1} \cdot B_n^{r-1} 
+
(B_n+1) \cdot \|o^{(r-1)}-\bar{o}^{(r-1)}\|_\infty
  \right)
  \\
  &&
  \leq
  r \cdot
  k_{max}^{r} \cdot (M_{max}^2+1)^{r} \cdot (B_n+1)^{r-1}
  \cdot   \|\btheta - \bar{\btheta}\|_\infty.
  \end{eqnarray*}
From this and
\[
|
\max\{a_1, \dots, a_n\} -
\max\{b_1, \dots, b_n\}
|
\leq
\max\{|a_1-b_1|, \dots, |a_n-b_n|\} 
\]
we conclude
\begin{eqnarray*}
  &&
  |f_{\bw,\bw_{bias},\bw_{out}}(x)
  -
  f_{\bar{\bw},\bar{\bw}_{bias},\bar{\bw}_{out}}(x)|
  \\
  &&
  \leq
  k_{max} \cdot  \|\btheta - \bar{\btheta}\|_\infty \cdot B_n
  \cdot
  k_{max}^{L_n^{(1)}} \cdot (M_{max}^2+1)^{L_n^{(1)}} \cdot B_n^{L_n^{(1)}}
  \\
  &&
  \quad
  +  k_{max} \cdot (B_n+1) \cdot
    L_n^{(1)} \cdot
  k_{max}^{L_n^{(1)}} \cdot (M_{max}^2+1)^{L_n^{(1)}} \cdot (B_n+1)^{L_n^{(1)}-1}
  \cdot   \|\btheta - \bar{\btheta}\|_\infty
  \\
  &&
  \leq
    (L_n^{(1)}+1) \cdot
  k_{max}^{L_n^{(1)}+1} \cdot (M_{max}^2+1)^{L_n^{(1)}} \cdot (B_n+1)^{L_n^{(1)}+1}
  \cdot   \|\btheta - \bar{\btheta}\|_\infty.
  \end{eqnarray*}
This implies
\begin{eqnarray*}
  &&
  |g(x)-\bar{g}(x)|
  \\
  &&
  \leq
  (L_n^{(2)}+1) \cdot  \|\btheta - \bar{\btheta}\|_\infty \cdot
  2 \cdot B_n \cdot k_{max} \cdot B_n \cdot
  k_{max}^{L_n^{(1)}} \cdot (M_{max}^2+1)^{L_n^{(1)}} \cdot (B_n+1)^{L_n^{(1)}}
  \\
  &&
  \quad
  +
  L_n^{(2)} \cdot (B_n+1)
  \cdot \Bigg(
  2 \cdot  \|\btheta - \bar{\btheta}\|_\infty
  \cdot
  k_{max} \cdot 
  B_n \cdot
  k_{max}^{L_n^{(1)}} \cdot (M_{max}^2+1)^{L_n^{(1)}} \cdot (B_n+1)^{L_n^{(1)}}
  \\
  &&
  \hspace*{1.5cm}
  +
  (B_n+1)
  \cdot
     (L_n^{(1)}+1) \cdot
  k_{max}^{L_n^{(1)}+1} \cdot (M_{max}^2+1)^{L_n^{(1)}} \cdot (B_n+1)^{L_n^{(1)}+1}
  \cdot   \|\btheta - \bar{\btheta}\|_\infty
  \Bigg)
  \\
  &&
  \leq
  7 \cdot
  L_n^{(2)} \cdot L_n^{(1)} \cdot k_{max}^{L_n^{(1)}+1} \cdot (M_{max}^2+1)^{L_n^{(1)}} \cdot (B_n+1)^{L_n^{(1)}+3}
  \cdot   \|\btheta - \bar{\btheta}\|_\infty.
  \end{eqnarray*}
\hfill $\Box$

\subsubsection{Approximation error}
In our next lemma
we present a bound on the error we make in case
that we replace the functions $g_{k,s}$ in a hierarchical model
by some approximations of them.
\begin{lemma}
  \label{le4}
  Let $d_1,d_2,t \in \N$ and $l\in\N$ with $2^{l}\leq \min\{d_1,d_2\}$.
  For $a\in\{1,\dots,t\}$, set $I=\{0, 1, \dots, 2^{l}-1\} \times \{0, 1, \dots, 2^{l}-1\}$ and define
  \[
m_a(\bx)=
\max_{
  (i,j) \in \Z^2 \, : \,
  (i,j)+I \subseteq \{1, \dots, d_1\} \times \{1, \dots, d_2\}
}
f_a\left(
x_{(i,j)+I}
\right)
\]
and
\[
\bar{m}_a(\bx)=
\max_{
  (i,j) \in \Z^2 \, : \,
  (i,j)+I\subseteq \{1, \dots, d_1\} \times \{1, \dots, d_2\}
}
\bar{f}_a\left(
x_{(i,j)+I}
\right),
\]
where $f_a$ and $\bar{f}_a$ satisfy
    \[
f_a=f_{l,1}^{(a)} \quad \mbox{and} \quad \bar{f}_a=\bar{f}_{l,1}^{(a)}
    \]
    for some
    $f_{k,s}^{(a)}, \bar{f}_{k,s}^{(a)} :\R^{\{1, \dots, 2^k\} \times \{1, \dots, 2^k\}} \rightarrow \R$ recursively defined by
    \begin{eqnarray*}
    f_{k,s}^{(a)}(\bx)&=&g_{k,s}^{(a)} \big(
    f_{k-1,4 \cdot (s-1) + 1}^{(a)}(x_{
\{1, \dots, 2^{k-1}\} \times \{1, \dots, 2^{k-1}\}
    })
    , \\
        &&
        \hspace*{1cm}
        f_{k-1,4 \cdot (s-1) + 2}^{(a)}(x_{
\{2^{k-1}+1, \dots, 2^k\} \times \{1, \dots, 2^{k-1}\}
        }), \\
        &&
        \hspace*{1cm}
        f_{k-1,4 \cdot (s-1) + 3}^{(a)}(x_{
\{1, \dots, 2^{k-1}\} \times \{2^{k-1}+1, \dots, 2^k\}
        }), \\
        &&
        \hspace*{1cm}
        f_{k-1,4 \cdot s}^{(a)}(x_{
\{2^{k-1}+1, \dots, 2^k\} \times \{2^{k-1}+1, \dots, 2^k\}
        })
    \big)
    \end{eqnarray*}
    and
    \begin{eqnarray*}
    \bar{f}_{k,s}^{(a)}(\bx)&=&\bar{g}_{k,s}^{(a)} \big(
    \bar{f}_{k-1,4 \cdot (s-1) + 1}^{(a)}(x_{
\{1, \dots, 2^{k-1}\} \times \{1, \dots, 2^{k-1}\}
    })
    , \\
        &&
        \hspace*{1cm}
        \bar{f}_{k-1,4 \cdot (s-1) + 2}^{(a)}(x_{
\{2^{k-1}+1, \dots, 2^k\} \times \{1, \dots, 2^{k-1}\}
        }),\\
        &&
        \hspace*{1cm}
        \bar{f}_{k-1,4 \cdot (s-1) + 3}^{(a)}(x_{
\{1, \dots, 2^{k-1}\} \times \{2^{k-1}+1, \dots, 2^k\}
        }), \\
        &&
        \hspace*{1cm}
        \bar{f}_{k-1,4 \cdot s}^{(a)}(x_{
\{2^{k-1}+1, \dots, 2^k\} \times \{2^{k-1}+1, \dots, 2^k\}
        })
    \big)
    \end{eqnarray*}
        for $k=2, \dots, l, s=1, \dots,4^{l-k}$,
        and
        \[
 f_{1,s}^{(a)}(x_{1,1},x_{1,2},x_{2,1},x_{2,2})= g_{1,s}^{(a)}(x_{1,1},x_{1,2},x_{2,1},x_{2,2})
 \]
and
    \[
 \bar{f}_{1,s}^{(a)}(x_{1,1},x_{1,2},x_{2,1},x_{2,2})= \bar{g}_{1,s}^{(a)}(x_{1,1},x_{1,2},x_{2,1},x_{2,2})
 \]
  for $s=1, \dots, 4^{l-1}$, where
  \[g_{k,s}^{(a)}:\R^4\rightarrow[0,1]\mbox{ and  }\bar{g}_{k,s}^{(a)}:\R^4\rightarrow\R\]   
are functions for $a\in\{1,\dots,t\}$, $k\in\{1, \dots, l\}$ and $s\in\{1, \dots,4^{l-k}\}$.
Furthermore,
let $g:\R^t\rightarrow[0,1]$ and $\bar{g}:\R^t\rightarrow\R$ be functions.
Assume that all restrictions
$g_{k,s}^{(a)}|_{[-2,2]^4}: [-2,2]^4 \rightarrow [0,1]$ and $g|_{[-2,2]^t}:[-2,2]^t\rightarrow[0,1]$ are Lipschitz continuous
with respect to the Euclidean distance with Lipschitz constant $C>0$ and for all $a\in\{1,\dots,t\}$, $k\in\{1,\dots,l\}$ and $s\in\{1,\dots,4^{l-k}\}$ we assume that
\begin{equation}
\left\|\bar{g}_{k,s}^{(a)}\right\|_{[-2,2]^4,\infty}\leq2.
\label{le3eq3}
\end{equation} Then for any
$\bx \in [0,1]^{\{1, \dots, d_1\} \times \{1, \dots, d_2\}}$ it holds:
\begin{align*}
&|g(m_1(\bx),\dots,m_t(\bx))-\bar{g}(\bar{m}_1(\bx),\dots,\bar{m}_t(\bx))| \\
&\leq\sqrt{t}\cdot(2C+1)^{l}
\\&~~~
\cdot\max_{a\in\{1,\dots,t\},j\in\{1,\dots,l\},s\in\{1,\dots,4^{l-j}\}}\left\{\|g_{j,s}^{(a)}-\bar{g}_{j,s}^{(a)}\|_{[-2,2]^4,\infty},\|g-\bar{g}\|_{[-2,2]^t,\infty}\right\}.
\end{align*}
\end{lemma}

\noindent
{\bf Proof.} See Lemma 1 in Kohler, Krzy\.zak and Walter (2022).
\hfill $\Box$

\begin{lemma}
  \label{le5}
    Let $d \in \N$,
  let $f:\Rd \rightarrow \R$ be $(p,C)$--smooth for some $p=q+s$,
  $q \in \N_0$  and $s \in (0,1]$, and $C>0$. Let $A \geq 1$
    and $M \in \N$ sufficiently large (independent of the size of $A$, but
     \begin{align*}
       M \geq 2 \ \mbox{and} \ M^{2p} \geq c_{13} \cdot \left(\max\left\{A, \|f\|_{C^q([-A,A]^d)}
       \right\}\right)^{4(q+1)},
     \end{align*}
     where
     \[
     \|f\|_{C^q([-A,A]^d)}
     =
     \max_{\alpha_1, \dots, \alpha_d \in \N_0, \atop \alpha_1 + \dots + \alpha_d \leq q}
     \left\|
\frac{
\partial^q f
}{
\partial x_1^{\alpha_1}
\dots
\partial x_d^{\alpha_d}
}
     \right\|_{\infty, [-A,A]^d}
     ,
     \]
 must hold for some sufficiently large constant $c_{13} \geq 1$).
 \\
Let $L, r \in \N$ be such that
\begin{enumerate}
\item $L \geq 5M^d+\left\lceil \log_4\left(M^{2p+4 \cdot d \cdot (q+1)} \cdot e^{4 \cdot (q+1) \cdot (M^d-1)}\right)\right\rceil \\
  \hspace*{4cm}
\cdot \lceil \log_2(\max\{q,d\}+1)\rceil+\lceil \log_4(M^{2p})\rceil$
\item $r \geq 132 \cdot 2^d\cdot   \lceil e^d\rceil
  \cdot \binom{d+q}{d} \cdot \max\{ q+1, d^2\}$
\end{enumerate}
hold.
There exists a feedforward neural network
$f_{net, deep}$ with ReLU activation function, $L$ hidden layers
and $r$ neurons per hidden layer 
where all weights are bounded in absolute value by
\[
e^{c_{14} \cdot (p+1) \cdot M^d}
\]
for some $c_{14}=c_{14}(f)>0 $, such that
\begin{align}
 \| f-f_{net, deep}\|_{\infty, [-A,A]^d} \leq
  c_{15} \cdot \left(\max\left\{A, \|f\|_{C^q([-A,A]^d)}\right\} \right)^{4(q+1)} \cdot M^{-2p}.
  \label{le3eq1}
\end{align}
 holds.
\end{lemma}

\noindent
    {\bf Proof.} This theorem is proven without the upper bound
on the absolute values of the weights in
Theorem 2 in Kohler
    and Langer (2021). It is explained in the supplement how the
above upper bound on the absolute value of the weights follows from
the proof given there.\hfill $\Box$

\begin{lemma}
  \label{le6}
  Let $d_1,d_2,l\in \N$ with $2^l\leq\min\{d_1,d_2\}$. For $k \in \{1, \dots,l\}$
  and $s \in \{1, \dots, 4^{l-k}\}$ let
  \[
\bar{g}_{net, k,s} : \R^4 \rightarrow \R
\]
be defined by a feedforward neural network with $L_{net}\in\N$
hidden layers and $r_{net}\in\N$ neurons per hidden layer and ReLU
activation function, where all the weights are bounded in absolute
value by some $B_n \geq 1$.
Set
\[I=\left\{0,\dots,2^{l}-1\right\} \times \left\{0, \dots, 2^{l}-1\right\}\]
 and
define
$\bar{m}:[0,1]^{\{1, \dots, d_1\} \times \{1, \dots, d_2\}} \rightarrow \R$
by
\[
\bar{m}(\bx)=
\max_{
  (i,j) \in \Z^2 \, : \,
  (i,j)+I \subseteq \{1, \dots, d_1\} \times \{1, \dots, d_2\}
}
\bar{f}\left(
x_{(i,j)+I}
\right),
\]
where $\bar{f}$ satisfies
    \[
\bar{f}=\bar{f}_{l,1}
    \]
    for some
    $\bar{f}_{k,s} :[-2,2]^{\{1, \dots, 2^k\} \times \{1, \dots, 2^k\}} \rightarrow \R$ recursively defined by
    \begin{eqnarray*}
    \bar{f}_{k,s}(\bx)&=&\bar{g}_{net,k,s} \big(
    \bar{f}_{k-1,4 \cdot (s-1)+1}(\bx_{
\{1, \dots, 2^{k-1}\} \times \{1, \dots, 2^{k-1}\}
    })
    , \\
        &&
        \hspace*{1cm}
        \bar{f}_{k-1,4 \cdot (s-1)+2}(\bx_{
\{2^{k-1}+1, \dots, 2^k\} \times \{1, \dots, 2^{k-1}\}
        }), \\
        &&
        \hspace*{1cm}
        \bar{f}_{k-1,4 \cdot (s-1)+3}(\bx_{
\{1, \dots, 2^{k-1}\} \times \{2^{k-1}+1, \dots, 2^k\}
        }), \\
        &&
        \hspace*{1cm}
        \bar{f}_{k-1,4 \cdot s}(\bx_{
\{2^{k-1}+1, \dots, 2^k\} \times \{2^{k-1}+1, \dots, 2^k\}
        })
    \big)
    \end{eqnarray*}
    for $k=2, \dots, l, s=1, \dots,4^{l-k}$,
 and
    \[
 \bar{f}_{1,s}(x_{1,1},x_{1,2},x_{2,1},x_{2,2})= \bar{g}_{net,1,s}(x_{1,1},x_{1,2},x_{2,1},x_{2,2})
 \]
 for $s=1, \dots, 4^{l-1}$.
 Set
 \[
l_{net}=\frac{4^{l}-1}{3} \cdot L_{net}+l,
\]
\[
k_s=\frac{2 \cdot 4^{l} + 4}{3} +r_{net}
\quad
(s=1, \dots, l_{net}),
\]
and set
\[M_s=2^{\pi(s)}\quad \mbox{for }s\in\{1,\dots,l_{net}\},\]
where the function $\pi:\{1,\dots,l_{net}\}\rightarrow\{1,\dots,l\}$ is defined by 
\[\pi(s)=\sum_{i=1}^{l}\IND_{\left\{s\geq i+\sum_{r=l-i+1}^{l-1}4^r\cdot L_{net}\right\}}.\]
 Then there exists some $m_{net} \in \F^{CNN}_{l_{net},\bk,\bM}$,
 where all the weights are bounded in absolute value by $B_n$, such that
 \[
\bar{m}(\bx) = m_{net}(\bx)
 \]
holds for all $\bx \in [-2,2]^{\{1, \dots, d_1\} \times \{1, \dots, d_2\}}$.
  \end{lemma}

\noindent
    {\bf Proof.} This theorem is proven without the upper bound
on the absolute values of the weights in
Lemma 2 in Kohler, Krzy\.zak
    and Walter (2022). 
In its proof the weights from the feedforward neural networks
are copied in the convolutional neural network, and all other
weights occurring are bounded in absolute value by $1$, therefore the
bound
on the absolute values of the weights holds.
\hfill $\Box$

\begin{lemma}
  \label{le7}
  Set
  \[
f(z)=\begin{cases}
\infty & ,z=1 \\
\log \frac{z}{1-z} & ,0<z<1 \\
- \infty & ,z=0,
\end{cases}
\]
let $K \in \N$ with $K \geq 6$, let $m:\R^{d \cdot l} \rightarrow [0,1]$
and let
$\bar{g}:\R^{d_1 \times d_2} \rightarrow \R$ such that \linebreak
$\|\bar{g}-m\|_{\infty,[0,1]^{d_1 \times d_2}} \leq \epsilon$ for some
\[
0 \leq \epsilon \leq \frac{1}{K}.
\]
Then there exists a neural
network $\bar{f}:\R \rightarrow \R$ with ReLU activation function,
and one hidden layer with $3 \cdot K + 9$ neurons,
where all the weights are bounded in absolute value by $K$,
such that for each network
$\tilde{f}:\R \rightarrow \R$
which has the same structure and which
has weights which are in supremum norm not more than
\[
0 \leq  \bar{\epsilon} \leq 1
\]
away from the weights of the above network
we have that 
$\tilde{f} \circ\bar{g}$
satisfies
\[
\| \tilde{f} \circ\bar{g}\|_{\infty,[0,1]^{d_1 \times d_2}}
\leq
132 \cdot K^2 \cdot \bar{\epsilon} + \log K
\]
and
\begin{eqnarray*}
  &&
\sup_{x \in [0,1]^{d_1 \times d_2}}
\Bigg(
\left|
m(x) \cdot
\left(
\varphi(\tilde{f}(\bar{g} ( x)))-
\varphi(f(m(x)))
\right)
\right|
\\
&&
\hspace*{4cm}
+
\left|
(1-m(x)) \cdot
\left(
\varphi(-\tilde{f}(\bar{g}(x)))-
\varphi(-f(m(x)))
\right)
\right|
\Bigg)
\\
&&
\leq
c_{16} \cdot \left(\frac{\log K}{K} + \epsilon\right) +
132 \cdot K^2 \cdot \bar{\epsilon}.
\end{eqnarray*}
\end{lemma}

\noindent
    {\bf Proof.} See Lemma 13 in Kohler and Krzy\.zak (2023).
\hfill $\Box$

\begin{lemma}
  \label{le8}
  Let $p \geq 1$ and $C>0$ be arbitrary.
  Assume that $\eta: \R^{d_1 \times d_2} \rightarrow \R$ satisfies
  a $(p,C)$--smooth hierarchical max-pooling model. Let
  $\F_n$ be the set of all CNNs with
  ReLU activation function, which have $L_n^{(1)}$ convolutional
  layers with $c_7$ channels in each layer, where $c_7$ is sufficiently
  large, one max pooling layer and one additional layer with
  $L_n^{(2)}$  neurons. Furthermore  assume $(L_n^{(1)})^{2p/4}
  \geq c_{17} \cdot L_n^{(2)}$.
Then there exists a network $f \in \F_n$ where all the weights
are bounded in absolute value by
\[
\max\left\{ L_n^{(2)}, e^{c_{18}(\eta) \cdot (p+1) \cdot L_n^{(1)} } \right\} 
\]
such that
\begin{eqnarray*}
&&
\EXP\left\{
\varphi( Y \cdot f(\bX))
\right\}
-
\EXP\left\{
\varphi( Y \cdot f_\varphi^* (\bX))
\right\}
\leq c_{19} \cdot \left(
\frac{\log L_n^{(2)}}{L_n^{(2)}} + \frac{1}{(L_n^{(1)})^{2p/4}}
\right)
\end{eqnarray*}
and such that $f$ is bounded in absolute value by $\log L_n^{(2)}$. 
  \end{lemma}

\noindent
{\bf Proof.} 
Use Lemma \ref{le5}
(applied with $d=4$ and
\[
M= \left\lfloor
\left(
\frac{3}{4^l-1} \cdot (L_n^{(1)}-l)
\right)^{\frac{1}{4}}
\right\rfloor,
\]
where $l$ is the level of the hierarchical max--pooling model for $\eta$)
and Lemma \ref{le6} to construct a convolutional
neural network $\bar{g}_{NN}$ built on the basis of feedforward neural
networks which approximate the functions in the hierarchical model
of the a posteriori probability $\eta$ in supremum norm up to an
error of order
\begin{equation}
  \label{ple8eq1}
\frac{1}{(L_n^{(1)})^{2p/4}}.
\end{equation}
Application of Lemma \ref{le4} with $t=1$ and $g(x)=\bar{g}(x)=x$
yields that $\bar{g}_{NN}$ approximates $\eta$ in supremum norm
by an error of order (\ref{ple8eq1}).
Next apply Lemma \ref{le7} (with $\epsilon=c_{20} \cdot
\frac{1}{(L_n^{(1)})^{2p/4}}$
and $\bar{\epsilon}=0$) to construct a neural network $\tilde{f}$
with
one hidden layer and $L_{n}^{(2)}$ neurons which takes on function values bounded in
absolute
value by $\log L_{n}^{(2)}$ and which satisfies
\begin{eqnarray*}
  &&
\sup_{x \in [0,1]^{d_1 \times d_2}}
\Bigg(
\left|
\eta(x) \cdot
\left(
\varphi(\tilde{f}(\bar{g}_{NN} ( x)))-
\varphi(f(\eta(x)))
\right)
\right|
\\
&&
\hspace*{4cm}
+
\left|
(1-\eta(x)) \cdot
\left(
\varphi(-\tilde{f}(\bar{g}_{NN}(x)))-
\varphi(-f(\eta(x)))
\right)
\right|
\Bigg)
\\
&&
\leq
c_{21} \cdot \left(\frac{\log L_{n}^{(2)}}{L_{n}^{(2)}} + \frac{1}{(L_n^{(1)})^{2p/4}}\right),
\end{eqnarray*}
where $f$ is the function defined in Lemma \ref{le7}.
Because of
\[
f_\varphi^*(x)=f(\eta(x))
\]
this implies
\begin{eqnarray*}
&&
\EXP\left\{
\varphi( Y \cdot \tilde{f}(\bar{g}_{NN} (\bX)))
\right\}
-
\EXP\left\{
\varphi( Y \cdot f_\varphi^* (\bX))
\right\}
\\
&&
=
\EXP\left\{
(1_{\{Y=1\}}+1_{\{Y=-1\}}) \cdot
\left(
\varphi( Y \cdot \tilde{f}(\bar{g}_{NN} (\bX)))
-
\varphi( Y \cdot f(\eta(\bX)))
\right)
\right\}
\\
&&
=
\EXP \Bigg\{
\eta(\bX) \cdot
\left(
\varphi(\tilde{f}(\bar{g}_{NN} (\bX)))-
\varphi(f(\eta(\bX)))
\right)
\\
&&
\hspace*{4cm}
+
(1-\eta(\bX)) \cdot
\left(
\varphi(-\tilde{f}(\bar{g}_{NN} (\bX)))-
\varphi(-f(\eta(\bX)))
\right)
\Bigg\}
\\
&&
\leq
\sup_{x \in [0,1]^{d_1 \times d_2}}
\Bigg(
\left|
\eta(x) \cdot
\left(
\varphi(\tilde{f}(\bar{g}_{NN} ( x)))-
\varphi(f(\eta(x)))
\right)
\right|
\\
&&
\hspace*{4cm}
+
\left|
(1-\eta(x)) \cdot
\left(
\varphi(-\tilde{f}(\bar{g}_{NN} (x)))-
\varphi(-f(\eta(x)))
\right)
\right|
\Bigg)
\\
&&
\leq
c_{21} \cdot \left(\frac{\log L_{n}^{(2)}}{L_{n}^{(2)}} + \frac{1}{(L_n^{(1)})^{2p/4}}\right).
\end{eqnarray*}
\hfill
 $\Box$

\subsubsection{Generalization error}

\begin{lemma}
\label{le9}
Let $\{f_\btheta \, : \, \btheta \in \bar{\bTheta}\}$ be defined
as in Subsection \ref{se3sub3} and let $\beta_n = c_1 \cdot \log n$.
Then we have
\begin{eqnarray*}
&&
 \sup_{x_1, \dots, x_n \in [0,1]^{d_1 \times d_2}}
  \EXP \left\{
  \left|
  \sup_{\btheta \in \bar{\bTheta}}
  \frac{1}{n} \sum_{i=1}^n \epsilon_i \cdot T_{\beta_n} f_\btheta(x_i)
  \right|
  \right\}
\\
&&
 \leq
c_{22} \cdot (\log n)^2 \cdot \frac{
\sqrt{
((L_n^{(1)})^2 + L_n^{(1)} \cdot L_n^{(2)})
\cdot
\log(\max\{L_n^{(1)},L_n^{(2)}\})
}
}{\sqrt{n}}.
\end{eqnarray*}
\end{lemma}

In the proof of Lemma \ref{le8} we will need the following bound
on the VC dimension of the class $\{f_\btheta \, : \, \btheta \in
\bar{\bTheta}\}$ of functions.

\begin{lemma}
\label{le10}
The VC dimension of the class $\{f_\btheta \, : \, \btheta \in
\bar{\bTheta}\}$ of functions in Lemma \ref{le9}
is bounded from above by
\[
c_{23} \cdot ((L_n^{(1)})^2 + L_n^{(1)} \cdot L_n^{(2)}) \cdot
\log(\max\{L_n^{(1)},L_n^{(2)}\}).
\]
\end{lemma}

\noindent
{\bf Proof.} The result follows from the proof of Lemma 7 in Kohler, Krzy\.zak and
Walter (2022).
\hfill $\Box$

\noindent
{\bf Proof of Lemma \ref{le9}.}
The results follows from Lemma \ref{le10} by an easy application
of standard techniques from VC theory. For the sake of completeness
we present nevertheless a detailed proof here.

Set $\F = \{f_{\btheta} \, : \, \btheta \in \bar{\Theta} \}$.
For $\delta_n>0$ and $x_1, \dots, x_n \in [0,1]^{d_1 \times d_2}$
we have
\begin{eqnarray*}
  &&
  \EXP \left\{
    \sup_{f \in \F}
    \left|
    \frac{1}{n} \sum_{i=1}^n \epsilon_i \cdot T_{\beta_n}(f(x_i))
    \right|
    \right\}
  \\
    &&
    =
    \int_0^{\beta_n}
    \PROB \left\{
    \sup_{f \in  \F}
    \left|
    \frac{1}{n} \sum_{i=1}^n \epsilon_i \cdot T_{\beta_n}(f(x_i))
    \right|
    > t
    \right\}
    \, dt
    \\
    &&
    \leq
    \delta_n + \int_{\delta_n}^{\beta_n}
    \PROB \left\{    
    \sup_{f \in  \F}
    \left|
    \frac{1}{n} \sum_{i=1}^n \epsilon_i \cdot T_{\beta_n}(f(x_i))
    \right|
    > t
    \right\}
    \, dt
.
\end{eqnarray*}
    Using a standard covering argument from empirical process theory we see that
    for any $t \geq \delta_n$ we have
    \begin{eqnarray*}
      &&
    \PROB \left\{    
    \sup_{f \in  \F}
    \left|
    \frac{1}{n} \sum_{i=1}^n \epsilon_i \cdot T_{\beta_n}(f(x_i))
    \right|
    > t
    \right\}
    \\
    &&
    \leq
    \Mu_1 \left(
    \frac{\delta_n}{2},
    \left\{
T_{\beta_n} f : f \in \F
    \right\}, x_1^n
    \right)
    \\
    &&
    \hspace*{3cm}
    \cdot
    \sup_{f \in \F}
    \PROB \left\{
    \left|
    \frac{1}{n} \sum_{i=1}^n \epsilon_i \cdot T_{\beta_n}(f(x_i))
    \right|>\frac{t}{2}
    \right\}.
    \end{eqnarray*}
Application of Lemma \ref{le10} and Theorem 9.4 in Gy\"orfi et
al. (2002)
yields
    \[
        \Mu_1 \left(
    \frac{\delta_n}{2},
    \left\{
T_{\beta_n} f : f \in  \F
    \right\}, x_1^n
    \right)
    \leq c_{24} \cdot
    \left(
    \frac{c_{25}  \cdot \beta_n}{\delta_n}
    \right)^{c_{26}  \cdot 
((L_n^{(1)})^2 + L_n^{(1)} \cdot L_n^{(2)}) \cdot
\log(\max\{L_n^{(1)},L_n^{(2)}\})
}.
    \]    
    By the inequality of Hoeffding (cf., e.g., Lemma A.3 in Gy\"orfi et al. (2002))
    and
    \[
    |T_{\beta_n}(f(x))| \leq  \beta_n
    \quad (x \in \R^d)
    \]
we have for any $f \in  \F$
\[
\PROB \left\{
     \left|
    \frac{1}{n} \sum_{i=1}^n \epsilon_i \cdot 
    T_{\beta_n}(f(x_i))
    \right|
    > t
    \right\}
    \leq 2 \cdot \exp \left(
- \frac{2 \cdot n \cdot t^2}{4  \cdot \beta_n^2}
    \right).
\]
Hence we get
\begin{eqnarray*}
  &&
 \sup_{x_1, \dots, x_n \in [0,1]^{d_1 \times d_2}}
  \EXP \left\{
    \sup_{f \in  \F}
    \left|
    \frac{1}{n} \sum_{i=1}^n \epsilon_i \cdot T_{\beta_n}(f(X_i))
    \right|
    \right\}
  \\
    &&
  \leq
  \delta_n+ \int_{\delta_n}^{\beta_n}
c_{24} \cdot
    \left(
    \frac{c_{25}  \cdot \beta_n}{\delta_n}
    \right)^{c_{26}  \cdot 
((L_n^{(1)})^2 + L_n^{(1)} \cdot L_n^{(2)}) \cdot
\log(\max\{L_n^{(1)},L_n^{(2)}\})
    }
    \\
    &&
    \hspace*{7cm}
 \cdot    2 \cdot \exp \left(
- \frac{ n \cdot \delta_n \cdot t}{2  \cdot \beta_n^2}
\right) \, dt
\\
&&
\leq
\delta_n+
c_{24} \cdot
    \left(
    \frac{c_{25}  \cdot \beta_n}{\delta_n}
    \right)^{c_{26}  \cdot 
((L_n^{(1)})^2 + L_n^{(1)} \cdot L_n^{(2)}) \cdot
\log(\max\{L_n^{(1)},L_n^{(2)}\})
} 
    \frac{4  \cdot \beta_n^2}{ n \cdot \delta_n}
    \cdot
    \exp \left(
- \frac{ n \cdot \delta_n^2}{2  \cdot \beta_n^2}
\right). 
  \end{eqnarray*}
With
\[
\delta_n=
\sqrt{ ((L_n^{(1)})^2 + L_n^{(1)} \cdot L_n^{(2)}) \cdot
\log(\max\{L_n^{(1)},L_n^{(2)}\}) } \cdot \log n
\cdot
\sqrt{
\frac{2 \cdot  \beta_n^2 }{ n}
}
\]
we get the assertion.
\hfill $\Box$

\subsubsection{A bound on the gradient}
\begin{lemma}
  \label{le11}
  Let $A$, $\bar{\Theta}$ and $f_{(\bw,\vartheta)}$ be defined
  as in Section \ref{se3}, set
  \[
  M_{\max}=\max\{M_1, \dots, M_{L_n^{(1)}}\}
  \quad \mbox{and} \quad
  k_{max}=\max\{ k_1, \dots, k_{L_n^{(1)}} \}
  \]
  and assume that all weights in $f_{(\bw,\vartheta)}$ are bounded in absolute
  value by $B_n \geq 1$. Then
\begin{eqnarray*}
&&
\sup_{\bw \in A, \btheta \in \bar{\Theta}^{K_n}, y \in \{-1,1\}, x \in [0,1]^{d_1 \times d_2}} \| \nabla_{\vartheta}
\varphi( y \cdot f_{(\bw,\vartheta)}(x)) \|_\infty
\\
&&
\leq 
L_n^{(2)} \cdot k_{max}^{2 \cdot L_n^{(1)}+1}
\cdot (M_{max}^2+1)^{2 \cdot L_n^{(1)}+2}
\cdot B_n^{2 \cdot L_n^{(1)}+2}.
\end{eqnarray*}
\end{lemma}

\noindent
    {\bf Proof.}
    Let
    \[
f_{(\bw, \btheta)}(x)=
\sum_{k=1}^{K_n} w_k \cdot T_{\beta_n} f_{\btheta_k}(x)
\]
where
\[
 f_{\btheta_k}(x) = g_{\bw_k}(f_{\bw_k}(x)),
 \]
 \[
 g_{\bw_k}(z)
 =\sum_{i=1}^{L_n^{(2)}} (\bw_k)_{i}^{(1) } \sigma \left(
(\bw_k)_{i,1}^{(0)} \cdot z + (\bw_k)_{i,0}^{(0)}
\right)+(\bw_k)_{0}^{(1)},
\]
\begin{eqnarray*}
f_{\bw_k}(x)
&=&
\max\Bigg\{\sum_{s_2=1}^{k_{L_n^{(1)}}} (\bw_k)_{s_2} \cdot (o_{(\bw_k)})_{(i,j), s_2}^{(L_n^{(1)})}: i \in \{1, \dots, d_1-M_{L_n^{(1)}}+1\},\\
&& \hspace*{7cm} j \in \{1, \dots, d_2-M_{L_n^{(1)}}+1\}\Bigg\},
\end{eqnarray*}
\[
(o_{(\bw_k)})_{(i,j), s_2}^{(l)} = \sigma\left(\sum_{s_1=1}^{k_{l -1}} \sum_{\substack{t_1, t_2 \in \{1, \dots, M_l\} \\ (i+t_1-1, j+t_2-1) \in D}} (\bw_k)_{t_1, t_2, s_1, s_2}^{(l)} (o_{(\bw_k)})_{(i+t_1-1, j+t_2-1), s_1}^{(l-1)} + (\bw_k)_{s_2}^{(l)}\right)
\]
for $l \in \{1, \dots, L_n^{(1)}\}$,
and
\[
(o_{(\bw_k)})_{(i,j), 1}^{(0)} = x_{i,j} \quad \text{for} \ i \in \{1, \dots, d_1\} \ \text{and} \ j \in \{1, \dots, d_2\}.
\]
By the proof of Lemma \ref{le3} we know for any $l \in \{1, \dots, L_n^{(1)} \}$
\begin{equation}
  \label{ple11eq1}
|(o_{(\bw_k)})_{(i,j),s_2}^{(l)}(x)|
\leq
k_{max}^{L_n^{(1)}} \cdot (M_{max}^2+1)^{L_n^{(1)}} \cdot B_n^{L_n^{(1)}}. 
\end{equation}
Using the chain rule we get for any 
$l \in \{1, \dots, L_n^{(1)} \}$ and suitably chosen (random)
$(i_0,j_0) \in D$
\begin{eqnarray*}
  &&
  \frac{\partial}{\partial (\bw_k)_{\tilde t_1,\tilde t_2,\tilde s_1,\tilde s_2}^{(l)} }
  f_{(\bw, \btheta)}(x)
  \\
  &&
  =
  w_k
  \cdot
  \frac{\partial}{\partial z} T_{\beta_n} z \Big|_{z=f_{\btheta_k}(x)}
  \cdot
  \sum_{i=1}^{L_n^{(2)}} (\bw_k)_{i}^{(1) } \sigma^\prime \left(
(\bw_k)_{i,1}^{(0)} \cdot f_{\bw_k}(x) + (\bw_k)_{i,0}^{(0)}
  \right)
  \cdot
  (\bw_k)_{i,1}^{(0)} 
  \cdot
  \sum_{s_2=1}^{k_{L_n^{(1)}}} (\bw_k)_{s_2} 
  \\
  &&
  \quad
  \cdot
  \sigma^\prime
  \left(\sum_{s_1=1}^{k_{L_n^{(1)}-1}} \sum_{\substack{t_1, t_2 \in \{1, \dots, M_{L_n^{(1)}}\}
      \\ (i_0+t_1-1, j_0+t_2-1) \in D}} (\bw_k)_{t_1, t_2, s_1, s_2}^{(L_n^{(1)})} (o_{(\bw_k)})_{(i_0+t_1-1, j_0+t_2-1), s_1}^{(L_n^{(1)}-1)} + (\bw_k)_{s_2}^{(L_n^{(1)})}\right)
  \\
  &&
  \quad
  \cdot
\sum_{s_1^{(L_n^{(1)})}=1}^{k_{L_n^{(1)}-1}} \sum_{\substack{t_1^{(L_n^{(1)})}, t_2^{(L_n^{(1)})} \in \{1, \dots, M_{L_n^{(1)}}\}
    \\ (i_0+t_1^{(L_n^{(1)})}-1, j_0+t_2^{(L_n^{(1)})}-1) \in D}} (\bw_k)_{t_1^{(L_n^{(1)})}, t_2^{(L_n^{(1)})}, s_1^{(L_n^{(1)})}, s_2}^{(L_n^{(1)})}
\\
&&
\quad
\cdot
  \sigma^\prime
  \Bigg(\sum_{s_1=1}^{k_{L_n^{(1)}-2}} \sum_{\substack{t_1, t_2 \in \{1, \dots, M_{L_n^{(1)}-1}\}
      \\ (i_0+t_1^{(L_n^{(1)})}+t_1-2, j_0+t_2^{(L_n^{(1)})}+t_2-2) \in D}} (\bw_k)_{t_1, t_2, s_1, s_1^{(L_n^{(1)})}}^{(L_n^{(1)}-1)} \cdot
  \\
  &&
  \hspace*{4cm}
  (o_{(\bw_k)})_{(i_0+t_1^{(L_n^{(1)})}+t_1-2, j_0+t_2^{(L_n^{(1)})}+t_2-2), s_1}^{(L_n^{(1)}-2)} + (\bw_k)_{s_1^{(L_n^{(1)})}}^{(L_n^{(1)}-1)}\Bigg)
  \\
  &&
  \quad
  \cdot
  \sum_{s_1^{(L_n^{(1)}-1)}=1}^{k_{L_n^{(1)}-2} } \sum_{\substack{t_1^{(L_n^{(1)}-1)}, t_2^{(L_n^{(1)}-1)} \in \{1, \dots, M_{L_n^{(1)}-1}\}
      \\ (i_0+t_1^{(L_n^{(1)})}+t_1^{(L_n^{(1)}-1)}-2, j_0+t_2^{(L_n^{(1)})}+t_2^{(L_n^{(1)}-1)}-2) \in D}} (\bw_k)_{t_1^{(L_n^{(1)}-1)}, t_2^{(L_n^{(1)}-1)}, s_1^{(L_n^{(1)}-1)}, s_1^{(L_n^{(1)})}}^{(L_n^{(1)}-1)}
  \\
  &&
  \quad
  \cdot \ldots \cdot
  \sigma \Bigg(
  \sum_{s_1=1}^{k_{l -1}}
 \sum_{{t_1, t_2 \in \{1, \dots, M_l\}, (i_0
+t_1^{(L_n^{(1)})}+ \dots + t_1^{(l+1)}+ \tilde t_1
+
t_1-(L_n-(l-2)), \atop
\hspace*{3cm}
j_0
+t_2^{(L_n^{(1)})}+ \dots +t_2^{(l+1)}+ \tilde t_2
+t_2-(L_n^{(1)}-(l-2))) \in D}}
  \\
  &&
  (\bw_k)_{t_1, t_2, s_1, \tilde s_1}^{(l-1)}
\cdot
  (o_{(\bw_k)})_{(i_0
+t_1^{(L_n^{(1)})}+ \dots + t_1^{(l+1)}+\tilde t_1
+
t_1-(L_n^{(1)}-(l-2)),
\atop
j_0
+t_2^{(L_n^{(1)})}+ \dots +t_2^{(l+1)} + \tilde t_2+
t_2-(L_n^{(1)}-(l-2))), s_1}^{(l-2)}
+ (\bw_k)_{\tilde s_1}^{(l-1)}
  \Bigg).
  \end{eqnarray*}
Using (\ref{ple11eq1}), $|\sigma^\prime(z)| \leq 1$ and that all weights
are bounded in  the absolute
value by $B_n$ we get
\[
\left|
  \frac{\partial}{\partial (\bw_k)_{t_1,t_2,s_1,s_2}^{(l)} }
  f_{(\bw, \btheta)}(x)
  \right|
  \leq
L_n^{(2)} \cdot k_{max}^{2 \cdot L_n^{(1)}+1}
\cdot (M_{max}^2+1)^{2 \cdot L_n^{(1)}+2}
\cdot B_n^{2 \cdot L_n^{(1)}+2}.  
\]
Analogously we can derive bounds on all the other partial
derivatives occurring in the assertion.
    \hfill $\Box$

\subsubsection{Proof of Theorem \ref{th2}}
It suffices to show
\begin{equation}
\label{pth2eq1}
 \EXP \left\{ \varphi(Y \cdot f_n(X)) \right\}
  -
  \min_{f : [0,1]^{d_1 \times d_2} \rightarrow \bar{\R}}
  \EXP \left\{ \varphi(Y \cdot f(X)) \right\}
 \leq
c_{28} \cdot (\log n)^{4} \cdot n^{-\min\{\frac{p}{2p+4},
   \frac{1}{4} \}}
\end{equation}
for $n$ sufficiently large.

        This implies the assertion, because by Lemma \ref{le2} a)
        we conclude from (\ref{pth2eq1})
        \begin{eqnarray*}
          &&
          \PROB \left\{
Y \neq \hat{C}_n(X)  \right\}
-
\PROB \left\{
Y \neq f^*(X) \right\}
\\
&&
\leq
\EXP \left\{
\frac{1}{\sqrt{2}}
\cdot
\left(
\EXP\left\{
\varphi( Y \cdot f_n(X))
| \D_n
\right\}
-
\EXP\left\{
\varphi( Y \cdot f_{\varphi^*}(X))
\right\}
\right)^{1/2}
\right\}
\\
&&
\leq
\frac{1}{\sqrt{2}}
\cdot
\sqrt{
\EXP\left\{
\varphi( Y \cdot f_n(X))
\right\}
-
\EXP\left\{
\varphi( Y \cdot f_{\varphi^*}(X))
\right\}
}
\\
&&
\leq
 c_8 \cdot (\log n)^2 \cdot n^{-\min\{\frac{p}{2 \cdot (2p+4)},
   \frac{1}{8} \}}
          \end{eqnarray*}
And from Lemma \ref{le2} b), (\ref{th2eq1}) and Lemma \ref{le2} c)
we conclude from (\ref{pth2eq1})
\begin{eqnarray*}
  &&
        \PROB \left\{
Y \neq \hat{C}_n(X)\right\}
-
\PROB \left\{
Y \neq f^*(X) \right\}
\\
&&
\leq
2 \cdot
\left(
\EXP\left\{
\varphi( Y \cdot f_n(X))
\right\}
-
\EXP\left\{
\varphi( Y \cdot f_{\varphi^*}(X))
\right\}
\right)
+ 4 \cdot \frac{c_{29} \cdot \log n}{n^{1/4}}
\\
&&
\leq
c_9 \cdot (\log n)^{4} \cdot n^{-\min\{\frac{p}{2p+4},
   \frac{1}{4} \}}.
\end{eqnarray*}
Here we have used the fact that
\[
\max
\left\{
\frac{\PROB\{Y=1|X\}}{1-\PROB\{Y=1|X\}},
\frac{1-\PROB\{Y=1|X\}}{\PROB\{Y=1|X\}}
\right\}
> n^{1/4}
\]
is equivalent to
\[
|f_{\varphi^*}(X)|
=
\left|
\log \frac{\PROB\{Y=1|X\}}{1-\PROB\{Y=1|X\}}
\right|
> \frac{1}{4} \cdot \log n.
\]

In the remainder of the proof we apply
Theorem \ref{th1} in order to prove (\ref{pth2eq1}).
Here we assume throughout the proof that $n$ is sufficiently large.
Let $f_{\vartheta}$ be the network of Lemma \ref{le8}
which satisfies
\begin{eqnarray*}
&&
\EXP\left\{
\varphi( Y \cdot f_{\vartheta}(\bX))
\right\}
-
\min_{f : [0,1]^{d_1 \times d_2} \rightarrow \bar{\R}}
  \EXP \left\{ \varphi(Y \cdot f(X)) \right\}
\leq c_{30} \cdot \left(
\frac{\log L_n^{(2)}}{L_n^{(2)}} + \frac{\log L_n^{(2)}}{(L_n^{(1)})^{2p/4}}
\right).
\end{eqnarray*}
Since $f_{\vartheta}$ is bounded in supremum norm by $\log L_n^{(2)}
\leq \beta_n$ (for $c_3$ sufficiently large) this implies
\begin{eqnarray*}
&&
\EXP\left\{
\varphi( Y \cdot T_{\beta_n} f_{\vartheta}(\bX))
\right\}
-
\min_{f : [0,1]^{d_1 \times d_2} \rightarrow \bar{\R}}
  \EXP \left\{ \varphi(Y \cdot f(X)) \right\}
\leq c_{30} \cdot \left(
\frac{\log L_n^{(2)}}{L_n^{(2)}} + \frac{\log L_n^{(2)}}{(L_n^{(1)})^{2p/4}}
\right).
\end{eqnarray*}
Here we have used that we can decrease the value of $L_n^{(2)}$ in order to ensure that the assumption $(L_n^{(1)})^{2p/4}  \geq c_{17} \cdot L_n^{(2)}$ (cf., Lemma \ref{le8}) holds.\\
Set
\[
\delta_n=\frac{1}{n \cdot e^{n}}
\]
and choose
\[
\Theta^*=\left\{
\vartheta^* \in \Theta \, : \, \|\vartheta^*-\vartheta\|_\infty \leq \delta_n
\right\}.
\]
Then it follows from the Lipschitz continuity of the logistic loss and
Lemma \ref{le3} that we have for any $\vartheta^* \in \Theta^*$
\[
\left|
\varphi( Y \cdot T_{\beta_n} f_{\vartheta^*}(\bX))
-
\varphi( Y \cdot T_{\beta_n} f_{\vartheta}(\bX))
\right|
\leq
c_{31} \cdot e^{n} \cdot \| \vartheta^*-\vartheta\|_\infty \leq \frac{c_{32}}{n},
\]
from which we can conclude
\begin{eqnarray*}
&&
\sup_{\vartheta^* \in \Theta^*}
\EXP\left\{
\varphi( Y \cdot T_{\beta_n} f_{\vartheta^*}(\bX))
\right\}
-
\min_{f : [0,1]^{d_1 \times d_2} \rightarrow \bar{\R}}
  \EXP \left\{ \varphi(Y \cdot f(X)) \right\}\\
&&\leq c_{30} \cdot \left(
\frac{\log L_n^{(2)}}{L_n^{(2)}} 
+ \frac{\log L_n^{(2)}}{(L_n^{(1)})^{2p/4}}
\right) + \frac{c_{32}}{n}.
\end{eqnarray*}
Furthermore the definition of $\Theta^*$ implies that
\[
  \kappa_n = \PROB \left\{
\btheta^{(0)}_1 \in \bTheta^*
\right\}
\geq \left(
\frac{1}{n \cdot e^{n}}
\right)^{c_{33} \cdot (L_n^{(1)}+L_n^{(2)})} \geq e^{-n^{1.5}}. 
\]
So if we choose
\[
N_n= n^2 \cdot e^{2 \cdot n}, \quad I_n=n^2 \cdot e^{n^{1.5}}
\]
(which is possible because of $N_n \cdot I_n \leq K_n$)
then we have
\[
N_n \cdot (1- \kappa_n)^{I_n} \leq\frac{1}{n}, 
\]
so (\ref{th1eq3}) holds.

By Lemma \ref{le3} we know that (\ref{th1eq1}) is satisfied for
\[
C_n= c_{34} \cdot e^{n},
\]
and because of
\begin{eqnarray*}
&&
\| \nabla_{\bw} \varphi(y \cdot f_{(\bw,\vartheta)}(x))\|^2
\leq \sum_{k=1}^{K_n} |1 \cdot T_{\beta_n} f_{\vartheta_k}(x)|^2
\leq K_n \cdot \beta_n^2
\end{eqnarray*}
(\ref{th1eq4}) is satisfied for
\[
D_n= \sqrt{K}_n \cdot \beta_n.
\]
By Lemma \ref{le11} we know
\begin{eqnarray*}
&&
\sup_{\bw \in A, \btheta \in \bar{\Theta}^{K_n}, y \in \{-1,1\}, x \in [0,1]^{d_1 \times d_2}} \| \nabla_{\vartheta}
   \varphi( y \cdot f_{(\bw,\vartheta)}(x)) \|_\infty
\leq e^{n}.
\end{eqnarray*}
Application of Theorem \ref{th1} together with
Lemma \ref{le9} and
the above results
yields
\begin{eqnarray*}
  &&
  \EXP \left\{ \varphi(Y \cdot f_n(X)) \right\}
  -
  \min_{f : \R^{d_1 \times d_2} \rightarrow \bar{\R}}
  \EXP \left\{ \varphi(Y \cdot f(X)) \right\}
 \\
  &&
 \leq
 c_{35} \cdot
\Bigg(
\frac{\log n}{n}
+
(\log n)^3
\frac{
\sqrt{
\left(
n^{\frac{4}{2p+4}}+n^{\frac{2}{2p+4}} \cdot n^{\frac{1}{4}}
\right) \cdot \log n
}
}{\sqrt{n}}
+
\frac{e^{n}}{n \cdot e^{n}}
+
\frac{K_n \cdot \beta_n^2}{t_n}
\\
&&
\quad
+
\frac{
n \cdot \beta_n^2
\cdot \left(
K_n + e^{n} \cdot e^{n}
\right)
}{t_n}
+
\frac{\log L_n^{(2)}}{L_n^{(2)}} + \frac{\log L_n^{(2)}}{(L_n^{(1)})^{2p/4}}
\Bigg)
\\
&&
\leq
c_{28} \cdot (\log n)^{4} \cdot n^{-\min\{\frac{p}{2p+4},
   \frac{1}{4} \}}. 
\end{eqnarray*}

\hfill $\Box$


\section{Acknowledgment}
The second author would like to thank
Natural Sciences and Engineering Research Council of Canada for
funding this project under Grant RGPIN-2020-06793. 


\newpage

\begin{appendix}

\section{Supplement: Proof of Lemma \ref{le5}}
In the following we prove the weight constraints of Lemma \ref{le5} by
modifying the auxiliary results of the proof of Kohler and Langer (2021). 
\subsection{Further notation and definitions}
The following auxiliary notation is required for the statement of these resuts:

We introduce our framework for a fully-connected neural network $g:\Rd\to \R$ with ReLU activation function $\sigma(x)=\max\{x,0\}$:
Let 
\begin{equation}
g(x)=
\sum_{i=1}^{k_L} v_{1,i}^{(L)}g_i^{(L)}(x)
+v_{1,0}^{(L)},
\label{FNN-lastlayer}
\end{equation}
where output weights $v_{1,0}^{(L)},\dots,v_{1,k_L}^{(L)}\in\R$ denote the output weights of the network.\\
The outputs of the neurons $g_i^{(L)}$ are recursively defined by
\[
g_i^{(r)}(x)=
\sigma\left(
\sum_{j=1}^{k_{r-1}}v_{i,j}^{(r-1)}g_j^{(r-1)}(x)
+v_{i,0}^{(r-1)}
\right),
\]
with inner weights $v_{i,0}^{(r-1)},\dots,v_{i,k_{r-1}}^{(r-1)}\in\R$ for $i\in\{1,\dots,k_r\}$, $r\in\{1,\dots,L\}$, $k_0=d$ and
\[
g_j^{(0)}(x)=x^{(j)}.
\]
A fully-connected neural network of the form \eqref{FNN-lastlayer} is dependent on the number of layers $L$ and a width vector $\mathbf{k}=(k_1, \dots, k_L)$, hence we will denote the corresponding function class by  $\mathcal{F}^{FNN}_{ L, \bk}$. In case $k_1 = \ldots = k_L= r$, we write $\mathcal{F}^{FNN}_{ L, r}$ to indicate that all layers consist of a constant number of $r$ neurons.  Further, we denote by $\bv_f$ the vector that collects all weights required for the computation of $f\in \mathcal{F}^{FNN}_{ L, \bk}$: 
$$
\bv_f=
\left(
\left(v_{i,j}^{(l)}\right)_{
                            i\in\left\{1, \dots, k_{l+1}\right\},
                            \, 
                            j \in\left\{0, \dots, k_l\right\},
                            \,
                            l\in \{0,\dots, L-1\}
                            },
\left(v_{1,i}^L\right)_{
                        i\in\left\{1, \dots, k_L\right\}
                        }
\right).
$$
The proof is based on an approximation by a
piecewise Taylor polynomial, which is defined using a partition into equivolume cubes.
If $C$ is a cube, then $C_{left}$ is used to denote the "bottom left" of $C$. We can thus write each half-open cube $C$ with side length $s$ as a polytope defined by 
\begin{align*}
-x^{(j)} + C_{left}^{(j)} \leq 0 \ 
\mbox{and} 
\ x^{(j)} - C_{left}^{(j)}-s < 0 \quad (j \in \{1, \dots, d\}).
\end{align*}
Furthermore, we describe by $C_{\delta}^0 \subset C$ the cube, which contains all $x \in C$ that lie with a distance of at least $\delta$ to the boundaries of $C$, i.e. $C_{\delta}^0$ is the polytope defined by
\begin{align*}
-x^{(j)} + C_{left}^{(j)} \leq - \delta \ \mbox{and} \ x^{(j)} - C_{left}^{(j)}-s < -\delta \quad (j \in \{1, \dots, d\}).
\end{align*}
If $\P$ is a partition of cubes of $[-a,a)^d$
and $x \in [-a,a)^d$, then we denote the cube $C \in \P$, 
which satisfies $x \in C$, by $C_\P (x)$. 

Let $\mathcal{P}_N$ be the linear span of all monomials of the form 
\begin{align*}
\prod_{k=1}^d \left(x^{(k)}\right)^{r_k}
\end{align*}
for some $r_1, \dots, r_d \in \N_0$, $r_1+\dots+r_d \leq N$. Then, $\mathcal{P}_N$ is a linear vector space of functions of dimension 
\begin{align*}
dim \ \mathcal{P}_N = \left|\left\{(r_0, \dots, r_d) \in \N_0^{d+1}: r_0+\dots+r_d = N \right\}\right| = \binom{d+N}{d}.
\end{align*}
\subsection{Auxiliary results}
\begin{lemma}
\label{le_appendix1a}
Let $p=q+s$ for some $q \in \N_0$ and $s \in (0,1]$, and let $C > 0$. Let $f: \Rd \to \R$ be a $(p,C)$-smooth function, let $x_0 \in \Rd$ and let $T_{f,q,x_0}$ be the Taylor polynomial of total degree $q$ around $x_0$ defined by
  \begin{eqnarray*}
T_{f,q,x_0}(x) &=& \sum_{\substack{j \in \N_0^d: \|\bold{j}\|_1 \leq q}} (\partial^{\bold{j}} f) ( x_0)\cdot \frac{\left(x - x_0\right)^{\bold{j}}}{\bold{j}!}
  \end{eqnarray*}
Then for any $x \in \Rd$ 
\begin{align*}
\left|f(x) - T_{f,q,x_0}(x)\right| \leq c_{35} \cdot C \cdot \Vert x - x_0 \Vert^p
\end{align*}
holds for a constant $c_{35}=c_{35}(q,d)$ depending only on $q$ and $d$.
\end{lemma}

\noindent
{\bf Proof.}
See Lemma 1 in Kohler (2014).
\hfill $\Box$

In the proof of Lemma \ref{le5} we use Lemma \ref{le_appendix1a} and approximate our function by a piecewise Taylor polynomial. To define this piecewise Taylor polynomial, we partition  $[-a,a)^d$
  into $M^d$ and $M^{2d}$ half-open
  equivolume cubes of the form
  \[
[\bm{\alpha},\bm{\beta})=[\bm{\alpha}^{(1)},\bm{\beta}^{(1)}) \times \dots \times [\bm{\alpha}^{(d)},\bm{\beta}^{(d)}), \quad \bm{\alpha}, \bm{\beta} \in \Rd,
  \]
respectively.
Let 
\begin{align}
\label{partition}
\mathcal{P}_1=\{C_{k,1}\}_{k \in \{1, \dots, M^d\}} \ \mbox{and} \ \mathcal{P}_2=\{C_{j,2}\}_{j \in \{1, \dots, M^{2d}\}}
\end{align}
be the corresponding partitions. 
We denote for each $i \in \{1, \dots, M^d\}$
those cubes of $\mathcal{P}_2$ that are contained in $C_{i,1}$
by $\tilde{C}_{1, i}, \dots, \tilde{C}_{M^d, i}$ and order the cubes in such a way that the bottom left of $\tilde{C}_{1,i}$ and $C_{i,1}$ coincide, i.e. such that we have $(\tilde{C}_{1,i})_{left}=(C_{i,1})_{left}$ and that
\begin{align}
\label{tv}
(\tilde{C}_{k,i})_{left} = (\tilde{C}_{k-1,i})_{left}+\bold{\tilde{v}}_k
\end{align}
holds for all $k \in \{2, \dots, M^d\}, i \in \{1, \dots, M^d\}$ and some vector $\bold{\tilde{v}}_k$ with entries in $\{0, 2a/M^2\}$ where exactly one entry is different to zero. Here the vector $\bold{\tilde{v}}_k$ describes the position of $(C_{k,i})_{left}$ relative to $(C_{k-1, i})_{left}$ and we order the cubes in such a way that the position is independent of $i$.
Then Taylor expansion in Lemma \ref{le_appendix1a} can be used to define a piecewise Taylor polynomial on $\P_2$ by
\[
T_{f,q,(C_{\P_2}(x))_{left}}(x)
=
\sum_{k \in \{1, \dots, M^d\}, i \in \{1, \dots, M^d\} } 
T_{f,q,(\tilde{C}_{k,i})_{left}}(x)
\cdot 
\mathds{1}_{\tilde{C}_{k,i}}(x)
\]
and this piecewise Taylor polynomial satisfies
\begin{align*}
  \sup_{x \in [-a,a)^d} 
  \left|f(x) - T_{f,q,((C_{\P_2}(x))_{left}}(x)\right| 
  \leq
    c_{35} \cdot 
    C \cdot 
    (2 \cdot a \cdot d)^{p} 
    \cdot \frac{1}{M^{2p}}.
\end{align*}
To compute $T_{f,q,(C_{\P_2}(x))_{left}}(x)$ the very deep neural network of Lemma \ref{le5} b) proceeds in two steps: In a first step it computes
$(C_{\P_1}(x))_{left}$ and
the values of
\begin{align*}
(\partial^{\bll} f)((C_{i,1})_{left})
\end{align*}
for each $\bll \in \N_0^d$ with $\|\bll\|_1 \leq q$ and suitably defined numbers 
\begin{align*}
b_{k,i}^{(\bll)} \in \Z, \quad |b_{k,i}^{(\bll)}| \leq e^d+1 \quad (k \in \{1, \dots, M^d\}),
\end{align*}
which depend on $C_{i,1}$ for $i \in \{1, \dots, M^d\}$.
Assume that $x \in C_{i,1}$ for some $i \in \{1, \dots, M^d\}$. In the second step the neural network successively computes approximations 
\begin{align*}
(\partial^{\bll} \hat{f})((\tilde{C}_{k,i})_{left}), \quad k \in \{1, \dots, M^d\}
\end{align*}
of 
\begin{align*}
(\partial^{\bll} f)((\tilde{C}_{k,i})_{left})
\end{align*}
for each $\bll \in \N_0^d$ with $\|\bll\|_1 \leq q$. To do this we start with
\begin{align*}
(\partial^{\bll} \hat{f})((\tilde{C}_{1,i})_{left}) = (\partial^{\bll} f)((C_{\P_1}(x))_{left}).
\end{align*} 
By construction of the first step and since $(\tilde{C}_{1,i})_{left} = (C_{\P_1}(x))_{left}$ these estimates have error zero. As soon as we have computed the above estimates for some $k \in \{1, \dots, M^d-1\}$ we use the Taylor polynomials with these coefficients around $(\tilde{C}_{k,i})_{left}$ in order to compute 
\begin{align*}
&\sum_{\substack{\bj \in \N_0^d:\\ \|\bj\|_1 \leq q-\|\bll\|_1}} \frac{(\partial^{\bll+\bj} \hat{f})((\tilde{C}_{k,i})_{left})}{\bj!} \cdot \left((\tilde{C}_{k+1,i})_{left} - (\tilde{C}_{k,i})_{left}\right)^{\bj}
\end{align*}
for $\bll \in \N_0^d$ with $\|\bll\|_1 \leq q$ and we define
\begin{align*}
(\partial^{\bll} \hat{f})((\tilde{C}_{k+1,i})_{left}) =&\sum_{\substack{\bj \in \N_0^d:\\ \|\bj\|_1\leq q-\|\bll\|_1}} \frac{(\partial^{\bll+\bj}\hat{f})((\tilde{C}_{k,i})_{left})}{\bj!} \cdot \left((\tilde{C}_{k+1,i})_{left} - (\tilde{C}_{k,i})_{left}\right)^{\bj}\\
&+b_{k,i}^{(\bll)} \cdot c_{36} \cdot \left(\frac{2a}{M^2}\right)^{p-\|\bll\|_1}
\end{align*}
where
\[
c_{36}= C \cdot d^p \cdot \max\{
c_{35}(q,d),c_{35}(q-1,d), \dots, c_{35}(0,d)
\}
\]
(and $c_{35}$ is the constant of Lemma \ref{le_appendix1a}).
Assume that
\begin{align*}
\left|(\partial^{\bll}\hat{f})((\tilde{C}_{k,i})_{left}) - (\partial^{\bll} f) ((\tilde{C}_{k,i})_{left})\right| \leq c_{36} \cdot \left(\frac{2a}{M^2}\right)^{p-\|\bll\|_1}
\end{align*}
holds for all $\bll \in \N_0^d$ with $\|\bll\|_1 \leq q$ (which holds by construction for $k=1$). Then 
\begin{align*}
&\Bigg|\sum_{\substack{\mathbf{j}\in \N_0^d:\\ 
\|\bj\|_1\leq q-\|\bll\|_1}} \frac{(\partial^{\bll+\bj} \hat{f})\left((\tilde{C}_{k,i})_{left} \right)}{\bj!} \cdot  \left((\tilde{C}_{k+1,i})_{left} - (\tilde{C}_{k,i})_{left}\right)^{\bj}\\
& \quad  - (\partial^{\bll} f)\left((C_{k+1,i})_{left}\right)\Bigg|\\
\leq & \Bigg|\sum_{\substack{\mathbf{j}\in \N_0^d:\\
\|\bj\|_1\leq q-\|\bll\|_1}} \frac{(\partial^{\bll+\bj} \hat{f})((\tilde{C}_{k,i})_{left}) }{\bj!} \cdot \left((\tilde{C}_{k+1,i})_{left} - (\tilde{C}_{k,i})_{left}\right)^{\bj}
\\
& - \sum_{\substack{\mathbf{j}\in \N_0^d:\\ 
\|\bj\|_1\leq q-\|\bll\|_1}} \frac{(\partial^{\bll+\bj} f)((\tilde{C}_{k,i})_{left})}{\bj!} \cdot \left((\tilde{C}_{k+1,i})_{left} - (\tilde{C}_{k,i})_{left}\right)^{\bj}\Bigg|
\\
&+ \Bigg| \sum_{\substack{\mathbf{j}\in \N_0^d:\\ 
\|\bj\|_1\leq q-\|\bll\|_1}} \frac{(\partial^{\bll+\bj} f)(\tilde{C}_{k,i})_{left}) }{\bj!} \cdot \left((\tilde{C}_{k+1,i})_{left} - (\tilde{C}_{k,i})_{left}\right)^{\bj}\\
& \quad - (\partial^{\bll} f)((\tilde{C}_{k+1,i})_{left})\Bigg|\\
\leq & \sum_{\substack{\mathbf{j}\in \N_0^d:\\ 
\|\bj\|_1\leq q-\|\bll\|_1}} \frac{1}{\bj!} \cdot c_{36} \cdot \left(\frac{2a}{M^2}\right)^{p-\|\bll+\bj\|_1} \cdot  \left(\frac{2a}{M^2}\right)^{\|\bj\|_1} + c_{36} \cdot \left(\frac{2a}{M^2}\right)^{p-\|\bll\|_1}\\
\leq & (c_{36} \cdot e^d + c_{36}) \cdot \left(\frac{2a}{M^2}\right)^{p-\|\bll\|_1}.
\end{align*}
This implies that we can choose $b_{k,i}^{(\bll)} \in \Z$ such that
\begin{align*}
|b_{k,i}^{(\bll)}| \leq e^d +1
\end{align*}
and
\begin{align*}
\left|(\partial^{\bll} \hat{f})((\tilde{C}_{k+1,i})_{left}) - (\partial^{\bll} f)((\tilde{C}_{k+1,i})_{left})\right| \leq c_{36} \cdot \left(\frac{2a}{M^2}\right)^{p-\|\bll\|_1}.
\end{align*}
Observe that in this way we have defined the coefficients $b_{k,i}^{(\bll)}$ for each cube $C_{i,1}$. We will encode these coefficients for each $i \in \{1, \dots, M^d\}$ and each $\bll \in \N_0^d$ with $\|\bll\|_1\leq q$ in the single number
\begin{align*}
b_i^{(\bll)} = \sum_{k=1}^{M^d-1}\left(b_{k,i}^{(\bll)} + \lceil e^d \rceil +2\right) \cdot (4+2 \lceil e^d \rceil)^{-k} \in [0,1].
\end{align*}

In the last step the neural network then computes
\begin{align}
\label{hatT}
\hat{T}_{f,q,(C_{\P_2}(x))_{left}}(x) := \sum_{\substack{\bll \in \N_0^d:\\ \|\bll\|_1 \leq q}} \frac{(\partial^{\bll} \hat{f})((C_{\P_2}(x))_{left})}{\bll!} \cdot \left(x - (C_{\P_2}(x))_{left}\right)^{\bll},
\end{align}
where by construction we have $C_{\P_2}(x) = \tilde{C}_{k,i}$ for some $k \in \{1, \dots, M^d\}$. Since 
\begin{align}
\label{2bhatT}
&\left|\hat{T}_{f,q,(C_{\P_2}(x))_{left}}(x)  - T_{f,q,(C_{\P_2}(x))_{left}}(x)\right|\notag\\
&\leq \sum_{\substack{\bll \in \N_0^d: \\ \|\bll\|_1 \leq q}} \frac{\left|(\partial^{\bll} \hat{f}-\partial^{\bll} f)((C_{\P_2}(x))_{left})\right|}{\bll!} \cdot \left|x - (C_{\P_2}(x))_{left}\right|^{\bll} \notag\\
& \leq e^d \cdot c_{36} \cdot \left(\frac{2a}{M^2}\right)^p
\end{align}
the network approximating $\hat{T}_{f,q,(C_{\P_2}(x))_{left}}(x)$ is also a good approximation for $T_{f,q,(C_{\P_2}(x))_{left}}(x)$. 
\\
\\
To approximate $f(x)$ by neural networks the proof of Kohler and Langer follows 
\textit{four} key steps:
\begin{enumerate}
\item[1.] Compute $\hat{T}_{f,q,(C_{\P_2}(x))_{left}}(x)$ by recursively defined functions. 
\item[2.] Approximate the recursive functions by neural networks. The resulting network is a good approximation for $f(x)$ in case that
\[
x \in \bigcup_{k \in \{1, \dots, M^{2d}\}} (C_{k,2})_{1/M^{2p+2}}^0.
\]
\item[3.] Approximate the function $w_{\P_2}(x) \cdot f(x)$ by deep neural networks, where
\begin{equation}
\label{w_vb}
w_{\P_2}(x) = \prod_{j=1}^d \left(1- \frac{M^2}{a} \cdot \left|(C_{\mathcal{P}_{2}}(x))_{left}^{(j)} + \frac{a}{M^2} - x^{(j)}\right|\right)_+
\end{equation}
is a linear tensor product B-spline
which takes its maximum value at the center of $C_{\P_{2}}(x)$, which
is nonzero in the inner part of $C_{\P_{2}}(x)$ and which
vanishes
outside of $C_{\P_{2}}(x)$. 
\item[4.] Apply those networks to $2^d$ slightly shifted partitions of $\P_2$ to approximate $f(x)$ in supremum norm. 
\end{enumerate}
We focus on step 2 and 3 and modify the construction of the auxiliary neural networks by deriving constraints for the required weights.
\subsubsection{Key step 1: A recursive definition of $\hat{T}_{f,q,(C_{\P_2}(x))_{left}}(x)$}
\label{subsubrecursiveTaylordef}
To derive a recursive definition of  $\hat{T}_{f,q,(C_{\P_2}(x))_{left}}(x)$, we set
\begin{align*}
&\bm{\phi}_{1,0} = \left(\phi_{1,0}^{(1)}, \dots, \phi_{1,0}^{(d)}\right) =x\\
&\bm{\phi}_{2,0} = \left(\phi_{2,0}^{(1)}, \dots, \phi_{2,0}^{(d)}\right) =\mathbf{0}
\end{align*}
and 
\begin{align*}
\phi_{3, 0}^{(\bll)}=0 \ \mbox{and} \ \phi_{4, 0}^{(\bll)}=0
\end{align*}
for each $\bll \in \N_0^d$ with $\|\bll\|_1 \leq q$. For $j \in \{1, \dots, M^d\}$ set
\begin{align*}
\bm{\phi}_{1, j} = \bm{\phi}_{1, j-1},
\end{align*}
\begin{align*}
\bm{\phi}_{2, j} = (C_{j,1})_{left} \cdot \mathds{1}_{C_{j,1}}(\bm{\phi}_{1, j-1}) + \bm{\phi}_{2, j-1},
\end{align*}
\begin{align*}
\phi_{3, j}^{(\bll)} = (\partial^{\bll} f)((C_{j,1})_{left}) \cdot \mathds{1}_{C_{j,1}}(\bm{\phi}_{1, j-1}) +\phi_{3, j-1}^{(\bll)}
\end{align*}
and
\begin{align*}
\phi_{4, j}^{(\bll)} = b_j^{(\bll)} \cdot \mathds{1}_{C_{j,1}}(\bm{\phi}_{1, j-1}) + \phi_{4, j-1}^{(\bll)}.
\end{align*}
Furthermore set
\begin{align*}
\bm{\phi}_{1, M^d+j} = \bm{\phi}_{1, M^d+j-1}, \quad j \in \{1, \dots, M^d\},
\end{align*}
\begin{align*}
\bm{\phi}_{2, M^d+j} = \bm{\phi}_{2, M^d+j-1}+\bold{\tilde{v}}_{j+1}, 
\end{align*}
\begin{eqnarray*}
  \phi_{3,M^d+j}^{(\bll)}=&&
  \sum_{\substack{\mathbf{s} \in \N_0^d\\ \|\mathbf{s}\|_1 \leq q-\|\bll\|_1}} \frac{\phi_{3,M^d+j-1}^{(\bll+\mathbf{s})}}{\mathbf{s}!} \cdot
  \left(\bold{\tilde{v}}_{j+1}\right)^{\mathbf{s}}
  \\
  &&
  +
  \left(
  \lfloor
  (4+2 \cdot \lceil e^d \rceil) \cdot \phi_{4,M^d+j-1}^{(\bll)}
  \rfloor
- \lceil e^d \rceil-2
  \right)
  \cdot c_{36} \cdot \left(
\frac{2a}{M^2}
  \right)^{p-\|\bll\|_1},
\end{eqnarray*}
\[
\phi_{4, M^d+j}^{(\bll)}=  (4+2 \cdot \lceil e^d\rceil) \cdot \phi_{4, M^d+j-1}^{(\bll)}
-
\lfloor
  (4+2 \cdot \lceil e^d \rceil) \cdot \phi_{4, M^d+j-1}^{(\bll)}
  \rfloor
\]
for $j \in \{1, \dots, M^d-1\}$ and each $\bll \in \N_0^d$ with $\|\bll\|_1 \leq q$ and
\begin{align*}
\bm{\phi}_{5, M^d+j} = \mathds{1}_{\mathcal{A}^{(j)}}(\bm{\phi}_{1, M^d+j-1}) \cdot \bm{\phi}_{2, M^d+j-1} + \bm{\phi}_{5, M^d+j-1}
\end{align*}
and 
\begin{align*}
\phi_{6, M^d+j}^{(\bll)} = \mathds{1}_{\mathcal{A}^{(j)}}(\phi_{1, M^d+j-1}) \cdot \phi_{3, M^d+j-1}^{(\bll)} + \phi_{6, M^d+j-1}^{(\bll)}
\end{align*}
for $j \in \{1, \dots, M^d\}$, where 
\begin{align*}
\bm{\phi}_{5, M^d} = \left(\phi_{5, M^d}^{(1)}, \dots, \phi_{5, M^d}^{(d)}\right)=\mathbf{0},  \ \phi_{6, M^d}^{(\bll)} =0
\end{align*}
and
\begin{align*}
&\mathcal{A}^{(j)}=\left\{x \in \R^d: -x^{(k)} + \phi_{2, M^d+j-1}^{(k)} \leq 0 \ \right. \notag\\
&
\quad \left.\mbox{und} \ x^{(k)} - \phi_{2, M^d+j-1}^{(k)}- \frac{2a}{M^2} < 0 \  \text{for all} \ k \in \{1, \dots, d\}\right\}.
\end{align*}
Finally define
\begin{align*}
\bm{\phi}_{1, 2M^d+1} = \sum_{\substack{\bll \in \N_0^d:\\ \|\bll\|_1\leq q}} &\frac{\phi_{6, 2M^d}^{(\bll)}}{\bll!} \cdot \left(\bm{\phi}_{1, 2M^d} - \bm{\phi}_{5, 2M^d}\right)^{\bll}.
\end{align*}
The next lemma shows that this recursion computes $\hat{T}_{f,q,(C_{\P_2}(x))_{left}}(x)$.
\begin{lemma}
\label{supple11}
  Let $p=q+s$ for some $q \in \N_0$ and $s \in (0,1]$, let $C > 0$ and $x \in [-a,a)^d$. Let $f: \Rd \to \R$ be a $(p,C)$-smooth function and let $\hat{T}_{f,q,(C_{\mathcal{P}_2}(x))_{left}}$ be defined as in \eqref{hatT}. Define $\bm{\phi}_{1, 2M^d+1}$ recursively as above. Then we have
  \[
\phi_{1, 2M^d+1}=\hat{T}_{f,q,(C_{\mathcal{P}_2}(x))_{left}}(x).
  \]
\end{lemma}
\noindent
{\bf Proof.}
    See Lemma 11 in Kohler and Langer (2021). \hfill
    $\Box$

\subsubsection{Key step 2: Approximating $\phi_{1, 2M^d+1}$ by neural networks}
In this step we show that a neural network approximates $\phi_{1, 2M^d+1}$ in case that
\begin{align*}
x \in \bigcup_{i \in \{1, \dots, M^{2d}\}} (C_{i,2})_{1/M^{2p+2}}^0.
\end{align*}
 
We define a composition neural network, which approximately computes the recursive functions in the definition of $\phi_{1, 2M^d+1}$. 
\begin{lemma}
\label{supple13}
Let $\sigma:\R \to \R$ be the ReLU activation function $\sigma(x) = \max\{x,0\}$. Let $\mathcal{P}_2$ be defined as in \eqref{partition}. Let $p = q+s$ for some $q \in \N_0$ and $s \in (0,1]$, and let $C >0$. Let $f: \Rd \to \R$ be a $(p,C)$-smooth function.
    Let $1 \leq a < \infty$. Then there exists for $M \in \N$ sufficiently large (independent of the size of $a$, but 
    \begin{eqnarray}
      \label{supple13eq1}
      M^{2p} &\geq& 2^{4(q+1)+1}
      \max\{ c_{37} \cdot (6+ 2 \lceil e^d \rceil)^{4(q+1)}, c_{36} \cdot e^d \}
      \nonumber     \\
      &&
      \hspace*{3cm}
      \cdot \left(\max\left\{a, \|f\|_{C^q([-a,a]^d)} \right\}\right)^{4(q+1)}
    \end{eqnarray}
     must hold), a neural network
$\hat{f}_{deep,\P_2} \in \mathcal{F}(L,r)$ with
\begin{itemize}
\item[(i)] $L= 4M^d+\left\lceil \log_4\left(M^{2p+4 \cdot d \cdot (q+1)} \cdot e^{4 \cdot (q+1) \cdot (M^d-1)}\right)\right\rceil\\
\hspace*{0.8cm} \cdot \lceil \log_2(\max\{q+1, 2\})\rceil$
\item[(ii)] $r= \max\left\{10d+4d^2+2 \cdot \binom{d+q}{d} \cdot \left(2 \cdot (4+2\lceil e^d\rceil)+5+2d\right), \right.\\
\left. \hspace*{1.5cm} 18 \cdot (q+1) \cdot \binom{d+q}{d}\right\}$
\end{itemize}
such that 
\begin{align*}
|\hat{f}_{deep, \mathcal{P}_2}(x) - f(x)|\leq c_{38} \cdot \left(\max\left\{2a, \|f\|_{C^q([-a,a]^d)}\right\}\right)^{4(q+1)} \cdot \frac{1}{M^{2p}}
\end{align*}
holds for all $x \in \bigcup_{i \in \{1, \dots, M^{2d}\}} \left(C_{i,2}\right)_{1/M^{2p+2}}^0$. The network value is bounded by 
\begin{align*}
|\hat{f}_{deep, \mathcal{P}_2}(x)| & \leq 1 + \Bigg(\|f\|_{C^q([-a,a]^d)} \cdot e^{(M^d-1)}\\
 & \quad + (4+2 \cdot \lceil e^d \rceil) \cdot (M^d-1)\cdot e^{(M^d-2)}\Bigg)\cdot e^{2ad}
\end{align*}
for all $x \in [-a,a)^d$.\\
$\hat{f}_{deep, \mathcal{P}_2}$ satisfies the weight constraint: 
$$
    \norm{\bv_{\hat{f}_{deep, \mathcal{P}_2}}}_\infty
    \leq  e^{c_{39}\cdot (M^d+d)\cdot 2(q+1)},
$$
where $c_{39}=c_{39}(f)$.
\end{lemma}

As in Kohler and Langer (2021), auxiliary networks are required to prove these results. 
We introduce the auxiliary networks with weight constraints and modify parts of the proofs accordingly. 

In the construction of our network we will compose smaller subnetworks to successively build the final network. Here instead of using an additional layer, we "merge" the weights of both networks $f$ and $g$ to define $f \circ g$. The following lemma clarifies this idea and derives appropriate weight constraints:
\begin{lemma}
\label{lemma:composition-weight-bound}
 Let $f_0:\mathbb{R}^k\to \mathbb{R}$ be a neural network of the class $\mathcal{F}(L, r)$ with weight vector $\mathbf{v}_0$ and let $f_1, \ldots, f_k:\mathbb{R}^d\to \mathbb{R}$ be neural networks of class $\mathcal{F}(\bar{L},\bar{r})$ with weight vectors $\bv_1, \ldots, \bv_k.$ Denote by $\bar{\bv}$ the vector that contains $(\bv_j)_{j\in\{1, \ldots, k\}}$. \\
Then the network $f= f_0 (f_1, \ldots, f_k)$ has $L+\bar{L}$ layers and at most $\max\{k\cdot \bar{r}, r\}$ neurons.

\begin{enumerate}
\item [a)] In general $\bv$ satisfies the constraints
$$
    \norm{\bv}_\infty 
\leq 
\max\left\{
\norm{\bv_0}_\infty, 
\norm{\bar{\bv}}_\infty, 
\norm{(\bv_0)^{(0)}}_\infty
\cdot
\left(k\norm{(\bar{\bv})^{(\bar{L})}}_\infty+1\right)
\right\}.
$$
    \item [b)] If $(\bv_j)_{1,0}^{(\bar{L})}=0$ for all $j\in \{1, \ldots, k\}$, $\bv$ satisfies :
    $$
    \norm{\bv}_\infty
    \leq 
    \max\left\{\norm{\bv_0}_\infty, 
    \norm{\bar{\bv}}_\infty,
    \norm{(\bv_0)^{(0)}_{i,j>0}}_\infty \cdot
    \norm{(\bar{\bv})^{(\bar{L})}_{1,j>0}}_\infty\right\}.
    $$
    \item [c)] If additionally to $(\bv_j)_{1,0}^{(\bar{L})}=0$ for all $j\in \{1, \ldots, k\}$, we have $\norm{(\bv_0)^{(0)}_{i,j>0}}_\infty \leq 1$ or $
    \norm{(\bar{\bv})^{(\bar{L})}_{i,j>0}}_\infty\le 1$, then
    $\bv$ satisfies
    $$
    \norm{\bv}_\infty
    \leq 
    \max\left\{\norm{\bv_0}_\infty, 
    \norm{\bar{\bv}}_\infty
    \right\}.
    $$
\end{enumerate}
\end{lemma}

\noindent
{\bf Proof.}
The network $f= f_0 (f_1, \ldots, f_k)$ is recursively defined as follows
$$
f(x)=\sum_{i=1}^{r}\left(\bv_0\right)_{1,i}^{(L)}f_i^{(\bar{L}+L)}(x)+\left(\bv_0\right)^{(L)}_{1,0},
$$
where for $l\in\{2,\ldots, L\}$ the outputs of the neurons $f_i^{(\bar{L}+l)}$ are recursively defined by
\[f_i^{(\bar{L}+l)}(x)=\sigma\left(\sum_{j=1}^{r}\left(\bv_0\right)_{i,j}^{(l-1)}f_j^{(\bar{L}+l-1)}(x)+\left(\bv_0\right)_{i,0}^{(l-1)}\right)\]
for $i\in \{1,\ldots, r\}$.
In layer $\bar{L}$, the effect of "merging" the networks together becomes apparent and we can see that layer $\bar{L}$ the network $f=f_0(f_1,\ldots, f_k)$ consists of $k\cdot \bar{r}$ neurons:
 \begin{equation}
     \label{compositon-layer}
 \begin{aligned}
  f_i^{(\bar{L}+1)}(x) &= \sigma\left(\sum_{j=1}^k (\bv_0)^{(0)}_{i,j} f_j(x)+ (\bv)^{(0)}_{i,0}\right)\\ 
  &= \sigma\left(\sum_{j=1}^k (\bv_0)^{(0)}_{i,j} \cdot\left(\sum_{l=1}^{\bar{r}} (\bv_j)^{(\bar{L})}_{1,l} \cdot f^{(\bar{L})}_{j,l}(x)
  + (\bv_j)_{1,0}^{(\bar{L})}\right)
  + (\bv_0)^{(0)}_{i,0}\right)\\
  &=\sigma\left(\sum_{j=1}^k \sum_{l=1}^{\bar{r}} (\bv_0)^{(0)}_{i,j} \cdot(\bv_j)^{(\bar{L})}_{1,l}\cdot f^{(\bar{L})}_{j,l}(x)
  + \sum_{j=1}^k (\bv_0)^{(0)}_{i,j} \cdot(\bv_j)_{1,0}^{(\bar{L})}+ (\bv_0)^{(0)}_{i,0}\right)
 \end{aligned}
 \end{equation}
 for $i\in\{1, \ldots, r\}$, where the $f^{(s)}_{j,i}(x)$ are defined by
 $$
 f^{(s)}_{(j-1)\cdot i+i}(x)= f^{(s)}_{j,i}(x) = \sigma
 \left(\sum_{l=1}^{\bar r}\left(\bv_j\right)_{i,l}^{(s-1)}f_{j,l}^{(s-1)}(x)
 +\left(\bv_j\right)_{i,0}^{(s-1)}\right)\
 $$
 for $j\in\{1,\ldots,k\}, s\in \{1, \ldots,\bar L\}$ and $i\in\{1,\ldots,\bar r\}$.
 Finally we have 
 \[
f^{(s)}_{(j-1)\cdot l+l}(x)=f_{j, l}^{(0)}(x)=x^{(l)}.
\]
for $j\in\{1,\ldots, k\}, l\in \{1,\ldots,d\}$.

In layers $l\in\{1,\ldots \bar{L}-1\}$, the weights of $f$ satisfy the same constraints as $f_1, \ldots, f_k$, in layers $l\in \{\bar{L}+1, \ldots, L+\bar{L}\}$ the weight constraints of $f$ correspond to the constraints of $f_0$. However, in layer $\bar{L}$ we have to consider the product of the output weights of $f_i$ for $i\in\{1,\ldots, k\}$ and the weights of the input layer of $f_0$, as shown in \eqref{compositon-layer}.
The weights there satisfy
\begin{eqnarray*}
  &&
  \lvert
(\bv_0)^{(0)}_{i,j} \cdot(\bv_j)^{(\bar{L})}_{1,l}
  \rvert
 \leq 
 \max_{j\in \{1,\ldots, k\}, l\in \{1,\ldots, r\}}
 \left|(\bv_0)^{(0)}_{i,j}\right|\cdot\left|(\bv_j)^{(\bar{L})}_{1,l}\right|\leq 
 \norm{(\bv_0)^{(0)}}_\infty \cdot
    \norm{(\bar{\bv})^{(\bar{L})}}_\infty
 \end{eqnarray*}
 and
 $$
 \left|\sum_{j=1}^k (\bv_0)^{(0)}_{i,j} \cdot(\bv_j)_{1,0}^{(\bar{L})}+ (\bv_0)^{(0)}_{i,0}\right|
 \leq 
 \norm{(\bv_0)^{(0)}}_\infty
\cdot
\left(\sum_{j=1}^k\norm{(\bv_j)^{(\bar{L})}}_\infty+1\right),
 $$
 which implies part a) of the assertion.\\
 If, additionally, $(\bv_j)_{1,0}^{(\bar{L})}=0$ for all $j\in \{1, \ldots, k\}$, the latter bound can be refined to
  $$
 \lvert(\bv)\rvert^{(\bar{L})}_{i,0}
 =
 \left|(\bv_0)^{(0)}_{i,0}\right|
 \leq 
 \norm{(\bv_0)^{(0)}}_\infty .
 $$
 Since $\norm{(\bv_0)^{(0)}}_\infty \leq \norm{(\bv_0)}_\infty$, this case no longer influences the upper bound and we can reduce the third argument of the maximum to $\norm{(\bv_0)^{(0)}_{i,j>0}}_\infty \cdot
    \norm{(\bar{\bv})^{(\bar{L})}_{i,j>0}}_\infty$, which implies the upper bound of part b) of the assertion.\\
Part c) follows directly from part b).
\hfill $\Box$

The following lemma presents a neural network, that approximates the square function, which is essential to build neural networks for more complex tasks.

\begin{lemma}
\label{square-fct}
Let $\sigma: \R \to \R$ be the ReLU activation function $\sigma(x) = \max\{x,0\}$. Then for any $R \in \N$ and any $a \geq 1$ a neural network
  \begin{equation*}
  \hat{f}_{sq} \in \mathcal{F}(R,9)
  \end{equation*}
  exists with weight constraints
$$
\norm{\bv_{\hat{f}_{sq}}}_\infty 
=
  \norm{(\bv_{\hat{f}_{sq}})^{(R)}}_\infty 
  \leq 
  4\cdot a^2
  \qquad\text{and} \qquad
  \norm{(\bv_{\hat{f}_{sq}})^{(0)}}_\infty 
  \leq 
  1,
  $$
  such that
  \begin{equation*}
  \left|\hat{f}_{sq}(x) - x^2\right| \leq a^2 \cdot 4^{-R}
  \end{equation*}
  holds for  $x\in [-a,a]$.
\end{lemma}

\noindent
{\bf Proof.}
This proof follows as in Kohler and Langer (2021). We only modify the parts required to show the constraints on the weights.
In Kohler and Langer (2021) it was shown that linear combinations of the "tooth" function $g: [0,1] \to [0,1]$
\begin{equation*}
g(x) = 
\begin{cases}
2x &\quad, x \leq \frac{1}{2}\\
2 \cdot (1-x) &\quad, x > \frac{1}{2}
\end{cases}
\end{equation*}
and the iterated composition function 
\begin{equation*}
g_s(x) = \underbrace{g \circ g \circ \dots \circ g}_{s}(x).
\end{equation*}
can be used to approximate $f(x) = x^2$ for  $x \in [0,1]$. Let $S_R$ denote the piecewise linear interpolation of $f$ with $2^R+1$ uniformly distributed breakpoints. 

It can be shown that $S_R(x)$ is given by
\begin{equation*}
S_R(x) = x- \sum_{s=1}^R \frac{g_s(x)}{2^{2s}}
\end{equation*}
and that it satisfies
\begin{equation*}
|S_R(x) - x^2| \leq 2^{-2R-2}
\end{equation*}
for $x \in [0,1]$.
\newline
In a \textit{third step of their proof} Kohler and Langer show, that there exists a feedforward neural network that computes $S_R(x)$ for $x \in [0,1]$. 
In order to derive the weight constraints we include the construction of this network.\\
The function $g(x)$ can be implemented by the network:
\begin{equation*}
\hat{f}_g(x) = 2\cdot \sigma(x) - 4 \cdot \sigma(x-\frac{1}{2}) + 2\cdot \sigma(x-1)
\end{equation*}
and the function $g_s(x)$ can be implemented by a network
\begin{equation*}
\hat{f}_{g_s}  \in \mathcal{F}(s, 3)
\end{equation*}
with
\begin{equation*}
\hat{f}_{g_s}(x) = \underbrace{\hat{f}_g(\hat{f}_g(\dots(\hat{f}_g}_s(x))).
\end{equation*}
Thus we have 
$$
\norm{\bv_{\hat{f}_g}}_\infty = \norm{\bv_{\hat{f}_{g_s}}}_\infty= 4\quad \text{ for all } s\in \mathbb{N}.
$$
Let 
\begin{equation*}
\hat{f}_{id}(z) = \sigma(z) - \sigma(-z),
\end{equation*}
and define $\hat{f}^t_{id}$ recursively by
\begin{align*}
\hat{f}_{id}^0(z) &= z \quad &(z \in \R)\\
\hat{f}_{id}^{t+1}(z) &= \hat{f}_{id}(\hat{f}_{id}^t(z)) \quad &(z \in \R, t \in \N_0),
\end{align*}
which implies
\begin{equation*}
\hat{f}_{id}^t(z)=z.
\end{equation*}
It is easy to see that these networks satisfy 
$$\norm{\bv_{\hat{f}_{id}}}_\infty= \norm{\bv_{\hat{f}^t_{id}}}_\infty= 1 \quad \text{ for all } t\in \mathbb{N}.$$
By combining the networks above we can implement the function $S_R(x)$ by a network
\begin{equation*}
\hat{f}_{sq_{[0,1]}} \in \mathcal{F}(R,7)
\end{equation*}
recursively defined as follows:
We set
$\hat{f}_{1,0}(x)=\hat{f}_{2,0}(x)=x$ and $\hat{f}_{3,0}(x)=0$.

Then we set
\[
\hat{f}_{1, i+1}(x)=\hat{f}_{id}(\hat{f}_{1,i}(x)),
\]
\[
\hat{f}_{2, i+1}(x)={\hat{f}_g(\hat{f}_{2, i}(x))}
\]
and
\[
\hat{f}_{3,i+1}(x)=f_{id}(\hat{f}_{3,i}(x))-\frac{\hat{f}_g(\hat{f}_{2,i}(x))}{2^{2(i+1)}}
\]
for $i \in \{0,1, \dots,R-2\}$ and
\[
\hat{f}_{sq_{[0,1]}}(x)
=
\hat{f}_{id}(\hat{f}_{1, R-1}(x)) -  \frac{\hat{f}_{2, R}(x)}{2^{2R}} + \hat{f}_{id}(\hat{f}_{3, R-1}(x)).
\]

Using the positive homogeneity of the ReLU function, this implies
\begin{align*}
\hat{f}_{sq_{[0,1]}}(x) 
= & \hat{f}_{i d}^R(x)-\frac{1}{2^{2 R}} \hat{f}_{g_R}(x)-\hat{f}_{i d}\left(\frac{1}{2^{2(R-1)}} \hat{f}_{g_{R-1}}(x)\right. \\ & \left.-\hat{f}_{i d}\left(\frac{1}{2^{2(R-2)}} \hat{f}_{g_{R-2}}(x)-\cdots-\hat{f}_{i d}\left(\frac{1}{2^2} \hat{f}_{g_1}(x)\right)\right)\right) \\
=& S_R(x),
\end{align*}
hence $\hat{f}_{sq_{[0,1]}}(x)$ satisfies
\begin{equation}\label{le1aeq1}
|\hat{f}_{sq_{[0,1]}}(x) - x^2| \leq 2^{-2R-2}
\end{equation}
for $x \in [0,1]$.
We can construct the network $\hat{f}_{sq_{[0,1]}}$ such that its weights satisfy the constraint
$$
\norm{\bv_{\hat{f}_{sq_{[0,1]}}}}_\infty = \frac{1}{4}\cdot\norm{\bv_{\hat{f}_{g}}}_\infty
= 1.
$$
To show this, we need to show that the absolute value of each weight of the network $\tilde{f}_{2,i}(x):= \frac{\hat{f}_{2,i}(x)}{2^{2i}}$ is smaller or equal to one for every $i$. 

We proceed by induction. For $i=1 $ we have 
$$\frac{\tilde{f}_{2,1}(x)}{2^2} = \frac{1}{4}\hat{f}_g(x)= \frac{1}{2}\sigma(x)-\sigma\left(x-\frac{1}{2}\right)+\frac{1}{2}\sigma(x-1)$$
and thus $\norm{\bv_{\tilde{f}_{2,1}}}_\infty =1$. \\
Assume that $\norm{\bv_{\tilde{f}_{2,i}}}_\infty =1$ holds for arbitrary but fixed $i>0$. \\
Using the positive homogeneity of the ReLU function we can conclude 
$$
\begin{aligned}
    \tilde{f}_{2,i+1}(x)= \frac{\hat{f}_{2,i+1}(x)}{2^{2(i+1)}}&= \frac{1}{4}\cdot \frac{1}{2^{2i}}\hat{f}_g(\hat{f}_{2,1}(x))\\
    &= \frac{1}{2^{2i}} \cdot\left(\frac{1}{2}\sigma\left(\hat{f}_{2,i}(x)\right)-\sigma\left(\hat{f}_{2,i}(x)-\frac{1}{2}\right)+\frac{1}{2}\sigma\left(\hat{f}_{2,i}(x)-1\right)\right)\\
    &=\frac{1}{2}\sigma\left(\frac{1}{2^{2i}}\cdot\hat{f}_{2,i}(x) \right)-\sigma\left(\frac{1}{2^{2i}}\cdot\hat{f}_{2,i}(x)-\frac{1}{2^{2i}}\cdot\frac{1}{2}\right)+\frac{1}{2}\sigma\left(\frac{1}{2^{2i}}\cdot\hat{f}_{2,i}(x)-\frac{1}{2^{2i}}\right)\\
     &=\frac{1}{2}\sigma\left(\tilde{f}_{2,i}(x) \right)-\sigma\left(\tilde{f}_{2,i}(x)-\frac{1}{2^{2i}}\cdot\frac{1}{2}\right)+\frac{1}{2}\sigma\left(\tilde{f}_{2,i}(x)-\frac{1}{2^{2i}}\right).
\end{aligned}$$
and thus $\norm{\bv_{\tilde{f}_{2,i+1}}}_\infty =1$.\\
In a \textit{last step} $\hat{f}_{sq_{[0,1]}}$ is extended to approximate the function $f(x) = x^2$ on the domain $ [-a,a]$. Therefore $f_{tran}: [-a,a] \to [0,1]$ is defined by
\begin{equation*}
f_{tran}(z) = \frac{z}{2a}+\frac{1}{2}
\end{equation*}
to transfers the value of $x \in [-a,a]$ in the interval, where \eqref{le1aeq1} holds. Set 
\begin{equation*}
\hat{f}_{sq}(x) = 4a^2\hat{f}_{sq_{[0,1]}}(f_{tran}(x)) - (2a \cdot \hat{f}_{id}^R(x) +a^2)
\end{equation*}
The extension to the domain $[-a,a]$ only increases the weights of the last layer of the network, which results in the constraints 
$$
\norm{(\bv_{\hat{f}_{sq}})}_\infty 
=
  \norm{(\bv_{\hat{f}_{sq}})^{(R)}}_\infty 
  \leq 
  4\cdot a^2 
  \qquad
  \text{and}
  \qquad
  \norm{(\bv_{\hat{f}_{sq_{[0,1]}}})^{(0)}}_\infty= \norm{(\bv_{\hat{f}_{sq}})^{(0)}}_\infty
  \leq 1.
  $$
Since 
\begin{equation*}
x^2 = 4a^2 \cdot \left(\frac{x}{2a}+\frac{1}{2}\right)^2 - 2ax -a^2
\end{equation*}
we have
\begin{align*}
&|\hat{f}_{sq}(x) - x^2|\\
=& |4a^2\hat{f}_{sq_{[0,1]}}(f_{tran}(x)) - (2a \cdot \hat{f}_{id}^R(x) +a^2) - (4a^2 \cdot (f_{tran}(x))^2 - 2ax -a^2)|\\
\leq & 4a^2 \cdot |\hat{f}_{sq_{[0,1]}}(f_{tran}(x)) - (f_{tran}(x))^2| + 2a|\hat{f}_{id}^R(x) - x|\\
\leq & 4a^2 \cdot 2^{-2R-2} = a^2 \cdot 4^{-R}.
\end{align*}
  
\hfill $\Box$
\begin{lemma}
\label{lemma:fmult}
Let $\sigma: \mathbb{R} \to \R$ be the ReLU activation function $\sigma(x) = \max\{x,0\}$. Then for any $R \in \N$ and any $a \geq 1$ a neural network
\begin{equation*}
\hat{f}_{mult} \in \mathcal{F}(R,18)
\end{equation*}
 exists, whose weights satisfy
 $$
\norm{\bv_{\hat{f}_{mult}}}_\infty \leq 4\cdot a^2, \qquad  \left(\bv_{\hat{f}_{mult}}\right)_{1,0}^{(R)}=0\qquad \text{and}$$
 $$\norm{\left(\bv_{\hat{f}_{mult}}\right)^{(0)}}_\infty \leq 1,
 $$
  such that
\begin{equation*}
|\hat{f}_{mult}(x,y) - x \cdot y| \leq 2 \cdot a^2 \cdot 4^{-R}
\end{equation*}
holds for all $x, y \in [-a, a]$.
\end{lemma}

\noindent
{\bf Proof.}
Let 
\begin{equation*}
\hat{f}_{sq} \in \mathcal{F}(R, 9)
\end{equation*}
be the neural network from Lemma \ref{square-fct}, i.e.
\begin{equation*}
\hat{f}_{sq}(x) = 16a^2\hat{f}_{sq_{[0,1]}}\left(\frac{x}{4a}+\frac{1}{2}\right) - (4a \cdot \hat{f}_{id}^R(x) +4a^2).
\end{equation*}
which satisfies
\begin{equation*}
|\hat{f}_{sq}(x) - x^2| 
\leq 
4 \cdot a^2 \cdot 4^{-R}
\end{equation*}
for $x \in [-2a, 2a]$ and with weight constraints 
$\norm{\bv_{\hat{f}_{sq_{[0,1]}}}}_\infty = 
\norm{\bv_{\hat{f}_{id}^R}}_\infty = 1,$ and set
\[
\begin{aligned}
  \hat{f}_{mult}(x,y)
=&
\frac{1}{4} \cdot \left(\hat{f}_{sq}(x+y)-\hat{f}_{sq}(x-y) \right) \\
=&4a^2\hat{f}_{sq_{[0,1]}}\left(\frac{x+y}{4a}+\frac{1}{2}\right) - a \cdot \hat{f}_{id}^R(x+y)\\
&-\left(
4a^2\hat{f}_{sq_{[0,1]}}\left(\frac{x-y}{4a}+\frac{1}{2}\right) - a \cdot \hat{f}_{id}^R(x-y)\right).
\end{aligned}
\]
Note that since $\hat{f}_{sq}(x+y)$ and $\hat{f}_{sq}(x-y)$ have the same offset in the last layer, they cancel out and we have $\left(\bv_{\hat{f}_{mult}}\right)_{1,0}^{(R)}=0$.
The constraints $\norm{\bv_{\hat{f}_{mult}}}_\infty = 4\cdot a^2$ and $\norm{\left(\bv_{\hat{f}_{mult}}\right)^{(0)}}_\infty \leq 1
 $ follow directly from Lemma \ref{square-fct}.
Since 
\begin{equation*}
x \cdot y = \frac{1}{4} \left((x+y)^2 - (x-y)^2\right)
\end{equation*}
we have
\begin{align*}
\begin{split}
|\hat{f}_{mult}(x,y) - x \cdot y|
&\leq \frac{1}{4} \cdot \left|\hat{f}_{sq}(x+y) - (x+y)^2\right| + \frac{1}{4} \cdot \left|(x-y)^2 - \hat{f}_{sq}(x-y)\right|\\
&\leq \frac{1}{4} \cdot 2 \cdot 4 \cdot a^2 \cdot 4^{-R}\\
&\leq 2 \cdot a^2 \cdot 4^{-R}
\end{split}
\end{align*}
for $x,y \in [-a, a]$. 
\hfill $\Box$
\begin{lemma}
\label{lemma:fmultd}
Let $\sigma: \R \to \R$ be the ReLU activation function $\sigma(x) = \max\{x,0\}$. Then for $R \in \N, R \geq \log_4 \left(2 \cdot 4^{2 \cdot d} \cdot a^{2 \cdot d}\right)$ and any $a \geq 1$ a neural network 
\begin{align*}
\hat{f}_{mult, d} \in \mathcal{F}(R \cdot \lceil \log_2(d) \rceil, 18d),
\end{align*}
which satisfies 
$$
\norm{\bv_{\hat{f}_{mult,d}}}_\infty\leq 4\cdot 4^{2d}\cdot a^{2d},
\quad \left(\bv_{\hat{f}_{mult,d}}\right)_{1,0}^{(R \cdot \lceil \log_2(d) \rceil)} =0
\quad 
\text{and}$$
$$
\norm{\left(\bv_{\hat{f}_{mult,d}}\right)^{(0)}}_\infty\leq 1,
$$
exists such that
\begin{align*}
\left|
\hat{f}_{mult, d}(x) - \prod_{i=1}^dx^{(i)}
\right| 
\leq 4^{4d+1}\cdot a^{4d} \cdot d \cdot 4^{-R}
\end{align*}
holds for all $x \in [-a,a]^d$. 
\end{lemma}

\noindent
{\bf Proof.}
We set $q=\lceil \log_2(d)\rceil$. The feedforward neural network $\hat{f}_{mult, d}$ with $L=R \cdot q$ hidden layers and $r=18d$ neurons in each layer is constructed as follows: Set 
\begin{equation}
\label{neq1}
(z_1, \dots, z_{2^q})=
  \left(x^{(1)}, x^{(2)}, \dots, x^{(d)}, \underbrace{1, \dots,1}_{2^q-d} \right).
\end{equation}
In the construction of our network we will  use the network $\hat{f}_{mult}$ of Lemma \ref{lemma:fmult},  
which satisfies
\begin{equation}
  \label{ple4eq1}
|\hat{f}_{mult}(x,y) - x \cdot y| \leq 2\cdot (4^{d} a^{d})^2 \cdot 4^{-R}
\end{equation}
and  
 $$
\norm{\bv_{\hat{f}_{mult}}}_\infty \leq 4\cdot 4^{2d} a^{2d}, 
\quad  
\left(\bv_{\hat{f}_{mult}}\right)_{1,0}^{(R)}=0
\quad \text{and}\quad 
\norm{\left(\bv_{\hat{f}_{mult}}\right)^{(0)}}_\infty \leq 1
 $$
for $x,y \in [-4^{d} a^{d},4^{d} a^{d}]$. In the first $R$ layers we compute
\[
\hat{f}_{mult}(z_1,z_2), 
\hat{f}_{mult}(z_3,z_4), 
\dots,
\hat{f}_{mult}(z_{2^q-1},z_{2^q}), 
\]
which can be done by $R$ layers of $18 \cdot 2^{q-1} \leq 18 \cdot d$
neurons. E.g., in case
in case $z_l=x^{(d)}$ and $z_{l+1}=1$ we have 
\[
\hat{f}_{mult}(z_l,z_{l+1})
=\hat{f}_{mult}(x^{(d)},1).
\]
As a result of the first $R$ layers we get a vector of outputs
which has length $2^{q-1}$. Next we pair these outputs and apply $\hat{f}_{mult}$ again. This procedure is continued until there is only one output left.
Therefore we need $L =R q$ hidden layers and
at most $18d$
neurons in each layer. 
\newline
\newline
By  (\ref{ple4eq1}) 
and $R \geq \log_4 \left(2 \cdot 4^{2 \cdot d} \cdot a^{2 \cdot d}\right)$
we get for any $l \in \{1,\dots,d\}$ and any
$z_1,z_2 \in [-(4^l-1) \cdot a^l,(4^l-1) \cdot a^l]$
\[
|\hat{f}_{mult}(z_1,z_2)| \leq
|z_1 \cdot z_2| + |\hat{f}_{mult}(z_1,z_2)-z_1 \cdot z_2|
\leq
(4^l-1)^2 a^{2l} + 1
\leq
(4^{2l}-1) \cdot a^{2l}.
\]
From this we get successively that all outputs
of 
layer $l \in \{1,\dots,q-1\}$
  are contained in the interval
$[-(4^{2^l}-1) \cdot a^{2^l},(4^{2^l}-1) \cdot a^{2^l}]$, hence in particular they
  are contained in the interval
$[-4^{d} a^{d},4^{d} a^{d}]$
where inequality  (\ref{ple4eq1}) does hold.
\newline
\newline
Define $\hat{f}_{2^q}$ recursively by
\[
\hat{f}_{2^q}(z_1,\dots,z_{2^q})
=
\hat{f}_{mult}(\hat{f}_{2^{q-1}}(z_1,\dots,z_{2^{q-1}}),\hat{f}_{2^{q-1}}(z_{2^{q-1}+1},\dots,z_{2^q}))
\]
and
\[
\hat{f}_2(z_1,z_{2})= \hat{f}_{mult}(z_1,z_{2}).
\]
The constraints $\norm{\left(\bv_{\hat{f}_{mult,d}}\right)^{(0)}}_\infty \leq 1
 $ and $\left(\bv_{\hat{f}_{mult}}\right)_{1,0}^{(R)}=0$ follow directly from Lemma \ref{lemma:fmult}, since $\hat{f}_{mult,d}$ is a repeated composition of $\hat{f}_{mult}$. Note that our construction of $\hat{f}_{mult}$ satisfies the special case of Lemma \ref{lemma:composition-weight-bound}, i.e. 
$\left(\bv_{\hat{f}_{mult}}\right)_{1,0}^{(R)}=0$ and further $\norm{\left(\bv_{\hat{f}_{mult}}\right)^{(0)}}_\infty \leq 1$.

Applying Lemma \ref{lemma:composition-weight-bound} b) we get for the repeated composition of $\hat{f}_{mult}$: 
$$\norm{\bv_{\hat{f}_{mult,d}}}_\infty\leq\norm{\left(\bv_{\hat{f}_{mult}}\right)^{(0)}}_\infty\cdot \norm{\left(\bv_{\hat{f}_{mult}}\right)^{(R)}}_\infty \leq  4\cdot 4^{2d}\cdot a^{2d}.$$

The rest of the proof follows analogously to the proof of Lemma 8 in Kohler and Langer (2021).

\hfill $\Box$

\begin{lemma}
\label{le:polynomial-approx}
Let $m_1, \dots, m_{\binom{d+N}{d}}$ denote all monomials in $\mathcal{P}_N$ for some $N \in \N$. Let $r_1, \dots, r_{\binom{d+N}{d}} \in \R$, define 
\begin{align*}
p\left(x, y_1, \dots, y_{\binom{d+N}{d}}\right) = \sum_{i=1}^{\binom{d+N}{d}} r_i \cdot y_i \cdot m_i(x), \quad x \in [-a,a]^d, y_i \in [-a,a]
\end{align*}
and set $\bar{r}(p) = \max_{i \in \left\{1, \dots, \binom{d+N}{d}\right\}} |r_i|$. Let $\sigma: \mathbb{R} \to \R$ be the ReLU activation function $\sigma(x) = \max\{x,0\}$. Then for any $a \geq 1$ and 
\begin{equation}
\label{le_appendix3eq1}
R \geq \log_4 (2 \cdot 4^{2 \cdot (N+1)} \cdot a^{2 \cdot (N+1)})
\end{equation}
 a neural network 
$
\hat{f}_{p} \in \mathcal{F}(L,r)$
with $L=R \cdot \lceil\log_2(N+1)\rceil$ and $r =18 \cdot (N+1) \cdot \binom{d+N}{d}$ exists, whose weights satisfy 
$$\norm{\left(\bv_{\hat{f}_{p}}\right)^{(0)}}_\infty \leq 1
 ,\quad 
 \left(\bv_{\hat{f}_{p}}\right)_{1,0}^{(R)}=0
 \quad\text{and}$$
$$\norm{\bv_{\hat{f}_p}}_\infty\leq 4\cdot \bar{r}(p)\cdot 4^{2(N+1)}\cdot a^{2(N+1)},$$ such that
\begin{align*}
\left|\hat{f}_{p}\left(x, y_1, \dots, y_{\binom{d+N}{d}}\right) - p\left(x, y_1, \dots, y_{\binom{d+N}{d}}\right) \right| \leq c_{40} \cdot \bar{r}(p) \cdot a^{4(N+1)} \cdot 4^{-R}
\end{align*}
for all $x \in [-a,a]^d$, $y_1, \dots, y_{\binom{d+N}{d}} \in [-a,a]$, where $c_{40}$ depends on $d$ and $N$.
\end{lemma}

\noindent
{\bf Proof.}
A neural network $\hat{f}_m$ is constructed in order to approximate 
\begin{align*}
y \cdot m(x) = y \cdot \prod_{k=1}^d \left(x^{(k)}\right)^{r_k}, \quad x \in [-a,a]^d, y \in [-a,a],
\end{align*}
where $m \in \mathcal{P}_N$ and $r_1, \dots, r_d \in \N_0$ with $r_1+\dots+r_d \leq N$. Note that Lemma \ref{lemma:fmultd} can easily be extended to monomials. We set $d$ by $N+1$ and thus get a network 
\begin{align*}
\hat{f}_{m} \in \mathcal{F}(R \cdot \lceil \log_2(N+1)\rceil, 18 \cdot (N+1))
\end{align*} whose weights satisfy 
$$\norm{\bv_{\hat{f}_m}}_\infty \leq 4\cdot 4^{2(N+1)}\cdot a^{2(N+1)}.$$

We then set $\hat{f}_p = \sum_{i=1}^{\binom{d+N}{d}} r_i \cdot \hat{f}_{m_i}(x,y_i)$, which increases the weight constraint by a factor $\bar{r}(p)$. The weight constraints
$\norm{\left(\bv_{\hat{f}_{p}}\right)^{(0)}}_\infty \leq 1
 $ and $\left(\bv_{\hat{f}_{p}}\right)_{1,0}^{(R)}=0$ follow directly from Lemma \ref{lemma:fmultd}. The rest of the proof follows as in the proof of Lemma 5 of Kohler and Langer (2021).
\hfill $\Box$
\begin{lemma}
\label{le:f_ind_f_test}
Let $\sigma: \R \to \R$ be the ReLU activation function $\sigma(x) = \max\{x,0\}$. Let $R \in \N$. Let $\mathbf{a}, \mathbf{b} \in \Rd$ with
\begin{align*}
b^{(i)} - a^{(i)} \geq \frac{2}{R} \ \mbox{for all} \ i \in \{1, \dots, d\}
\end{align*}
and let
\begin{align*}
&K_{1/R} = \big\{x \in \Rd: x^{(i)} \notin [a^{(i)}, a^{(i)}+1/R) \cup (b^{(i)} - 1/R, b^{(i)})\\
& \hspace*{8cm}  \mbox{for all} \ i \in \{1, \dots, d\}\big\}.
\end{align*}
a) Then the network
\begin{align*}
\hat{f}_{ind, [\bold{a}, \bold{b})}(x) &= \sigma\bigg(1-R \cdot \sum_{i=1}^d \left(\sigma\left(a^{(i)} + \frac{1}{R} - x^{(i)}\right)\right.\\
& \hspace*{3cm} \left. + \sigma\left(x^{(i)} - b^{(i)} + \frac{1}{R}\right) \right)\bigg)
\end{align*}
of the class $\mathcal{F}(2, 2d)$ satisfies
the weight constraint 
$$
\norm{\bv_{\hat{f}_{ind, [\bold{a}, \bold{b})}}}_\infty
\leq \max\left\{\norm{\bold{a}}_\infty+\frac{1}{R}, \norm{\bold{b}}_\infty+\frac{1}{R}, R\right\},$$
as well as 
$$
\norm{\left(\bv_{\hat{f}_{ind, [\bold{a}, \bold{b})}}\right)^{(0)}_{i,j>0}}_\infty=1, 
\quad
\left(\bv_{\hat{f}_{ind, [\bold{a}, \bold{b})}}\right)^{(2)}_{1,0}=0
\quad \text{ and } \quad
\norm{\left(\bv_{\hat{f}_{ind, [\bold{a}, \bold{b})}}\right)^{(2)}_{1,i>0}}_\infty=1.
$$
For $x \in K_{1/R}$ we have
\begin{align*}
\hat{f}_{ind, [\bold{a}, \bold{b})}(x) = \mathds{1}_{ [\bold{a}, \bold{b})}(x)
\end{align*}
and 
\begin{align*}
\left|\hat{f}_{ind, [\bold{a}, \bold{b})}(x) - \mathds{1}_{[\bold{a}, \bold{b})}(\mathbf{x})\right| \leq 1
\end{align*}
for $x \in \Rd$.
\\
b) Let $|s| \leq R$. Then the network
\begin{align*}
\hat{f}_{test}(x, \mathbf{a}, \mathbf{b}, s) &= \sigma\bigg(\hat{f}_{id}(s)-R^2 \cdot \sum_{i=1}^d \left(\sigma\left(a^{(i)} + \frac{1}{R} - x^{(i)}\right)\right.\\
& \left. \hspace*{4cm} + \sigma\left(x^{(i)} - b^{(i)} + \frac{1}{R}\right)\right)\bigg)\\
& \quad - \sigma\bigg(-\hat{f}_{id}(s)-R^2 \cdot \sum_{i=1}^d \left(\sigma\left(a^{(i)} + \frac{1}{R} - x^{(i)}\right)\right.\\
& \left. \hspace*{4cm} + \sigma\left(x^{(i)} - b^{(i)} + \frac{1}{R}\right)\right)\bigg)
\end{align*}
of the class $\mathcal{F}(2, 2 \cdot (2d+2))$
satisfies
the weight constraint 
$$\norm{\bv_{\hat{f}_{test} }}_\infty\leq R^2,$$
as well as 
$$
\norm{\left(\bv_{\hat{f}_{test}}\right)^{(0)}_{i,j>0}}_\infty=1,
\quad
\left(\bv_{\hat{f}_{test}}\right)^{(2)}_{1,0}=0 
\quad \text{ and } \quad
\norm{\left(\bv_{\hat{f}_{test}}\right)^{(2)}_{1,i>0}}_\infty=1.
$$
For $x \in K_{1/R}$ $\hat{f}_{test}(x, \mathbf{a}, \mathbf{b}, s)$ satisfies 
\begin{align*}
  \hat{f}_{test}(x, \mathbf{a}, \mathbf{b}, s)
  =
  s \cdot \mathds{1}_{[\bold{a}, \bold{b})}(x)
\end{align*}
 and 
\begin{align*}
  \left|\hat{f}_{test}(x, \mathbf{a}, \mathbf{b}, s) -
  s \cdot \mathds{1}_{[\bold{a}, \bold{b})}(x)\right| \leq |s|
\end{align*}
for $x \in \Rd$.
\end{lemma}
\noindent
{\bf Proof.}
The weight constraints can easily be
seen in the definition of $\hat{f}_{ind, [\bold{a}, \bold{b})}$ and $\hat{f}_{test}$, the proof of the approximation bounds can be found in the proof of Lemma 6 in Kohler and Langer (2021).
\hfill $\Box$

\begin{lemma}
\label{supple12}
Let $\sigma:\R \to \R$ be the ReLU activation function $\sigma(x) = \max\{x,0\}$. Let $R >0$, $B \in \N$ and
\begin{align*}
\hat{f}_{ind, [j, \infty)}(z) = R \cdot \sigma(z-j) - R \cdot \sigma\left(z-j-\frac{1}{R}\right) \in \mathcal{F}(1, 2)
\end{align*}
for $j \in \{1, \dots, B\}$. Then the neural network
\begin{align*}
\hat{f}_{trunc}(z)= \sum_{j=1}^B \hat{f}_{ind, [j, \infty)}(z) \in \mathcal{F}(1, 2B)
\end{align*}
satisfies
$$\norm{\bv_{\hat{f}_{trunc}}}_\infty\leq \max\left\{R,B+\frac{1}{R}\right\},$$ more specifically the network has no offset in its last layer, i.e. $\left(\bv_{\hat{f}_{trunc}}\right)_{1,0}^{(1)}=0$, and satisfies $ $$\norm{\left(\bv_{\hat{f}_{trunc}}\right)_{i,j>0}^{(0)}}_\infty=1$ and $\norm{\left(\bv_{\hat{f}_{trunc}}\right)_{i,j>0}^{(1)}}_\infty=R$.\\
Further, $\hat{f}_{trunc}$ satisfies
\begin{align*}
\hat{f}_{trunc}(z) = \lfloor z \rfloor
\end{align*}
for $z \in [0, B+1)$ and $ \min\{ |z-j| \, : \, j \in \N\} \geq 1/R$. 
\end{lemma}

\noindent
{\bf Proof.}
     The weight constraints can easily be derived from the definition of $\hat{f}_{trunc}$, the proof of the approximation bounds can be found in the proof of Lemma 13 in Kohler and Langer (2021).
\hfill $\Box$

In the proof of Lemma \ref{supple13} every function of $\bm{\phi}_{1, 2M^d+1}$ is computed by a neural network. In particular, the indicator functions in $\bm{\phi}_{2, j}$, $\phi^{(\bll)}_{3,j}$ and $\phi_{4,j}^{(\bll)}$ $(j \in \{1, \dots, M^d\}, \bll \in \N_0^d, \|\bll\|_1 \leq q)$ are computed by Lemma \ref{le:f_ind_f_test} a), while we apply the identity network to shift the computed values from the previous step. The functions $\phi_{3, M^d+j}^{(\bll)}$ and $\phi_{4, M^d+j}^{(\bll)}$ are then computed according to their definition above Lemma \ref{supple11}, while we again use the identity network to shift values in the next hidden layers. For the functions $\bm{\phi}_{5, M^d+j}$ and $\phi_{6, M^d+j}^{(\bll)}$ we use the network of Lemma \ref{le:f_ind_f_test} b) to successively compute $(C_{\mathcal{P}_2}(x))_{left}$ and the derivatives on the cube $C_{\mathcal{P}_2}(x)$. The final Taylor polynomial in $\phi_{1, 2M^d+1}$ is then approximated with the help of Lemma \ref{le:polynomial-approx}.

\noindent
{\bf Proof of Lemma \ref{supple13}.}
In the \textit{first step of the proof} we describe how $\phi_{1, 2M^d+1}$ of Lemma \ref{supple11} can be approximated by neural networks. In the construction we will use the network
$
\hat{f}_{ind, [\bold{a}, \bold{b})} \in \mathcal{F}(2, 2d)$
and the network
$
\hat{f}_{test}  \in \mathcal{F}(2, 2 \cdot (2d+2))
$
of Lemma \ref{le:f_ind_f_test}.
Here we set $R=B_M=M^{2p+2}$ in Lemma \ref{le:f_ind_f_test}, such that the weights of the network satisfy the constraints 
$$\norm{\bv_{\hat{f}_{ind} }}_\infty\leq M^{2p+2}\quad \text{and}\quad \norm{\bv_{\hat{f}_{test} }}_\infty\leq M^{4p+4}.$$

For some vector $\bold{v} \in \R^d$ we set
\begin{align*}
&\bold{v} \cdot \hat{f}_{ind, [\bold{a}, \bold{b})}(x) = \left(v^{(1)} \cdot \hat{f}_{ind, [\bold{a}, \bold{b})}(x), \dots, v^{(d)} \cdot \hat{f}_{ind, [\bold{a}, \bold{b})}(x)\right).
\end{align*}
Furthermore we use the networks 
\begin{align*}
\hat{f}_{trunc,i} \in \mathcal{F}(1, 2 \cdot (4+2 \lceil e^d\rceil)), \quad (i \in \{1, \dots, M^d-1\})
\end{align*}
of Lemma \ref{supple12}.
  Here we choose
  \[
  R=R_{M,i} = {(4 + 2 \lceil e^d \rceil)^{M^d-i-1}}
  \quad \mbox{and} \quad
  B=4+ 2 \lceil e^d\rceil
  \]
 in Lemma \ref{supple12}, which implies
 \begin{equation}
 \label{ftruncNorm_phiProof}
     \norm{\bv_{\hat{f}_{trunc}}}_\infty\leq \max\left\{R,B+\frac{1}{R}\right\}\leq 4+ 2 \lceil e^d\rceil+(4 + 2 \lceil e^d \rceil)^{M^d}.
 \end{equation}
\\
To compute the final Taylor polynomial we use the network
\begin{align*}
\hat{f}_{p} \in \mathcal{F}\left(B_{M,p} \cdot \lceil \log_2(\max\{q+1, 2\}) \rceil, 18 \cdot (q+1) \cdot \binom{d+q}{d}\right)
\end{align*} 
from Lemma \ref{le:polynomial-approx} satisfying 
$$
\begin{aligned}
    \norm{\bv_{\hat{f}_p}}_\infty
\leq 4\cdot \bar{r}(p)
\cdot 4^{2(q+1)}\cdot 
\bigg(
&2 \cdot \max\left\{\|f\|_{C^q([-a,a]^d)}, a\right\} 
\cdot e^{(M^d-1)} \\
&+  (4 + 2 \lceil e^d \rceil) 
\cdot (M^d-1) 
\cdot e^{(M^d-2)}\bigg)^{2(q+1)},
\end{aligned}
$$
and
\begin{align}
\label{2bfpeq}
&\left|\hat{f}_{p}\left(\bold{z}, y_1, \dots, y_{\binom{d+q}{q}}\right) - p\left(\bold{z}, y_1, \dots, y_{\binom{d+q}{q}}\right)\right| \notag\\
& \leq c_{41}
\cdot (6 + 2 \lceil e^d \rceil)^{4(q+1)}
\cdot \bar{r}(p) \cdot \left(\max\left\{\|f\|_{C^{q}([-a,a]^d)}, a \right\}\right)^{4(q+1)}  \notag\\
& \quad \cdot M^{d \cdot 4 \cdot (q+1)} \cdot e^{4(q+1) \cdot (M^d-1)}  \cdot 4^{-B_{M,p}}
\end{align}
for all $z^{(1)}, \dots, z^{(d)}, y_1, \dots, y_{\binom{d+q}{d}}$ contained in
\begin{align*}
&\left[- 2 \cdot \max\left\{\|f\|_{C^q([-a,a]^d)}, a\right\} \cdot e^{(M^d-1)} + (4 + 2 \lceil e^d \rceil) \cdot (M^d-1) \cdot e^{(M^d-2)}, \right.\\
&\left. \quad 2 \cdot \max\left\{\|f\|_{C^q([-a,a]^d)}, a\right\} \cdot e^{(M^d-1)} + (4 + 2 \lceil e^d \rceil) \cdot (M^d-1) \cdot e^{(M^d-2)} \right], 
\end{align*}
where 
\begin{align*}
R=B_{M,p} = \left\lceil \log_4\left(M^{2p+4 \cdot d \cdot (q+1)} \cdot e^{4 \cdot (q+1) \cdot (M^d-1)}\right)\right\rceil.
\end{align*}
A polynomial of degree zero is treated as a polynomial of degree $1$, where we choose $r_i = 0$ for all coefficients greater than zero. Thus we substitute $\log_2(q+1)$ by
$\log_2(\max\{q+1, 2\})$ in the definition of $L$ in Lemma \ref{le:polynomial-approx}.
To compute $\bm{\phi}_{1, j}, \bm{\phi}_{2, j}, \phi_{3, j}^{(\bll)}$ and $\phi_{4, j}^{(\bll)}$ for $j \in \{0, \dots, M^d\}$ and each $\bll \in \N_0^d$ with $\|\bll\|_1 \leq q$ we use the networks
\begin{align*}
&\bm{\hat{\phi}}_{1,0} = \left(\hat{\phi}_{1,0}^{(1)}, \dots, \hat{\phi}_{1,0}^{(d)}\right) = x\\
&\bm{\hat{\phi}}_{2,0} = \left(\hat{\phi}_{2,0}^{(1)}, \dots, \hat{\phi}_{2,0}^{(d)}\right) = \mathbf{0}, \\
& \hat{\phi}_{3, 0}^{(\bll)} =0 \ \mbox{and} \ \hat{\phi}_{4, 0}^{(\bll)} =0.
\end{align*}
for $\bll \in \N_0^d$ with $\|\bll\|_1 \leq q$. For $j \in \{1, \dots, M^d\}$ we set
\begin{align*}
&\bm{\hat{\phi}}_{1,j} = \hat{f}_{id}^2\left(\bm{\hat{\phi}}_{1, j-1}\right),\\
&\bm{\hat{\phi}}_{2, j} = (C_{j,1})_{left} \cdot \hat{f}_{ind, C_{j,1}}(\bm{\hat{\phi}}_{1,j-1}) + \hat{f}_{id}^2(\bm{\hat{\phi}}_{2, j-1}),\\
&\hat{\phi}_{3, j}^{(\bll)} = (\partial^{\bll} f)((C_{j,1})_{left}) \cdot \hat{f}_{ind, C_{j,1}}(\bm{\hat{\phi}}_{1, j-1})+ \hat{f}_{id}^2(\hat{\phi}_{3, j-1}^{(\bll)}),\\
&\hat{\phi}_{4, j}^{(\bll)} = b_j^{(\bll)} \cdot \hat{f}_{ind, C_{j,1}}(\bm{\hat{\phi}}_{1, j-1})+ \hat{f}_{id}^2(\hat{\phi}_{4, j-1}^{(\bll)})\\
\end{align*}
for $\bll \in \N_0^d$ with $\|\bll\|_1 \leq q$.\\

It is easy to see that this parallelized network needs $2M^d$ hidden layers and $2d+d \cdot (2d+2)+2 \cdot \binom{d+q}{d} \cdot (2d+2)$ neurons per layer, where we have used
the fact that we have $\binom{d+q}{d}$ different vectors $\bll \in \N_0^d$ satisfying $\|\bll\|_1 \leq q$. \\
\noindent
To compute $\bm{\phi}_{1, M^d+j}, \bm{\phi}_{5, M^d+j}$ and $\phi_{6, M^d+j}^{(\bll)}$ for $j \in \{1, \dots, M^d\}$ and \linebreak $\bm{\phi}_{2, M^d+j}, \phi_{3,  M^d+j}^{(\bll)}$ and $\phi_{4, M^d+j}^{(\bll)}$ for $j \in \{1, \dots, M^d-1\}$ we use the networks
\begin{align*}
&\bm{\hat{\phi}}_{1, M^d+j} = \hat{f}_{id}^2\left(\bm{\hat{\phi}}_{1, M^d+j-1}\right), \quad j \in \{1, \dots, M^d\}\\
&\bm{\hat{\phi}}_{2, M^d+j} = \hat{f}_{id}^2\left(\bm{\hat{\phi}}_{2, M^d+j-1} + \bold{\tilde{v}}_{j+1}\right)\\
& \hat{\phi}_{3,M^d+j}^{(\bll)} = \hat{f}_{id}\Bigg(\hat{f}_{id}\bigg(\sum_{\substack{\mathbf{s}\in \N_0^d\\ \|\bold{s}\|_1 \leq q-\|\bold{l}\|_1}} \frac{\hat{\phi}_{3, M^d+j-1}^{(\bll+\mathbf{s})}}{\bold{s}!} \cdot \left(\bold{\tilde{v}}_{j+1}\right)^{\bold{s}}\bigg)\\
& \hspace{3cm}+ \left(\hat{f}_{trunc,j}\left((4+2 \cdot \lceil e^d \rceil) \cdot \hat{\phi}_{4, M^d+j-1}^{(\bll)}\right)\right.\\
&\hspace{3.5cm} \left. - \lceil e^d \rceil -2\right) \cdot c_{36} \cdot \left(\frac{2a}{M^2}\right)^{p-\|\bll\|_1}\Bigg), \\
&\hat{\phi}_{4, M^d+j}^{(\bll)} = \hat{f}_{id}\left(\hat{f}_{id}\left((4+2 \cdot \lceil e^d \rceil) \cdot \hat{\phi}_{4, M^d+j-1}^{(\bll)}\right)\right.\\
& \hspace{3cm} \left. - \hat{f}_{trunc,j}\left((4+2 \cdot \lceil e^d \rceil) \cdot \hat{\phi}_{4, M^d+j-1}^{(\bll)}\right)\right)
\end{align*}
for $j \in \{1, \dots, M^d-1\}$.\\
Further we set 
\begin{align}
\label{2bneur5}
\hat{\phi}_{5, M^d+j}^{(k)} &= \hat{f}_{test}\left(\bm{\hat{\phi}}_{1, M^d+j-1}, \bm{\hat{\phi}}_{2, M^d+j-1}, \right. \notag\\
& \hspace*{1.5cm} \left.\bm{\hat{\phi}}_{2, M^d+j-1}+\frac{2a}{M^2} \cdot \mathbf{1}, \hat{\phi}_{2, M^d+j-1}^{(k)}\right) \notag\\
& \quad + \hat{f}_{id}^2\left(\hat{\phi}_{5, M^d+j-1}^{(k)}\right)
\end{align}
and
\begin{align}
\label{2bneur6}
\hat{\phi}_{6, M^d+j}^{(\bll)} &= \hat{f}_{test}\left(\bm{\hat{\phi}}_{1, M^d+j-1}, \bm{\hat{\phi}}_{2, M^d+j-1}, \right. \notag\\
& \hspace*{1.5cm} \left. \bm{\hat{\phi}}_{2, M^d+j-1}+\frac{2a}{M^2} \cdot \mathbf{1}, \hat{\phi}_{3, M^d+j-1}^{(\bll)}\right) \notag\\
& \quad + \hat{f}_{id}^2\left(\hat{\phi}_{6, M^d+j-1}^{(\bll)}\right),
\end{align}
where $\bm{\hat{\phi}}_{5, M^d} = \left(\hat{\phi}_{5, M^d}^{(1)}, \dots, \hat{\phi}_{5, M^d}^{(d)}\right) = \mathbf{0}$ and $\hat{\phi}_{6, M^d}^{(\bll)} = 0$ for each $\bll \in \N_0^d$ with $\|\bll_1\| \leq q$.

Again it is easy to see, that this parallelized and composed network needs $4M^d$ hidden layers and has width $r$ with 
\begin{equation*}
r= 10d+4d^2+2 \cdot \binom{d+q}{d} \cdot \left(2 \cdot (4+2\lceil e^d\rceil)+5+2d\right).
\end{equation*}
Choose $\bll_1, \dots, \bll_{\binom{d+q}{d}}$ such that
\begin{align*}
\left\{\bll_1, \dots, \bll_{\binom{d+q}{d}}\right\} = \left\{\bold{s} \in \N_0^d: \|\bold{s}\|_1 \leq q \right\}
\end{align*}
holds. 
The value of $\phi_{1, 2M^d+1}$ can then be computed by 
\begin{align}
\label{2bfp}
\hat{\phi}_{1, 2M^d+1} = \hat{f}_p\left(\bold{z}, y_1, \dots, y_{\binom{d+q}{d}}\right),
\end{align}
where 
\begin{align*}
\bold{z}= \bm{\hat{\phi}}_{1, 2M^d} - \bm{\hat{\phi}}_{5, 2M^d}
\end{align*}
and 
\begin{align*}
y_v = \hat{\phi}_{6, 2M^d}^{(\bll_v)} 
\end{align*}
for $v \in \left\{1, \dots, \binom{d+q}{d}\right\}$. The coefficients $r_1, \dots, r_{\binom{d+q}{d}}$ in Lemma \ref{le:polynomial-approx} are chosen as 
\begin{align*}
r_i = \frac{1}{\bll_i!}, \quad i \in \left\{1, \dots, \binom{d+q}{d}\right\},
\end{align*}
i.e. $\bar{r}(p)\leq 1$.
The final network $\hat{\phi}_{1, 2M^d+1}$ is then contained in the class
\begin{align*}
\mathcal{F}(4M^d+B_{M,p} \cdot \lceil \log_2(\max\{q+1, 2\})\rceil, r)
\end{align*}
with 
\begin{align*}
r=&\max\left\{10d+4d^2+2 \cdot \binom{d+q}{d} \cdot \left(2 \cdot (4+2\lceil e^d\rceil)+5+2d\right), \right.\\
& \left. \hspace*{1cm} 18 \cdot (q+1) \cdot \binom{d+q}{d}\right\}
\end{align*}
and we set
\begin{align*}
\hat{f}_{deep, \P_2}(x) = \hat{\phi}_{1, 2M^d+1}.
\end{align*}
Note that the repeated composition of $\hat{f}_{id}$ and $\hat{f}_{ind, C_{j,1}}$ does not affect any weight constraints, since both networks have no offset in their respective output layers and the weights used in the input and output (i.e. last) layers are bounded by 1 respectively, such that Lemma \ref{lemma:composition-weight-bound} c) is applicable. This also holds for the repeated composition of $\hat{f}_{id}$ and $\hat{f}_{test}$.
Since $\norm{\bv_{\hat{f}^t_{id}}}_\infty=1$, we have $\norm{\bv_{\bm{\hat{\phi}}_{1,j}}}_\infty = 1$ for $j\in\{1,\ldots, 2M^d\}$. \\
From Lemma \ref{le:f_ind_f_test}, $((C_{j,1})_{left})<a$, the conditions on $M$ and since by definition $\bold{\tilde{v}}_{j+1}$ has entries in $\{0, \frac{2a}{M^2}\}$, we conclude that $\bm{\hat{\phi}}_{2,j}$ satisfies the weight constraint $\norm{\bv_{\bm{\hat{\phi}}_{2,j}}}_\infty \leq \max\left\{a+\frac{1}{B_M}, B_M\right\}\leq M^{2p+2}$ for $j\in\{1,\ldots, 2M^d\}$.\\
Analogously, we can derive the weight constraints for $\hat{\phi}_{3,j}^{(\mathbf{l})}$ and $\hat{\phi}_{4, j}^{(\bll)}$:
$$
\norm{\bv_{\hat{\phi}_{3,j}^{(\mathbf{l})}}}_\infty\leq  
\max
\left\{
\|f\|_{C^q([-a,a]^d)},
a+\frac{1}{B_M}, 
B_M\right\}\leq M^{2p+2}+c_{42}(f),
$$
respectively for $\hat{\phi}_{4, j}^{(\bll)}$, we have
$$\norm{\bv_{\bm{\hat{\phi}}_{4,j}^{(\mathbf{l})}}}_\infty \leq \max\left\{a+\frac{1}{B_M}, B_M\right\}\leq M^{2p+2} \quad \text{and}$$
$\norm{\left(\bv_{\hat{\phi}_{4, j}^{(\bll)}}\right)^{(0)}_{1,i>0}}_\infty=1$,
$\left(\bv_{\hat{\phi}_{4, j}^{(\bll)}}\right)^{(2j)}_{1,0}=0$ 
and 
$\norm{\left(\bv_{\hat{\phi}_{4, j}^{(\bll)}}\right)^{(2j)}_{1,i>0}}_\infty\leq 1$,
where we have used Lemma \ref{le:f_ind_f_test} and that by definition
\begin{equation*}
   b_i^{(\bll)} \in [0,1].
\end{equation*}
Note that by construction $\left(\bv_{\bm{\hat{\phi}}_{2, 2M^d}}\right)^{(4M^d)}_{1,0}=0$ and $\norm{\left(\bv_{\bm{\hat{\phi}}_{2, 2M^d}}\right)^{(4M^d)}}_\infty\leq 1$.
Since we have 
$\norm{\left(\bv_{\hat{\phi}_{4, j}^{(\bll)}}\right)^{(2M^d)}_{1,i>0}}_\infty\leq 1$, 
$\norm{\left(\bv_{\hat{f}_{trunc}}\right)_{1,i>0}^{(0)}}_\infty=1$ and 
$\left(\bv_{\hat{\phi}_{4, j}^{(\bll)}}\right)^{(2j)}_{1,0}=0$ 
by Lemma \ref{lemma:composition-weight-bound} c) the composition of $\hat{f}_{trunc}$ and $\hat{\phi}_{4, j}^{(\bll)}$  does not affect the weight constraints. Thus  $\hat{\phi}_{3,M^d+j}^{(\bll)}$ and  $\hat{\phi}_{4,M^d+j}^{(\bll)}$ satisfy the constraint given in \eqref{ftruncNorm_phiProof}, i.e. 
$$
\norm{\bv_{\hat{\phi}_{3,M^d+j}^{(\bll)}}}_\infty=
\norm{\bv_{\hat{\phi}_{4,M^d+j}^{(\bll)}}}_\infty 
\leq 4+ 2 \lceil e^d \rceil+(4 + 2 \lceil e^d \rceil)^{M^d}+
(c_{36}+1) \cdot M^{2p+2}.
$$
By Lemma \ref{le:f_ind_f_test} a) and the previously stated weight constraints for $\bm{\hat{\phi}}_{1, M^d+j}$, $\bm{\hat{\phi}}_{2, M^d+j}$ and $\bm{\hat{\phi}}_{3, M^d+j}$, we have
$$
\begin{aligned}
   &\norm{\bv_{\bm{\hat{\phi}}_{5,M^d+j}}}_\infty
\leq 
M^{4p+4},
\quad
\norm{\left(\bv_{\hat{\phi}_{5,M^d+j}}\right)^{(2M^d+2j)}}_\infty\leq 1,
\quad 
\left(\bv_{\hat{\phi}_{5,M^d+j}}\right)^{(2M^d+2j)}_{1,0}=0 
 \quad \text{and}\\
&\norm{\bv_{\hat{\phi}_{6,M^d+j}}}_\infty
\leq 
\max
\left\{
M^{4p+4},
\norm{\bv_{\hat{\phi}_{3,M^d+j}^{(\bll)}}}_\infty
\right\}
\leq
4+ 2 \lceil e^d\rceil
+(4 + 2 \lceil e^d \rceil)^{M^d}+ (c_{36}+1) \cdot M^{4p+4},\\
&
\norm{\left(\bv_{\hat{\phi}_{6,M^d+j}}\right)^{(2M^d+2j)}}_\infty\leq 1,
\quad 
\left(\bv_{\hat{\phi}_{6,M^d+j}}\right)^{(2M^d+2j)}_{1,0}=0 
. 
\end{aligned}
$$
To derive the weight constraints for the final network $\hat{f}_{deep, \P_2}(x)$, note that 
using the previously stated constraints for the weights of $\bm{\hat{\phi}}_{1, 2M^d}$, $\bm{\hat{\phi}}_{5,2M^d}$ and $\bm{\hat{\phi}}_{6,2M^d}$ as well as Lemma \ref{le:polynomial-approx}, we know that
the composition of $\hat{f}_p$ and $\bm{\hat{\phi}}_{1, 2M^d} - \bm{\hat{\phi}}_{5, 2M^d}$ fulfills the conditions of Lemma \ref{lemma:composition-weight-bound} c).
This results in 
$$
\begin{aligned}
    \norm{\bv_{\hat{f}_{deep, \P_2}}}_\infty
    &\leq \max\{4\cdot 4^{2(q+1)}\cdot \left(2 \cdot \max\left\{\|f\|_{C^q([-a,a]^d)}, a\right\} \cdot e^{(M^d-1)} + (4 + 2 \lceil e^d \rceil)\right. \\
    &\hspace{40pt}\left.\left.\cdot (M^d-1) \cdot e^{(M^d-2)}\right)^{2(q+1)},
    (c_{36}+1) \cdot M^{4p+4}\right\}\\
    & \leq  \left(c_{43}(f)\cdot e^{M^d-1}\cdot\left(6+2\lceil e^d\rceil\right)\cdot M^d \right)^{2(q+1)}+ (c_{36}+1) \cdot M^{4p+4}\\
    &\leq e^{c_{44}(f) \cdot(p+1) \cdot M^d}
\end{aligned}
$$

The rest of the proof follows as in Kohler and Langer (2021).
\hfill $\Box$

\subsubsection{Key step 3: Approximating $w_{\P_2}(x) \cdot f(x)$ by deep neural networks} 
In order to approximate $f(x)$ in supremum norm, a neural network that approximates $w_{\P_2}(x) \cdot f(x)$, where $w_{\P_2}(x)$ is defined as in \eqref{w_vb}, is required. The construction of such a network is given in the following three results.
\begin{lemma}
\label{supple15}
Let $\sigma: \R \to \R$ be the ReLU activation function $\sigma(x) = \max\{x,0\}$. Let $1 \leq a < \infty$ and $M \in \N_0$ sufficiently large (independent of the size of $a$, but 
    \begin{eqnarray*}
      M^{2p} &\geq&
      2^{4(q+1)} \cdot
      \max\{ c_{45} (6+2 \lceil e^d \rceil)^{4(q+1)}, c_{36} \cdot e^d \}
      \\
      &&
      \hspace*{5cm}
      \cdot \left(\max\left\{a, \|f\|_{C^q([-a,a]^d)}\right\}\right)^{4(q+1)}
    \end{eqnarray*}
    must hold).
Let $p=q+s$ for some $q \in \N_0$, $s \in (0,1]$ and let $C>0$.
    Let $f: \Rd \to \R$ be a $(p,C)$-smooth function and let $w_{\P_2}$ be defined as in
\eqref{w_vb}. Then there exists a network
\begin{align*}
\hat{f} \in \mathcal{F}\left(L, r\right)
\end{align*}
with
\begin{align*}
L=&5M^d+\left\lceil \log_4\left(M^{2p+4 \cdot d \cdot (q+1)} \cdot e^{4 \cdot (q+1) \cdot (M^d-1)}\right)\right\rceil  \cdot \lceil \log_2(\max\{q,d\}+1)\rceil \\
&+ \lceil \log_4(M^{2p})\rceil
\end{align*}
and
\begin{align*}
r=&\max\left\{10d+4d^2+2 \cdot \binom{d+q}{d} \cdot \left(2 \cdot (4+2\lceil e^d\rceil)+5+2d\right), \right.\\
&\left. \hspace*{1cm} 18 \cdot (q+1) \cdot \binom{d+q}{d}\right\} + 6d^2+20d+2.
\end{align*}
which satisfies 
$$
\norm{\bv_{\hat{f}}}_\infty \leq e^{c_{46}(f) \cdot(p+1) \cdot M^d}
$$
such that
\begin{align*}
&\left|\hat{f}(x) - w_{\P_2}(x) \cdot f(x)\right| \leq c_{47} \cdot \left(\max\left\{2a, \|f\|_{C^q([-a,a]^d)}\right\}\right)^{4(q+1)} \cdot \frac{1}{M^{2p}}
\end{align*}
for $x \in [-a,a)^d$. 
\end{lemma}
Further auxiliary lemmata are required to show this result. First, it is shown 
that each weight $w_{\P_2}(x)$ can also be approximated by a very deep neural network. 
\begin{lemma}
\label{supple14}
Let $\sigma: \R \to \R$ be the ReLU activation function $\sigma(x) = \max\{x,0\}$. Let $1 \leq a < \infty$ and $M \geq 4^{4d+1} \cdot d$. Let $\mathcal{P}_{2}$
be the partition defined in (\ref{partition}) and let $w_{\P_2}(x)$ be
defined by \eqref{w_vb}. Then there exists a neural network
\begin{align*}
\hat{f}_{w_{\P_2},deep} \in \mathcal{F}\left(L, r\right),
\end{align*}
with
\begin{align*}
L=4M^d+1+\lceil \log_4(M^{2p})\rceil \cdot \lceil \log_2(d)\rceil\quad \text{and} \quad r= \max\left\{18d, 4d^2+10d\right\}
\end{align*}
which satisfies 
$$
\norm{\bv_{\hat{f}_{w_{\P_2,deep}}}}_\infty
\leq  
M^{4p+4},
\quad
\left(\bv_{\hat{f}_{w_{\P_2},deep}}\right)^{(L)}_{1,0}= 0
\quad
\text{and}
\quad
\norm{\left(\bv_{\hat{f}_{w_{\P_2},deep}}\right)^{(L)}_{1,j>0}}_\infty
\leq  
4^{3d+1}
$$
such that
\begin{align*}
\left|\hat{f}_{w_{\P_2},deep}(x) - w_{\P_2}(x)\right| \leq 4^{4d+1} \cdot d \cdot \frac{1}{M^{2p}}
\end{align*}
for $x \in \bigcup_{i \in \{1, \dots, M^{2d}\}} (C_{i,2})_{1/M^{2p+2}}^0$ and 
\begin{align*}
|\hat{f}_{w_{\P_2, deep}}(x)| \leq 1
\end{align*}
for $x \in [-a,a)^d$.
\end{lemma}

\noindent
{\bf Proof.}
The first $4M^d$ hidden layers of $\hat{f}_{w_{\P_2}}$ compute the value of 
\begin{align*}
(C_{\mathcal{P}_{2}}(x))_{left}
\end{align*}
 using $\bm{\hat{\phi}}_{5, 2M^d}$ of Lemma \ref{supple13} (with $d \cdot (2 \cdot (2d+2)+2)+2d$ neurons per layer) and shift the value of $x$ in the next hidden layer using the network $\hat{f}_{id}^{4M^d}$.
As stated in the proof of Lemma \ref{supple13}, $\bm{\hat{\phi}}_{5, 2M^d}$ has the following weight constraints 
$\norm{\bv_{\bm{\hat{\phi}}_{5,M^d+j}}}_\infty
\leq 
M^{4p+4},
\quad
\norm{\left(\bv_{\hat{\phi}_{5,M^d+j}}\right)^{(2M^d+2j)}}_\infty\leq 1$
{and}
$
\left(\bv_{\hat{\phi}_{5,M^d+j}}\right)^{(2M^d+2j)}_{1,0}=0
$.
The next hidden layer then computes the functions
\begin{eqnarray*}
&&\left(1-\frac{M^2}{a} \cdot \left|(C_{\mathcal{P}_{2}}(x))_{left}^{(j)} + \frac{a}{M^2} -x^{(j)}\right| \right)_+\\
  &&=
  \left(
  \frac{M^2}{a} \cdot
  \left(
x^{(j)} - (C_{\mathcal{P}_{2}}(x))_{left}^{(j)} 
  \right)
  \right)_+
  \\
  &&
  \quad
  -
  2 \cdot \left(
  \frac{M^2}{a} \cdot
  \left(
  x^{(j)} - (C_{\mathcal{P}_{2}}(x))_{left}^{(j)}
  - \frac{a}{M^2}
  \right)
  \right)_+
 \\
  &&
  \quad
  +
  \left(
  \frac{M^2}{a} \cdot
  \left(
  x^{(j)} - (C_{\mathcal{P}_{2}}(x))_{left}^{(j)}
  - \frac{2 \cdot a}{M^2}
  \right)
  \right)_+, \quad j \in \{1, \dots, d\},
\end{eqnarray*}
 using the networks
\begin{align*}
  \hat{f}_{w_{{\P_2},j}}(x) &= \sigma\left(
  \frac{M^2}{a} \cdot
  \left(
\hat{f}_{id}^{4M^d}(x^{(j)}) - {\hat{\phi}}_{5, 2M^d}^{(j)} 
  \right)
  \right)\\
  & \quad -2 \cdot \sigma\left(
\frac{M^2}{a} \cdot
  \left(
  \hat{f}_{id}^{4M^d}(x^{(j)}) - {\hat{\phi}}_{5, 2M^d}^{(j)}
  - \frac{a}{M^2}
  \right) 
  \right)\\
  & \quad + \sigma\left(
  \frac{M^2}{a} \cdot
  \left(
  \hat{f}_{id}^{4M^d}(x^{(j)})- {\hat{\phi}}_{5, 2M^d}^{(j)}
  - \frac{2 \cdot a}{M^2}
  \right)
  \right)
\end{align*}
with $3d$ neurons in the last layer. 
Note that  $|\hat{f}_{w_{\P_2},j}(x)| \leq 1$ for $j \in \{1, \dots, d\}$. The product of $w_{\P_2,j}(x)$ $(j \in \{1, \dots, d\})$
can then be computed by the network $\hat{f}_{mult,d}(x)$ of Lemma \ref{lemma:fmultd} for values $x\in[-1,1]^d$, where we choose $x^{(j)} = \hat{f}_{w_{{\P_2},j}}(x)$ and $R= \lceil \log_4(M^{2p})\rceil$. Finally we set
\begin{align*}
\hat{f}_{w_{\P_2},deep}(x) = \hat{f}_{mult, d}\left(\hat{f}_{w_{{\P_2},1}}(x), \dots, \hat{f}_{w_{{\P_2},d}}(x)\right).
\end{align*}
Since by Lemma \ref{lemma:fmultd}, in this case $\hat{f}_{mult,d}$ satisfies
$$\norm{\bv_{\hat{f}_{mult,d}}}_\infty\leq 4\cdot 4^{2d},\qquad \left(\bv_{\hat{f}_{mult,d}}\right)_{1,0}^{(R \cdot \lceil \log_2(d) \rceil)} =0\qquad \text{and}$$
$$
\norm{\left(\bv_{\hat{f}_{mult,d}}\right)^{(0)}}_\infty\leq 1
$$
and by construction 
$$
\left(\bv_{\hat{f}_{w_{\P_2},j}}\right)^{(4M^d+1)}_{1,0}= 0, \quad \norm{\left(\bv_{\hat{f}_{w_{\P_2},j}}\right)^{(4M^d+1)}_{1,i>0}}_\infty\leq 2,$$
$$
\norm{\bv_{\hat{f}_{w_{\P_2},j}}}_\infty
= 
M^{4p+4}
$$
an application Lemma \ref{lemma:composition-weight-bound} b) gives us 
$$\norm{\bv_{\hat{f}_{w_{\P_2},deep}}}_\infty
\leq 
\max
\left
\{\norm{\bv_{\hat{f}_{w_{\P_2},j}}}_\infty, 4^{2d+1}, 2
\right\}
= 
M^{4p+4}
$$
The rest of the proof follows as in the proof of Lemma 16 in Kohler and Langer (2021).

\hfill $\Box$

\begin{lemma}
\label{supple17}
Let $\sigma: \R \to \R$ be the ReLU activation function $\sigma(x) = \max\{x,0\}$. Let $1 \leq a < \infty$. Let $C_{i,2}$ $(i \in \{1, \dots, M^{2d}\})$
be the cubes of partition $\mathcal{P}_{2}$ as described in \eqref{partition} and let $M \in \N$. Then there exists a neural network 
\begin{align*}
\hat{f}_{check, deep, \mathcal{P}_{2}}(x) \in \mathcal{F}\left(5M^d, 2d^2+6d+2\right)
\end{align*}
satisfying
$$
\norm{\bv_{\hat{f}_{check, deep, \mathcal{P}_{2}}}}_\infty \leq M^{4p+4}, 
\quad
\norm{\bv_{\hat{f}_{check, deep, \mathcal{P}_{2}}}^{(5M^d)}}_\infty\leq 1
$$
\begin{align*}
  \hat{f}_{check, deep,\mathcal{P}_{2}}(x) = \mathds{1}_{
    \bigcup_{i \in \{1, \dots, M^{2d}\}}
    C_{i,2} \setminus (C_{i,2})_{1/M^{2p+2}}^0
}(x)
\end{align*}
for $x \notin \bigcup_{i \in \{1, \dots, M^{2d}\}} (C_{i,2})_{1/M^{2p+2}}^0 \textbackslash (C_{i,2})_{2/M^{2p+2}}^0$ and 
\begin{align*}
\hat{f}_{check, deep, \mathcal{P}_{2}}(x) \in [0,1]
\end{align*}
for $x \in [-a,a)^d$. 
\end{lemma}

\noindent
{\bf Proof.}
The value of $(C_{\P_1}(x))_{left}$ is computed by the network $\bm{\hat{\phi}}_{2, M^d}$ of Lemma \ref{supple13} with $2M^d$ hidden layers and $d \cdot (2d+2)$ neurons per layer and $x$ is shifted in consecutive layers by successively applying $\hat{f}_{id} \in \mathcal{F}(1, 2)$. 
 Furthermore we compute 
\begin{align*}
f_1(x)=1-\sum_{i \in \{1, \dots, M^d\}} \mathds{1}_{(C_{i,1})_{1/M^{2p+2}}^0}(x)
\end{align*}
by a network 
\begin{align*}
\hat{f}_{1, j}(x)= \hat{f}_{id}^2(\hat{f}_{1, j-1})-\hat{f}_{ind, (C_{j,1})^0_{1/M^{2p+2}}}(\hat{f}_{id}^{2(j-1)}(x)), \quad j \in \{1,\dots, M^d\},
\end{align*}
where $\hat{f}_{1, 0} = 1$. Here we use again the network $\hat{f}_{ind,[\bold{a},\bold{b})]}$ from Lemma \ref{le:f_ind_f_test} a with $R= M^{2p+2}$.
Next we define
\[
\bm{\bar{\phi}}_{2, M^d+j} = \hat{f}_{id}^3(\bm{\hat{\phi}}_{2, M^d+j-1} + \bold{\tilde{v}}_{j+1}) \in \mathcal{F}(2M^d+3j, 2d)
\]
for $j \in \{1, \dots, M^d\}$.
It is easy to see that $\bm{\bar{\phi}}_{2, M^d+j}$ satisfies the same weight constraints as $\bm{\bar{\phi}}_{2, M^d+j}$, i.e. $\norm{\bv_{\bm{\hat{\phi}}_{2,M^d+j}}}_\infty \leq M^{2p+2}$,which we derived in the proof of Lemma \ref{supple13} and by construction we have $\left(\bv_{\bm{\bar{\phi}}_{2, 2M^d}}\right)^{(4M^d)}_{1,0}=0$ and $\norm{\left(\bv_{\bm{\bar{\phi}}_{2, 2M^d}}\right)^{(4M^d)}}_\infty\leq 1$.
The value of 
\begin{align*}
\mathds{1}_{\bigcup_{i \in \{1, \dots, M^{2d}\}} C_{i,2} \textbackslash (C_{i,2})_{1/M^{2p+2}}^0}(x)
\end{align*}
is then successively computed by 
\begin{align*}
&\hat{f}_{1, M^d+j}(x)\\
&= 1-\sigma\left(1-\hat{f}_{test}\left(\hat{f}_{id}^{2M^d+3(j-1)}(x), \bm{\bar{\phi}}_{2, M^d+j-1} + \bold{\tilde{v}}_j + \frac{1}{M^{2p+2}} \cdot \mathbf{1}, \right. \right.\\
& \hspace*{1cm} \left. \left. \bm{\bar{\phi}}_{2, M^d+j-1} + \bold{\tilde{v}}_j + \frac{2a}{M^{2}} \cdot \mathbf{1}-\frac{1}{M^{2p+2}}\cdot \mathbf{1}, 1\right) - \hat{f}_{id}^2\left(\hat{f}_{1, M^d+j-1}\right)\right)
\end{align*}
for $j \in \{1, \dots, M^d\}$, where we use networks $\hat{f}_{test}$ from Lemma \ref{le:f_ind_f_test} b with $R= M^{2p+2}$ and which thus satisfies $\norm{\bv_{\hat{f}_{test}}}_\infty \leq M^{4p+4})$. 

Finally we set
\begin{align*}
\hat{f}_{check, deep, \P_2}(x) = \hat{f}_{1, 2M^d}(x).
\end{align*}
The weight constraints follow again from Lemma \ref{lemma:composition-weight-bound}, arguing in the same way as in the proof of Lemma \ref{supple13}.

\hfill $\Box$

In the proof of Lemma \ref{supple15} we use Lemma \ref{supple14} to approximate $w_{\P_2}(x)$ and Lemma \ref{supple13} to compute $f(x)$. As in Lemma 7 of Kohler and Langer (2021) we apply a network, that \textit{checks} whether $x$ is close to the boundaries of the cubes of the partition. 

\noindent
{\bf Proof of Lemma \ref{supple15}.}
This result follows by a straightforward modification of the proof of Lemma 7 of Kohler and Langer (2021). Here we use the network $\hat{f}_{deep, \P_2}$ of Lemma \ref{supple13} and $\hat{f}_{check, deep, \P_2}$ of Lemma \ref{supple15} to define 
\begin{eqnarray*}
  \hat{f}_{\P_2, true}(x) &=& \sigma\left(\hat{f}_{deep, \P_2}(x) - B_{true} \cdot \hat{f}_{check, deep, \P_2}(x)\right)
  \\
  &&
  -
  \sigma\left(-\hat{f}_{deep, \P_2}(x) - B_{true} \cdot \hat{f}_{check, deep, \P_2}(x)\right),
\end{eqnarray*}
with 
\begin{align*}
B_{true} &= 1 + \Bigg|\Bigg(\|f\|_{C^q([-a,a]^d)}\cdot e^{(M^d-1)}
\\
 & \quad 
 + (4+2 \cdot \lceil e^d \rceil) \cdot (M^d-1)\cdot e^{(M^d-2)}\Bigg)\cdot e^{2ad}\Bigg|
\end{align*}
which satisfies 
$$
\begin{aligned}
    \norm{\bv_{\hat{f}_{\mathcal{P}_2,true}}}_\infty
&\leq 
\max
\left\{
B_{true}\cdot \norm{\bv_{\hat{f}_{check, deep, \mathcal{P}_{2}}}^{(5M^d)}}_\infty,
\norm{\bv_{\hat{f}_{deep,\mathcal{P}_2}}}_\infty,
\norm{\bv_{\hat{f}_{check,deep,\mathcal{P}_2}}}_\infty 
\right\}\\
&\leq e^{c_{48} \cdot(p+1) \cdot M^d}
\end{aligned}$$
and $\left(\bv_{\hat{f}_{\mathcal{P}_2,true}}\right)^{(L)}_{1,0}=0, 
\quad
\norm{\left(\bv_{\hat{f}_{\mathcal{P}_2,true}}\right)^{(L)}_{1,j>0}}_\infty\leq 1$.\\
Remark that by successively applying $\hat{f}_{id}$ to the output of the networks $\hat{f}_{deep, \P_2}$ and $\hat{f}_{check, deep, \P_2}$ we can achieve that both networks have depth
\begin{align*}
L=5M^d+\left\lceil \log_4\left(M^{2p+4 \cdot d \cdot (q+1)} \cdot e^{4 \cdot (q+1) \cdot (M^d-1)}\right)\right\rceil  \cdot \lceil \log_2(\max\{q+1, 2\})\rceil.
\end{align*}
Furthermore, it is easy to see that this networks needs at most
\begin{align*}
 \max\left\{10d+4d^2+2 \cdot \binom{d+q}{d} \cdot \left(2 \cdot (4+2\lceil e^d\rceil)+5+2d\right), \right.\\
\left. \hspace*{2.5cm} 18 \cdot (q+1) \cdot \binom{d+q}{d}\right\} + 2d^2+6d+2
\end{align*}
neurons per layer.
In the definition of the final network we use the network $\hat{f}_{w_{\P_2}, deep}$ of Lemma \ref{supple14}, and the network $\hat{f}_{mult}$ $$
\hat{f}_{\text {mult }} \in \mathcal{F}\left(\left\lceil\log _4\left(M^{2 p}\right)\right\rceil, 18\right)
$$
of Lemma \ref{lemma:fmult}
for
\[
a=1 + \Bigg(\|f\|_{C^q([-a,a]^d)} \cdot e^{(M^d-1)} + (4+2 \cdot \lceil e^d \rceil) \cdot (M^d-1)\cdot e^{(M^d-2)}\Bigg)\cdot e^{2ad},
\]
which satisfies the constraint 
$$
\norm{\bv_{\hat{f}_{mult}}}_\infty \leq 4\cdot 4^{2d}
\cdot a^2
\leq
e^{c_{49}(f)}.
$$
Again we synchronize the depth of $\hat{f}_{w_{\P_2}, deep}$ and $\hat{f}_{\P_2, true}$ to insure that both networks have
\begin{align*}
\bar L = 5M^d+\left\lceil \log_4\left(M^{2p+4 \cdot d \cdot (q+1)} \cdot e^{4 \cdot (q+1) \cdot (M^d-1)}\right)\right\rceil  \cdot \lceil \log_2(\max\{q,d\}+1)\rceil
\end{align*}
many layers. The final network is given by
\begin{align*}
\hat{f}(x) = \hat{f}_{mult}\left(\hat{f}_{\P_2, true}(x), \hat{f}_{w_{\P_2}, deep}(x)\right).
\end{align*}
By Lemma \ref{supple14} and by definition of $\hat{f}_{\P_2, true}(x)$, it is easy to see that both networks have no offset in their respective output layers, thus Lemma \ref{lemma:composition-weight-bound} b) is applicable
and we get
$$
\begin{aligned}
 \norm{\bv_{\hat{f}}}_\infty 
\leq
&\max
\bigg\{
        \norm{\bv_{\hat{f}_{mult}}}_\infty,
        \max
        \left\{ 
        \norm{\bv_{\hat{f}_{\P_2, true}}}_\infty,
        \norm{\bv_{\hat{f}_{w_{\P_2}, deep}}}_\infty
        \right\}, \\
        &\max
        \left\{
        \norm{\left(\bv_{\hat{f}_{\P_2, true}}\right)^{(\bar L)}}_\infty
        \cdot 
        \norm{\left(\bv_{\hat{f}_{mult}}\right)^{(0)}}_\infty,
        \norm{\left(\bv_{\hat{f}_{w_{\P_2}, deep}}\right)^{(\bar L)}}_\infty
        \cdot 
        \norm{\left(\bv_{\hat{f}_{mult}}\right)^{(0)}}_\infty
        \right\}\bigg\}\\
        \leq &
        \max \left\{ e^{c_{49}(f) \cdot M^d },  \max
        \left\{ 
        \norm{\bv_{\hat{f}_{\P_2, true}}}_\infty,
        \norm{\bv_{\hat{f}_{w_{\P_2}, deep}}}_\infty
        \right\}\right\}\\
\leq& e^{c_{50} \cdot(p+1) \cdot M^d}
\end{aligned}
$$

This network is contained in the network class $\mathcal{F}(L,r)$ with
\begin{align*}
L=&5M^d+\left\lceil \log_4\left(M^{2p+4 \cdot d \cdot (q+1)} \cdot e^{4 \cdot (q+1) \cdot (M^d-1)}\right)\right\rceil  \cdot \lceil \log_2(\max\{q,d\}+1)\rceil\\
&+\lceil \log_4(M^{2p})\rceil
\end{align*}
and 
\begin{align*}
r= &\max\left\{10d+4d^2+2 \cdot \binom{d+q}{d} \cdot \left(2 \cdot (4+2\lceil e^d\rceil)+5+2d\right), \right.\\
&\left. \hspace*{1cm} 18 \cdot (q+1) \cdot \binom{d+q}{d}\right\} + 2d^2+6d+2+\max\{18d, 4d^2+10d\}.
\end{align*}
With the same argument as in the proof of Lemma 10 of Kohler and Langer (2021) we can show the assertion.
\hfill $\Box$

In a last step of the proof, one has to apply $\hat{f}$ to slightly shifted partitions. This follows as in section A.1.12 of Kohler and Langer (2021) and has no effect on the weight bounds that were just derived.

\end{appendix}

\end{document}